\documentclass{amsart}
\usepackage{epsfig}
\usepackage{amsmath}
\usepackage{amsthm}
\usepackage{latexsym}
\usepackage{graphicx}
\usepackage{amsfonts}
\usepackage{amssymb}
\newtheorem{theorem}{Theorem}[section]
\newtheorem{lemma}[theorem]{Lemma}
\newtheorem{corollary}[theorem]{Corollary}
\newtheorem{proposition}[theorem]{Proposition}
\newtheorem{definition}[theorem]{Definition}

\newcommand{\BH}{\mathbb{H}}
\newcommand{\BR}{\mathbb{R}}
\newcommand{\BC}{\mathbb{C}}

\newcommand{\BZ}{\mathbb{Z}}

\newcommand{\mr}[1]{\mathrm{#1}}

\newcommand{\mc}[1]{\mathcal{#1}}

\newcommand{\Vol}{\operatorname{Vol}}
\newcommand{\Area}{\operatorname{Area}}
\newcommand{\lessvol}{\operatorname{lessvol}}
\newcommand{\overlapArea}{\operatorname{overlapArea}}
\newcommand{\overlapApprox}{\operatorname{overlapApprox}}
\newcommand{\acos}{\operatorname{acos}}
\newcommand{\tr}{\operatorname{tr}}

\begin{document}

\title{Minimum volume cusped hyperbolic three-manifolds}

\author{David Gabai}\footnote{Partially supported by NSF grants
DMS-0554374 and DMS-0504110.}\address{Department of
Mathematics\\Princeton University\\Princeton, NJ}
\author{Robert Meyerhoff}\footnote{Partially supported by NSF grants
DMS-0553787 and DMS-0204311.}\address{Department of Mathematics\\Boston
College\\Chestnut Hill, MA}
\author{Peter Milley}\footnote{Partially supported by NSF grant
DMS-0554624 and by ARC Discovery grant DP0663399.}\address{Department
of Mathematics and Statistics\\University of Melbourne\\Melbourne,
Australia}

\maketitle

\section{Introduction}

In this paper, we prove:
\begin{figure}[t]
\begin{center}
\begin{tabular}{|c|c|c|c|c|c|c|}
\hline\hline
m125 & m129 & m203 & m202 & m292 & m295 & m328 \\
m329 & m359 & m366 & m367 & m391 & m412 & s596 \\
s647 & s774 & s776 & s780 & s785 & s898 & s959 \\
\hline\hline
\end{tabular}
\end{center}
\caption{The cusped manifolds which generate all
   one-cusped hyperbolic $3$-manifolds with volume $\le 2.848$. The
   manifolds are denoted here as they appear in the SnapPea census.}
\end{figure}
\begin{theorem}
Let $N$ be a one-cusped orientable hyperbolic $3$-manifold with
$\Vol(N)$ $\le 2.848$. Then $N$ can be obtained by Dehn filling all but
one of the cusps of $M$, where $M$ is one of the 21 cusped hyperbolic
$3$-manifolds listed in the table in figure 1.
\label{thrm:big_thrm}\end{theorem}

In \cite{mm} the Dehn surgery spaces of the 21
manifolds listed in figure 1 are rigorously analyzed, producing a
complete list of one-cusped manifolds with volume no greater than
$2.848$ which result from the Dehn fillings described above. We 
therefore obtain:

\begin{corollary}\label{cor:low_cusped} Let N be a 1-cusped orientable 
hyperbolic 3-manifold with $\Vol(N)$ $\le 2.848$, then N is one of 
m003, m004, m006, m007, m009, m010, m011, m015, m016, or m017. 
(Notation as in the Snappea census.)\end{corollary}

This corollary extends work of Cao and Meyerhoff who had earlier 
shown that m003 and m004 were the smallest volume cusped manifolds.
Also, the above list agrees with the SnapPea census of one-cusped
manifolds produced by Jeff Weeks (\cite{w}), whose initial members are
conjectured to be an accurate list of small-volume cupsed
manifolds.

Let $N$ be a closed hyperbolic 3-manifold with simple closed geodesic
$\gamma$ and let $N_\gamma$ denote the manifold
$N\setminus\gamma$. Agol (\cite{ago}) discovered a formula
relating $\Vol(N)$ to $\Vol(N_\gamma)$ and the tube radius of
$\gamma$. Assuming certain results of Perelman, Agol and Dunfield
(see \cite{ast}) have
further strengthened that result. A straightforward calculation (see
\cite{acs}) using this stronger result, the $\log(3)/2$ theorem of
\cite{gmt}, plus bounds on the density of hyperbolic tube packings by
Przeworksi, shows that a compact hyperbolic manifold with volume less
than that of the Weeks manifold must be obtainable by Dehn filling on
a cusped manifold with volume less than or equal to $2.848$. The paper
\cite{mm} rigorously shows that the Weeks manifold is the unique
compact hyperbolic 3-manifold of smallest volume obtained by filling
any of the 10 manifolds listed in Corollary
\ref{cor:low_cusped}.  We therefore obtain,

\begin{corollary}  The Weeks manifold is the unique closed orientable 
hyperbolic 3-manifold of smallest volume.\end{corollary}

The Weeks manifold is obtained is obtained by $(5,1)$, $(5,2)$ filling
on the Whitehead link, or by $(2,1)$ filling on the manifold m003 in
the SnapPea census.


The proof of Theorem \ref{thrm:big_thrm} is based on the Mom 
technology introduced in \cite{gmm2}.  Indeed, Figure 1 lists the 
collection of Mom-$2$ and Mom-$3$ manifolds enumerated in \cite{gmm2}, 
thus we have the following equivalent formulation:

\begin{theorem} Let $N$ be a one-cusped orientable hyperbolic $3$-manifold with
$\Vol(N)$ $\le 2.848$. Then $N$ can be obtained by Dehn filling all
but one of the cusps of $M$, where $M$ is a hyperbolic Mom-$2$ or
Mom-$3$ manifold.\end{theorem}

Recall that a Mom-$n$ manifold is a 3-manifold $M$ obtained by
starting with $T\times [0,1]$ where $T$ is the 2-torus and attaching
$n$ 1-handles and $n$ valence-3 2-handles to the $T^2\times 1$ side.
Furthermore, $\partial M$ is a union of tori.  Given $N$ as in the
theorem, the goal is to find a hyperbolic Mom-3 embedded in $N$, or in
the terminology of \cite{gmm2} show that $N$ possesses an
\emph{internal Mom-n structure} for some $n\le 3$. I.e., $M$ satisfies
the condition that the interior of $M$ has a complete hyperbolic
structure of finite volume and each component of $\partial M$ bounds
(to the outside) either a solid torus or a cusp.

In practice, we think of $T\times 0$ as the torus cutting off a
maximal cusp neighborhood $U$.  In the universal covering $U$ lifts
to a collection of horoballs $\{B_i\}$. To first approximation, when
lifted to $\BH^3$, the cores of the 1-handles of the Mom-$n$ structure
will be geodesic arcs connecting two $B_i$'s.  (Being maximal, some
$B_i$'s will be tangent to each other and these points of tangency
will also be viewed as 1-handles.)  The cores of the 2-handles, when
lifted to $\BH^3$ will correspond to totally geodesic hexagons whose
sides alternately lie on the 1-handles and boundaries of $B_i$'s.

Using the $2.848$ volume bound we will show that $N$ possesses a
\emph{geometric} Mom-$n$ structure, where $n\le 3$.  This means that we 
will find $n$ $\pi_1(N)$-orbits of geodesic arcs and $n$
$\pi_1(N)$-orbits of geodesic hexagons with boundaries on the geodesic
arcs and $B_i$'s as in the previous paragraph.  With some luck, when
thickened up, these geodesics and hexagons will descend to an internal
Mom-$n$ structure on $N$.  In reality, when projected to $N$, these
hexagons may self-intersect in undesirable ways and/or the resulting
handle structures may be unsuitable for various technical reasons.
Much of this paper is devoted to the process of promoting a geometric
Mom-$n$ structure $n\le 3$ which is \emph{torus-friendly} to a
hyperbolic internal Mom-$k$ structure, $k\le n$. (See definition
\ref{def:MOM}.)

This paper is organized as follows.  In Section 2 we give a detailed
definition of geometric Mom-$n$ structure.  In Section 3 we present
several useful geometric lemmas that will be used extensively in the
rest of the paper. Then in Sections 4 and 5 we show that if $N$
satisfies the hypotheses of Theorem \ref{thrm:big_thrm}, then $N$ must
contain a geometric Mom-$n$ structure which is torus-friendly. This
part of the proof, while theoretically simple, is computationally
complicated and was completed with the use of computer assistance; the
use of rigorous floating-point computations is discussed in Section 5.

Sections 6, 7, and 8 are concerned with the process of promoting the
geometric Mom-$n$ structure produced in Section 4 and 5 to an internal
Mom-$n$ structure of the type described in \cite{gmm2}. There
is a natural geometric object associated to a geometric Mom-$n$
structure, consisting of a thickened copy of the cusp torus $\partial
U$ (which corresponds to our $T^2\times [0,1]$) together with a
``one-handle'' for every orthogonal geodesic arc in the Mom-$n$
structure and a ``2-handle'' for every hexagon. However there are
three key conditions that this geometric object must meet before it
fits the definition of an internal Mom-$n$ structure. Each of Sections
6, 7, and 8 are devoted to one of these three conditions, and to
showing that either the condition holds or else we can replace our
geometric Mom-$n$ with a ``simpler'' structure, for some appropriate
definition of ``simpler''. Section 6 is concerned with whether or not
the geometric object associated to the Mom-$n$ is embedded in $N$,
in particular whether the various handles have undesirable
intersections or self-intersections.
Section 7 is concerned with whether the components of the complement
have the correct topology, and Section 8 is concerned with whether or
not our Mom-$n$ structure has ``simply-connected lakes'' in the
language of Matveev.  Having completed these three sections, we find
that the geometric Mom-$n$ structure produced in Sections 4 and 5 will
have evolved into an internal Mom-$k$ structure for some $k\le n$.

At that point, the proof of Theorem \ref{thrm:big_thrm} reduces to an
application of Theorems 4.1 and 5.1 of \cite{gmm2}. Together those two
theorems imply that if $N$ has an internal Mom-$n$ structure with $n\le
3$ then $N$ contains an embedded submanifold $M$ which is a hyperbolic
manifold with boundary whose interior is homeomorphic to one of the
manifolds in figure 1 and such that $N-M$ is a disjoint union of solid
tori and cusps. This is the desired result. Section 9 summarizes this
argument formally.

It should be noted that while the list in Figure 1 is precisely the
list of manifolds produced by Theorem 5.1 of \cite{gmm2}, it is
somewhat redundant for the purposes of Theorem \ref{thrm:big_thrm}
of this paper. The manifold s776 is a three-cusped manifold from which
many of the two-cusped manifolds on the list--experimentally, everything
up to and including m391--can be recovered by Dehn filling. Hence the
21 manifolds in Figure 1 could be reduced to a list of 9 manifolds
with no effect on Theorem \ref{thrm:big_thrm}. However, we use the
longer list here to be consistent with \cite{gmm2}.

Finally in Section 9 we will discuss some of the ideas used in 
\cite{mm} to rigorously analyze various Dehn fillings of the 
manifolds of Figure 1.


\section{Definition of a geometric Mom-$n$}

For the rest of this paper, $N$ will refer
to an orientable one-cusped hyperbolic 3-manifold.  Suppose $N$ is such
a manifold; then $N$ possesses a \emph{maximal cusp neighborhood}
which is a closed set whose interior is homeomorphic to $T^2\times
(0,\infty)$, with the property that each torus $T^2\times\{x\}$ has
constant sectional curvature in $N$. The term ``maximal'' here means
that this cusp neighborhood is not a proper subset of any other closed
subset of $N$ with this property.

$T$ bounds a horoball when lifted to $\tilde T$ in the universal covering 
$\BH^3$ of $N$.  In practice we think of $T$ as a maximal cusp bounding a
horoball. In $\BH^3$, the cores of the various 1-handles will be 
orthogonal geodesic arcs (or points) connecting $\pi_1(N)$-translates 
of this horoball.  The cores of the 2-handles will be totally geodesic discs
with boundary alternately on horoballs and these geodesic arcs.

We prove Theorem \ref{thrm:big_thrm} using \emph{geometric Mom-$n$
structures}, defined below.

In the universal cover $\BH^3$ of $N$, the maximal cusp neighborhood
lifts to a collection of horoballs $\{B_i\}$; any two such horoballs
have disjoint interiors but maximality implies that some pairs will be
tangent at their boundaries. Choose one such horoball and denote it
$B_\infty$. In the upper half-space model $\{(x,y,z)|z>0\}$ of
$\BH^3$, we may assume after conjugation by some element of
$Isom^+(\BH^3)$ that $B_\infty$ will be precisely the set
$\{(x,y,z)|z\ge 1\}$; then every other $B_i$ will appear as a sphere
with center $(x_i,y_i,z_i)$ and radius $z_i$ for some
$0<z_i\le \frac{1}{2}$. Let $H\subset\pi_1(N)$ be the subgroup which
fixes $B_{\infty}$, so that $B_{\infty}/H$ is homeomorphic to the cusp
neighborhood. Let $d_{E}$ denote the distance function in the subspace
(Euclidean) metric along the boundary of $B_{\infty}$.

The \emph{center} of a horoball $B_j\not =B_\infty$ is the limiting
point of the horoball on the sphere at infinity. For example, if $B_j$
appears in the upper half-space model as a sphere centered at
$(x_j,y_j,z_j)$ with radius $z_j$, then the center of the horoball is
the point $(x_j,y_j,0)$. We will sometimes refer to the complex number
$x_j+iy_j$ as the center of $B_j$ in this case, and define the center
of $B_\infty$ to be $\infty$; then the center of a horoball is always
an element of $\hat{\BC}=\BC\cup\infty$.

In addition to this we will define the \emph{orthocenter} of a
horoball $B_j\not =B_\infty$ to be the point on $\partial B_\infty$
which is closest to $B_j$. Uniqueness is guaranteed by the fact that
the boundary of any horoball has positive sectional curvature in
$\BH^3$, being a Euclidean surface in a negatively curved space. In
the upper half-space model, the orthocenter of $B_j$ is just the point
on the surface $z=1$ directly above the center of $B_j$.

Given two horoballs $A$ and $B$, neither equal to $B_{\infty}$, we
will say that $A$ and $B$ are in the same \emph{orthoclass} if either
$A$ and $B$ lie in the same $H$-orbit or there exists some
$g\in\pi_1(N)$ such that $g(A)=B_{\infty }$ and $g(B_{\infty})=B$. In
the latter case we say that $A$ and $B$ lie in
\emph{conjugate} $H$-orbits. We denote the orthoclasses by $\mathcal{O}(1)$,
$\mathcal{O}(2)$, and so forth. For any $B\in\mathcal{O}(n)$ we call
$d(B,B_{\infty})$ the \emph{orthodistance} and denote it $o(n)$; this
is clearly well-defined. Order the orthoclasses $\mathcal{O}(1)$,
$\mathcal{O}(2)$, \ldots\ in such a way that the corresponding
orthodistances are non-decreasing: $0=o(1)\le o(2) \le\cdots$. We will
refer to this as the
\emph{orthodistance spectrum}. In addition we also define $e_{n}=\exp
(o(n)/2)$, and refer to the sequence $1=e_1\le e_2\le\cdots$ as the
\emph{Euclidean spectrum}. Note that if $A\in\mathcal{O}(n)$, then in
the upper half-space model the point on $\partial A$ which is closest
to $B_\infty$ must appear to be at a height of $\exp(-d(A,B_\infty))$,
which equals ${e_n}^{-2}$ since $d(A,B_\infty)=o(n)$. The choice of
the word ``Euclidean'' actually comes from Lemma
\ref{lem:eucl_dist}.

Closely related to the orthoclasses are another set of equivalence classes
which we will call the \emph{orthopair classes}. These are just the
equivalence classes of the action of $\pi_1(N)$ on the set of unordered
pairs of horoballs $\{A,B\}$. It follows immediately from the definition that
$A$ and $B$ are in the same orthoclass if and only if $\{A,B_{\infty}\}$ and
$\{B,B_{\infty}\}$ lie in the same orthopair class. Hence we will occasionally
abuse notation and denote the orthopair classes by $\mathcal{O}(1)$,
$\mathcal{O}(2)$, \ldots\ as well. The definition of orthodistance in this
context is clear.

\begin{definition}  A \emph{$(p,q,r)$-triple} (or equivalently a 
triple of \emph{type
$(p,q,r)$}) is a triple of horoballs $\{B_1,B_2,B_3\}$ with the
property that $\{B_1,B_2\}\in\mathcal{O}(p)$, $\{B_2,B_3%
\}\in\mathcal{O}(q)$, and $\{B_3,B_1\}\in\mathcal{O}(r)$, possibly
after re-ordering.\end{definition}

Now we come to the key definition of this paper.

\begin{definition}
A \emph{geometric Mom-$n$ structure} is a collection of $n$
triples of type $(p_1,q_1,r_1)$, \ldots, $(p_n,q_n,r_n)$, no two of
which are equivalent under the action of $\pi_1(N)$, and such that the
indices $p_i$, $q_i$, and $r_i$ all come from the same $n$-element
subset of $\BZ_+$.\label{def:MOM}

We will occasionally drop the word ``structure'' when our meaning is
otherwise clear.

A geometric Mom-$n$ will be said to be \emph{torus-friendly} if
$n=2$ or if $n=3$ and the Mom-$3$ does not possess exactly two triples
of type $(p,q,r)$ for any set of distinct positive indicies $p$, $q$,
and $r$. (The geometrical implications of this term will be explained
in Section 7).
\end{definition}

So, for example, a $(1,1,3)$-triple and a $(1,3,3)$-triple would
constitute a geometric Mom-2, while a $(1,1,2)$-triple and a
$(1,1,3)$-triple would not. A $(1,1,2)$-triple, a $(1,1,3)$-triple,
and a $(1,2,3)$-triple, however, would constitute a geometric
Mom-$3$ which furthermore is torus-friendly. A $(1,1,2)$-triple and two
$(1,2,3)$-triples which are not equivalent under the action of
$\pi_1(M)$ form a geometric Mom-$3$ which is not torus-friendly.
Although this definition can clearly be generalized, in this
paper we will only be discussing geometric Mom-$n$'s where $n=2$
or $3$, and where the indices all come from the set $\{1,2,3,4\}$.

The connection between geometric Mom-$n$'s and the internal Mom-$n$
structures of \cite{gmm2} is clear.
The term
``geometric'' is meant to highlight the fact that this definition
does not include any of the topological assumptions
that are part of the definition of an internal Mom-$n$ structure, such
as embeddedness. Nevertheless, the correspondence between
geometric and internal Mom-$n$ structures is the key to proving
Theorem \ref{thrm:big_thrm}.

\section{Geometrical lemmas}


Throughout this chapter we will be using the upper half-space model of
$\BH^3$. We will take certain facts, listed below, as given; a reader
who is interested may refer to \cite{fen} for more information.

First, the orientation-proving isometries of $\BH^3$ can
be identified with the matrix group $\mr{PSL}(2,\BC)$ in a natural
way. Each element of $\mr{PSL}(2,\BC)$ acts on the sphere at infinity
$\hat{\BC}=\BC\cup\{\infty\}$ by the corresponding Mobius
transformation, i.e.
\[
\left[\begin{array}{cc}
a & b \\
c & d
\end{array}\right]:\ z\mapsto \frac{az+b}{cz+d}
\]
if $z\in\BC$, and $\infty\mapsto \frac{a}{c}$. The action of
$\mr{PSL}(2,\BC)$ on $\BH^3$ itself can be expressed similarly, using
quaternions. In the upper half-space model, the point $(x,y,t)$ in
$\BH^3$ with $t>0$ can be associated to the quaternion $x+yi+tj$; then
the action of $\mr{PSL}(2,\BC)$ can be expressed as
\[
\left[\begin{array}{cc}
a & b \\
c & d
\end{array}\right]:\ x+yi+tj\mapsto (a(x+yi+tj)+b)(c(x+yi+tj)+d)^{-1}
\]
The resulting quaternion will always be equal to $u+vi+sj$ for some
real $(u,v,s)$ with $s>0$.

If $g\in\mr{PSL}(2,\BC)$ is not the identity we can determine from the
trace of $g$ whether or not $g$ is hyperbolic, parabolic, or elliptic
as an isometry (note that trace is only defined up to sign in
$\mr{PSL}(2,\BC)$). For example, $g$ is an elliptic isometry (that is,
$g$ is a rotation about a line in $\BH^3$) if and only if
$\tr(g)=\pm
2\cos(\theta/2)$, where $\theta\in[-\pi,\pi]$ is the angle of rotation
of $g$. Note that in contrast to the hyperbolic and parabolic cases,
an elliptic isometry cannot be an element of $\pi_1(N)$, where $N$ is
a $1$-cusped hyperbolic $3$-manifold.

One final fact about $\mr{PSL}(2,\BC)$ that we will use: if
\(
g=\left[\begin{array}{cc} a & b \\ c & d \end{array}\right]
\)
then the image of $B_\infty$ under $g$ will be a horoball which
appears as a Euclidean ball of diameter $|c|^{-2}$. More
generally, the subset $B(t)=\{(x,y,z)|z\ge t\}$ of $\BH^3$ will be mapped
to a horoball which appears as a Euclidean ball of diameter
$t^{-1}|c|^{-2}$. This can be demonstrated as follows. Since
$g(\infty)=\frac{a}{c}$, the image of $B(t)$ will be a horoball with
center $\frac{a}{c}\in\BC$. Therefore the orthocenter of $g(B(t))$
will correspond to the quaternion $\frac{a}{c}+\delta j$ where
$\delta$ is the diameter we seek. This is just the point where
$g(B(t))$ intersects the line in $\BH^3$ from $\infty$ to
$\frac{a}{c}$. Therefore the pre-image of this point under $g$ is the
point where $B(t)$ intersects the line from
$g^{-1}(\infty)=-\frac{d}{c}$ to $g^{-1}(\frac{a}{c})=\infty$. In
other words,
\[
\frac{a}{c}+\delta j = g\left(-\frac{d}{c}+tj\right)
\]
Direct calculation with quaternions then yields the desired result.

We now begin enumerating the geometrical lemmas that we will use in
the rest of the paper. The following lemma first appears (using
different language) in \cite{ada}.

\begin{lemma}
Every orthoclass consists of two $H$-orbits. \label{lem:two_orbits}
\end{lemma}

\noindent\emph{Proof:} It is fairly clear from the definition that each
orthoclass contains no more than two $H$-orbits. If an orthoclass contains
exactly one $H$-orbit, then we must have $g^{-1}(B_{\infty})=hg(B_{\infty})$
for some $h\in H$, $g\not \in H$. In other words, $ghg\in H$, and therefore
$(gh)^2\in H$. But as an isometry $gh$ must be either hyperbolic, parabolic,
or the identity. In the first case $(gh)^2$ would also be hyperbolic (and
hence not in $H$), and in the last two we would have $gh\in H$ and hence $g\in
H$. Either case is a contradiction. \qed

\bigskip
The following lemma appears in \cite{cm}; a related lemma for compact
manifolds appears in \cite{gmm}. It is reproduced here for the sake of
completeness.

\begin{lemma}
If $A$ and $B$ are both in $\mathcal{O}(n)$, and if $g\in\pi_1(N)$ is an
isometry such that $g(B)=B_{\infty}$, then $g(A)\not \in\mathcal{O}(n)$. \label{lem:no_bad_triples}
\end{lemma}

\noindent\emph{Proof:} There are two cases to consider: either $A$ and
$B$ lie in different $H$-orbits or else they lie in the same
$H$-orbit.

In the first case, by the definition of an orthoclass $g(B_{\infty})=h(A)$ for
some $h\in H$. By replacing $g$ with $h^{-1}g$, therefore, we may assume that
$g(B_{\infty})=A$. Suppose that $g(A)$ is in $\mathcal{O}(n)$. Then either
$g(A)=h_1(B)$ for some $h_1\in H$, in which case we have $g^2(B_{\infty
})=h_1g^{-1}(B_{\infty})$, or else $g(A)=h_1(A)$, in which case we have
$g^2(B_{\infty})=h_1g(B_{\infty})$.

Either way we have $g^2=h_1g^{\pm1}h_2$ for some $h_2\in H$. Let
\[
g = \left[
\begin{array}
[c]{cc}%
a & b\\
c & d
\end{array}
\right]  ,\ h_{i} = \left[
\begin{array}
[c]{cc}%
1 & k_{i}\\
0 & 1
\end{array}
\right]  ,\ i \in\{1,2\}
\]
Expanding both sides of the equation $g^2=h_1g^{\pm1}h_2$ and taking the
$(2,1)$-entry of the resulting matrix on each side, we get $c(a+d) = \pm c$.
Therefore the square of the trace of $g$ is 1, which implies that $g\in\pi
_1(N)$ is elliptic of order 3, a contradiction.

The other case is when $A$ and $B$ lie in the same $H$-orbit, i.e. $A=h(B)$
for some $h\in H$. Suppose that $g(A)\in\mathcal{O}(n)$ and that $g(A)$ lies
in the same $H$-orbit as $A$ and $B$, i.e. $g(A)=k(B)$ for some $k\in H$. Then
$g(A)$ and $g(B_{\infty})$ both lie in $\mathcal{O}(n)$ but must lie in
different $H$-orbits (see Lemma \ref{lem:two_orbits}). Furthermore
$gk^{-1}(g(A))=g(B)=B_{\infty}$, and $gk^{-1}(B_{\infty})=g(B_{\infty})$.
Therefore replacing $B$ with $g(A)$, $A$ with $g(B_{\infty})$, and $g$ with
$gk^{-1}$ reduces the problem to the previous case. On the other hand, suppose
$g(A)\in\mathcal{O}(n)$ but $g(A)$ does not lie in the same $H$-orbit as $A$
and $B$. Then $g(A)=kg(B_{\infty})$ for some $k\in H$. Since $A=h(B)=hg^{-1}%
(B_{\infty})$, we have $ghg^{-1}=kgl$ for some $l\in H$. This last equation,
after some manipulation (note that $h$ and $l$ commute), implies that
$k^{-1}(gl)h=(gl)^2$, which leads to a contradiction just as in the
first case. This completes the proof of the lemma. \qed\ 

\bigskip
Another way of phrasing the above result is to say that the unordered pair
$\{A,B\}$ does not lie in the same orthopair class as the pairs $\{A,B_{\infty
}\}$ and $\{B,B_{\infty}\}$. Hence we immediately get the following:

\begin{corollary}
There are no $(n,n,n)$-triples. \label{cor:no_bad_triples}
\end{corollary}\qed

\bigskip
The following lemma will be quite useful when studying horoball
diagrams on the surface of $B_{\infty}$:

\begin{lemma}
Let $A$, $B$ be two horoballs not equal to $B_{\infty}$ and let $p$, $q$ be
their orthocenters. Suppose $A\in\mathcal{O}(m)$ and $B\in\mathcal{O}(n)$, and
suppose that $d(A,B)=O(r)$. Then the distance between $p$ and $q$ along the
surface of $B_{\infty}$ is given by $d_{E}(p,q)=e_{r}/(e_{m}e_{n})$. \label{lem:eucl_dist}
\end{lemma}

\noindent\emph{Proof:} Using the upper half-space model, we may assume that
$p=(0,0,1)$ and $q=(x,0,1)$ where $x=d_{E}(p,q)$. Then $A$ is a ball of height
$e_{m}^{-2}$ tangent to $S_{\infty}^2$ at 0 and $B$ is a ball of height
$e_{n}^{-2}$ tangent to $S_{\infty}^2$ at $(x,0)$. Consider the following
isometry given as an element of $\mathrm{PSL}(2,\BC)$:
\[
\sigma=\left[
\begin{array}
[c]{cc}%
1 & 0\\
-x^{-1} & 1
\end{array}
\right]
\]
Since $\sigma$ is a parabolic isometry which fixes 0, $\sigma$ preserves $A$.
And $\sigma$ sends $x$ to $\infty$; therefore $\sigma$ sends $B$ to a horoball
centered at infinity. The height of $\sigma(B)$ will be the same as the height
of the image of the point $(x,0,{e_{n}}^{-2})$, which is readily determined to
be $x^2{e_{n}}^2$. Therefore the distance between $A=\sigma(A)$ and
$\sigma(B)$ is the logarithm of the ratio of $x^2{e_{n}}^2$ to the
diameter of $A$, which implies that
\[
{e_{r}}^2=\frac{x^2{e_{n}}^2}{{e_{m}}^{-2}}%
\]
which proves the lemma. \qed\ 

\bigskip
The next two lemmas concern pairs of lines joining the centers of
horoballs in $\BH^{3}$; they will be used extensively in
discussing when geometric Mom-$n$'s are embedded in Section
6.

\begin{lemma}
Suppose $A$, $B$, $C$, and $D$ are disjoint horoballs in $\BH^{3}$ with
no horoball contained in the interior of another (or equivalently, with no two
having the same center on the sphere at infinity). Let $\lambda_1$ be the
line joining the centers of $A$ and $C$, $\lambda_2$ be the line joining the
centers of $B$ and $D$, and let $x$ be the distance between $\lambda_1$ and
$\lambda_2$. Then
\[
e^{\frac{d(A,B)+d(C,D)}{2}}+e^{\frac{d(A,D)+d(B,C)}{2}}=e^{\frac
{d(A,C)+d(B,D)}{2}}\cosh x
\]
\end{lemma}

\begin{corollary}
If $\{A,B\}\in\mathcal{O}(h)$, $\{B,C\}\in\mathcal{O}\left(  j\right)  $,
$\{C,D\}\in\mathcal{O}(k)$, $\{D,A\}\in\mathcal{O}\left(  l\right)  $,
$\{A,C\}\in\mathcal{O}(m)$, and $\{B,D\}\in\mathcal{O}(n)$, then%
\[
e_{h}e_{k}+e_{j}e_{l}=e_{m}e_{n}\cosh x
\]
\label{cor:two_line_dist}
\end{corollary}

\noindent\emph{Proof:} Arrange the horoballs in the upper half-space model in
such a way that $A$ is the horoball at infinity with height $1$ and $C$ is
centered at $0$. Let $z$, $w\in\BC$ be the centers of $B$ and $D$
respectively. To compute the distance between the two lines $\lambda_1$ and
$\lambda_2$, we turn to \cite{fen} which says that
\[
\cosh(x+iy)=-\frac{1}{2}\tr\left(  \Lambda_1\Lambda_2\right)
\]
where
\[
\Lambda_1=\left[
\begin{array}
[c]{cc}%
i & 0\\
0 & -i
\end{array}
\right]  \in\mathrm{PSL}(2,\BC)
\]
is an elliptic element of order $2$ fixing the line from $0$ to $\infty$ while%
\[
\Lambda_2=\frac{i}{w-z}\left[
\begin{array}
[c]{cc}%
w+z & -2wz\\
2 & -w-z
\end{array}
\right]  \in\mathrm{PSL}(2,\BC)
\]
is an elliptic element of order $2$ fixing the line from $z$ to $w$, and where
$y$ is the relative angle between the two lines along the shortest arc between
them, which is only defined modulo $\pi$. Then by direct calculation,%
\begin{align*}
\cosh(x+iy)  &  =-\frac{i}{2(w-z)}\tr\left(  \left[
\begin{array}
[c]{cc}%
i & 0\\
0 & -i
\end{array}
\right]  \left[
\begin{array}
[c]{cc}%
w+z & -2wz\\
2 & -w-z
\end{array}
\right]  \right) \\
&  =-\frac{i}{2(w-z)}\left(  i(w+z)-i(-w-z)\right) \\
&  =\frac{w+z}{w-z}%
\end{align*}
To determine $\cosh x$ note that as $y$ varies, $\cosh(x+iy)=\cosh x\cos
y+i\sinh x\sin y$ varies along an ellipse in the complex plane. It is
straightforward to determine that the foci of this ellipse are at $-1$ and
$1$, and that therefore for any $y$,%
\[
\left|  \cosh(x+iy)-1\right|  +\left|  \cosh(x+iy)+1\right|  =2\cosh x
\]
Therefore,%
\begin{align*}
\cosh x  &  =\frac{1}{2}\left(  \left|  \frac{w+z}{w-z}-1\right|  +\left|
\frac{w+z}{w-z}+1\right|  \right) \\
&  =\frac{|w|+|z|}{|w-z|}%
\end{align*}
So $|w|+|z|=|w-z|\cosh x$. But by Lemma \ref{lem:eucl_dist}, $|w|^2%
=e^{d(C,D)-d(A,C)-d(A,D)}$, $|z|^2=e^{d(C,B)-d(A,C)-d(A,B)}$, and
$|w-z|^2=e^{d(B,D)-d(A,B)-d(A,D)}$. The result follows.\qed\ 

\begin{lemma}
Let $N$ be a cusped hyperbolic manifold with cusp neighborhood $T$,
and suppose $A$, $B$, $C$, and $D$ are all horoballs in $\BH^{3}$
which are lifts of $T$. Define $\lambda_1$, $\lambda_2$, and $x$ as in
the previous lemma. If there exists an element $g\in\pi_1(N)$ such
that $g(A)=B$ and $g(C)=D$, and if $d(A,C)=d(B,D)$ is less than or
equal to $2\log(1.5152)$ then $x\ge 0.15$. \label{lem:dist_line_self}
\end{lemma}

\noindent\emph{Proof:} Suppose $x<0.15$; we will establish a contradiction.
For sake of notation, suppose that $\{A,C\}$ and $\{B,D\}$ are both elements
of $\mathcal{O}(k)$. Note that this implies that $e_{k}\le 1.5152$. Arrange the
four horoballs in the upper half-space model so that $A$ is centered at
infinity with height $e_{k}$ and $C$ is centered at $0$ with height
$e_{k}^{-1}$. Then suppose that%
\[
g=\left[
\begin{array}
[c]{cc}%
a & b\\
c & d
\end{array}
\right]  \in\mathrm{PSL}(2,\BC)
\]
We wish to compute $x$ in terms of $a$, $b$, $c$, and $d$ in the same fashion
as in the previous lemma. Since $g$ sends $0$ to $\frac{b}{d}$ and $\infty$ to
$\frac{a}{c}$, by the same arguments as before we have that for some $y$,%
\begin{align*}
\cosh(x+iy)  &  =\frac{\frac{a}{c}+\frac{b}{d}}{\frac{a}{c}-\frac{b}{d}}\\
&  =\frac{ad+bc}{ad-bc}\\
&  =ad+bc\\
&  =2ad-1
\end{align*}
where the last two steps used $ad-bc=1$. Hence%
\[
ad=\frac{1}{2}\left(  \cosh(x+iy)+1\right)
\]
Since $x<0.15$, this implies that $ad$ lies strictly inside an ellipse in the
complex plane with foci at $0$ and $1$, whose boundary intersects the real
axis at the points $\frac{1}{2}(\cosh0.15+1)$ and
$\frac{1}{2}(1-\cosh0.15)$.

Now recall that if $c\not =0$ then the element $g$ sends a horoball of
height $t$ centered at infinity to a horoball of height
$t^{-1}|c|^{-2}$ (whereas if $c=0$ then $g$ fixes the point at
infinity). So since $g(A)=B$ and since $B$ and $A$ have disjoint
interiors, we must have $e_{k}^{-1}|c|^{-2}\leq e_{k}$, or in other
words $|c|\geq e_{k}^{-1}\ge 1.5152^{-1}$.  But $C$ and $g(A)=B$ also
have disjoint interiors, so let%
\[
h=\left[
\begin{array}
[c]{cc}%
0 & i\\
i & 0
\end{array}
\right]  \in\mathrm{PSL}(2,\BC)
\]
Then $h$ swaps $A$ and $C$; therefore $hg(A)$ and $h(C)=A$ have disjoint
interiors. Since%
\[
hg=\left[
\begin{array}
[c]{cc}%
ic & id\\
ia & ib
\end{array}
\right]
\]
this implies that $e_{k}^{-1}|ia|^{-2}\leq e_{k}$, or in other words $|a|\geq
e_{k}^{-1}\ge 1.5152^{-1}$. Similarly, $A$ and $D=g(C)=gh(A)$ have disjoint
interiors,
which implies that $|d|\ge 1.5152^{-1}$, and $C$ and $D$ have disjoint
interiors, so
$h(C)=A$ and $h(D)=hgh(A)$ have disjoint interiors, which implies that
$|b|\ge 1.5152^{-1}$. Combining these facts we have $|ad|\ge
1.5152^{-2}$ and $|ad-1|=|bc|\ge 1.5152^{-2}$. Hence we can conclude
that $ad$ lies somewhere in the shaded region
indicated in
\begin{figure}[tb]
\begin{center}
\includegraphics{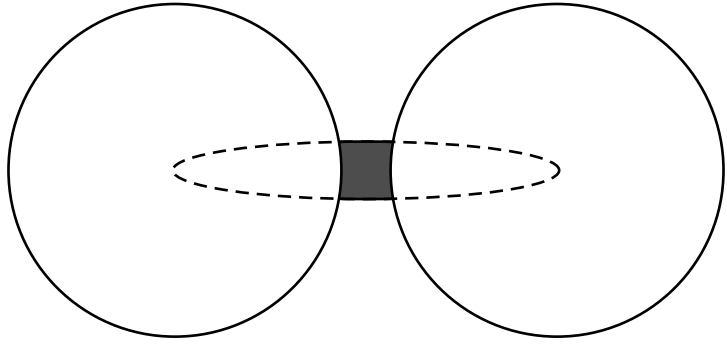}
\end{center}
\caption{The possible range of $ad$ in the complex plane. The shaded
region is defined by the equation $|z-0|+|z-1|\le\cosh 0.15$
(i.e. the solid ellipse bounded by the curve
$\frac{1}{2}(\cosh(0.15+it)+1)$ where $t$ is real) and the equations
$|z|\ge 1.5152^{-2}$ and $|z-1|\ge 1.5152^{-2}$.}
\end{figure}
Figure 2.

It's worthwhile at this point to sketch out the motivation for the
argument that follows. Roughly speaking, we've shown that $ad$ is
approximately equal to $1/2$. Given the stated lower bounds on $|a|$
and $|d|$, this implies that $|a|$ and $|d|$ are each approximately
equal to $\sqrt{1/2}$. If it were the case that $ad=1/2$ and
$|a|=|d|=\sqrt{1/2}$, then $a+d$ would have to be a real number
between $-\sqrt{2}$ and $\sqrt{2}$ which would imply that
$g\in\pi_1(N)$ is elliptic, a contradiction. Since we only have
approximate equality in the previous statement, we wish to show that
$g$ is ``approximately elliptic''. More rigorously, we wish to show
that $a+d$ lies close enough to the real interval from $-\sqrt{2}$ to
$\sqrt{2}$ to ensure that $g^n(A)$ or $g^n(C)$ intersects $A$ for some $n$
(specifically $n=2$, $3$, or $4$), a contradiction which will complete
the proof of the lemma.

To make this argument work we need to break the problem into two cases
depending on the value of $ad$. The first case will be when
$\Re(ad)\ge 1/2$ and $ad$ lies in the region in the above figure; the
second case is when $\Re(ad)\le 1/2$. We proceed with the proof in the
first case.

So suppose $\Re(ad)\ge 1/2$ and $ad$ lies in the region in the above
figure;
we wish to determine a domain for the trace $\tau=a+d$ of $g$. Suppose
that $a=|a|e^{i\theta}$ and $d=|d|e^{i\phi}$; then%
\[
\tau=e^{i\left(  \frac{\theta+\phi}{2}\right)  }\left(  \left(
|a|+|d|\right)  \cos  \frac{\theta-\phi}{2}  +i\left(
|a|-|d|\right)  \sin  \frac{\theta-\phi}{2}  \right)
\]
Hence $\tau$ lies in an ellipse, centered at the origin,
whose major axis has length $2(|a|+|d|)$ and whose minor axis has length
$2||a|-|d||$, and whose major axis is tilted away from the real axis by
half the argument of $ad$. 

Note that we have upper bounds on all three of these
quantities. First, the argument of $ad$ is bounded since $ad$ is
contained in a bounded simply-connected region which does not contain
$0$. Second, since $|ad|$ is bounded above and $|a|$ and $|d|$ are
both bounded below, the point $(|a|,|d|)$ is contained in a region in
the first quadrant of the real plane which is bounded by the lines
$x=1.5152^{-1}$ and $y=1.5152^{-1}$ and the hyperbola $xy=D$ where $D$
is the maximum possible value of $|ad|$.  This implies that $|a|+|d|$
and $||a|-|d||$ are both bounded as well.  Specifically, we can state
the following:
\begin{align*}
\frac{\theta+\phi}{2} & \le 0.07473 \\
|a|+|d| & \le 1.5323 \\
\left||a|-|d|\right| & \le 0.2124
\end{align*}
These facts together imply that $\tau$ must be contained in a bounded
region near the origin, as sketched out in
\begin{figure}[tb]
\begin{center}
\includegraphics{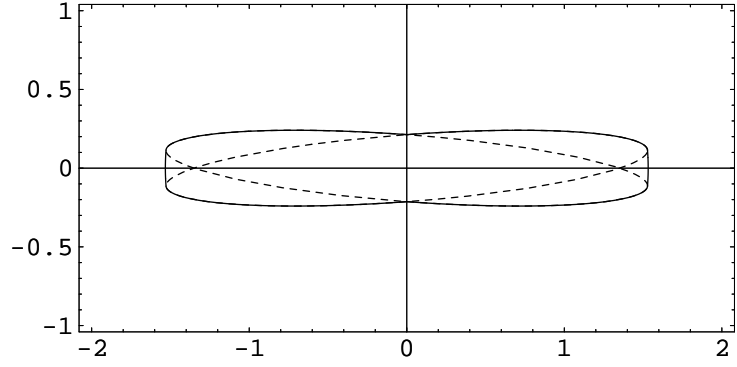}
\end{center}
\caption{The region of possible values for $a+d$ (solid contour),
generated by sweeping out ellipses with constant major axis and minor
 axis length, parameterized by the angle of inclination with the
 $x$-axis. The dashed contours represent the beginning and ending
 ellipses.}
\end{figure} figure 3.

Now consider $g^2$, $g^3$, and $g^4$; these elements do not fix
$\infty$ since $g$ does not. Similarly they do not fix $0$.
By direct calculation,%
\begin{align*}
g^2  &  =\left[
\begin{array}
[c]{cc}%
\ast & b\tau\\
c\tau & \ast
\end{array}
\right] \\
g^{3}  &  =\left[
\begin{array}
[c]{cc}%
\ast & b(\tau^2-1)\\
c(\tau^2-1) & \ast
\end{array}
\right] \\
g^{4}  &  =\left[
\begin{array}
[c]{cc}%
\ast & b(\tau^3-2\tau)\\
c(\tau^{3}-2\tau) & \ast
\end{array}
\right]
\end{align*}
(Here $\ast$ is used to denote entries whose value is unimportant.)
Since $A$ and $g^n(A)$ (respectively $C$ and $g^n(C)$) have disjoint interiors,
by the same arguments as before all three of the quantities $|c\tau|$,
$|c(\tau^2-1)|$, and $|c(\tau^3-2\tau)|$ (respectively $|b\tau|$,
$|b(\tau^2-1)|$, and $|b(\tau^3-2\tau)|$) must be no less than
$1.5152^{-1}$.
Since $|bc|=|ad-1|$ and since $ad$ lies in a bounded simply connected
region which does not contain $1$, $\sqrt{|bc|}$ is bounded above and
hence $|\tau|$, $|\tau^2-1|$, and $|\tau^3-2\tau|$ are bounded below.
Specifically we have
\begin{align*}
\left|  \tau\right|   &  \geq 0.9281\\
\left|  \tau^2-1\right|   &  \geq 0.9281\\
\left|  \tau^{3}-2\tau\right|   &  \geq 0.9281
\end{align*}
But as demonstrated in
\begin{figure}[tb]
\begin{center}
\includegraphics{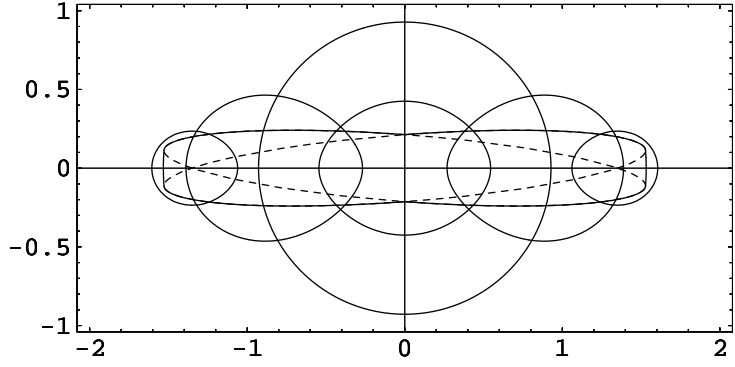}
\end{center}
\caption{The region from figure 3, contained inside the union of the
regions $|\tau|<0.9281$, $|\tau^2-1|<0.9281$, and
$|\tau^3-2\tau|<0.9281$.}
\end{figure}
figure 4, the intersection of the three
regions described by the above inequalities is disjoint from the
previously determined range of possible values of $\tau$, a
contradiction. This completes the proof in the first case.

The proof in the second case is nearly identical, only the specific
numbers differ: replace $0.07473$ with $0.08599$, $1.5323$ with
$1.4262$, $0.2124$ with $0.1062$, and $0.9281$ with $0.8698$.
This completes the proof of the lemma.\qed\

\bigskip
The remaining lemmas in this chapter will be quite arbitrary and
geometrical, but they will be quite useful in the next chapter when we
begin to compute the volume of manifolds which do not possess
geometric Mom-$n$'s. We start with a lemma and definition which first
appeared in \cite{cm}:

\begin{lemma}
Let $B$ be a horoball centered at infinity in the upper half-space
model of $\BH^3$, whose boundary has Euclidean height $1/b$. Let $A$
be a half-space, i.e. one of the two connected components of the
complement of a plane; assume that $A$ does not contain the point at
infinity, and that the plane which defines $A$ appears in the
upper half-space model as a Euclidean hemisphere with radius
$1/a$. Assume that $a<b$, so that $A\cap B$ is non-empty. Then the
volume of $A\cap B$ is
\[
\pi\left( \frac{b^2}{2a^2} - \frac{1}{2} + \log\frac{a}{b} \right)
\]
\label{lem:lessvol}\end{lemma}
\begin{definition}
Denote the above quantity by $\lessvol(a,b)$.\label{def:lessvol}
\end{definition}
\noindent\emph{Proof of Lemma \ref{lem:lessvol}:} Without loss of
generality, assume that the plane defining $A$ appears as a Euclidean
hemisphere centered at the origin. Then convert to cylindrical
coordinates; the desired volume is given by the following integral:
\[
\int_{\frac{1}{b}}^{\frac{1}{a}} \int_0^{2\pi}
\int_0^{\sqrt{\frac{1}{a^2}-z^2}} \frac{r}{z^3} \,dr\,d\theta\,dz
\]
The result follows immediately.\qed

\bigskip
The following lemma is a straightforward exercise in
trigonometry; the proof is left to the reader:
\begin{lemma}
Let $A$ and $B$ be two circular disks in the plane with radius $a$ and
$b$ respectively, such that the centers of $A$ and $B$ are $c$
units apart. Suppose that $|a-b|\le c\le a+b$, i.e. suppose that the
two disks overlap but neither disk is entirely contained in the
interior of the other. Then
\[
\Area(A\cap B) = a^2 f \left(\frac{a^2-b^2+c^2}{2ac}\right)
               + b^2 f \left(\frac{-a^2+b^2+c^2}{2bc}\right)
\]
where $f(x)=\cos^{-1}(x)-x\sqrt{1-x^2}$, if $a\not =b$. If $a=b$, then
\[
\Area(A\cap B) = 2 a^2 f\left(\frac{c}{2a}\right)
\]
\label{lem:overlapArea}\end{lemma}
\begin{definition}
Define the function $\overlapArea(a,b,c)$ to equal the right-hand side
of the first equation above if $a\not =b$, and the right-hand side of
the second equation above if $a=b$.
\label{def:overlapArea}\end{definition}

Finally, we provide a companion to the above lemma. The function
$f(x)$ described above, while relatively simple, turns out to be
unsuitable for rigorous floating-pointing computation. As will be
explained in more detail in Section 5, a polynomial approximation to
$f(x)$ will be extremely useful. Hence we provide the following:

\begin{lemma}
Let $A$ and $B$ be two circular disks in the plane with radius $a$ and
$b$ respectively, such that the centers of $A$ and $B$ are $c$ units
apart. Suppose that $|a-b|\le c\le a+b$, i.e. suppose that the two
disks overlap but neither disk is entirely contained in the interior
of the other. Then
\[
\Area(A\cap B) \le a^2 g \left(\frac{a^2-b^2+c^2}{2ac}\right)
               + b^2 g \left(\frac{-a^2+b^2+c^2}{2bc}\right)
\]
where
$g(x)=(\frac{5}{3}-\frac{\pi}{2})x^5+\frac{1}{3}x^3-2x+\frac{\pi}{2}$,
if $a\not =b$. If $a=b$, then
\[
\Area(A\cap B) \le 2 a^2 g\left(\frac{c}{2a}\right)
\]
\label{lem:overlapApprox}\end{lemma}
\begin{definition}
Define the function $\overlapApprox(a,b,c)$ to equal the right-hand
side of the first inequality above if $a\not =b$, and the right-hand
side of the second inequality above if $a=b$.
\label{def:overlapApprox}\end{definition}

\noindent\emph{Proof:} Comparing the above lemma to Lemma
\ref{lem:overlapArea}, clearly it would suffice to show that $g(x)\ge
f(x)$ for all $x\in[-1,1]$ where
$f(x)=\cos^{-1}(x)-x\sqrt{1-x^2}$. Unfortunately, this cannot be the case
as both $f(x)-\frac{\pi}{2}$ and $g(x)-\frac{\pi}{2}$ are odd
functions. However, it is true that $g(x)\ge f(x)$ for all $x\in
[0,1]$; we will use this fact in what follows. Let $h(x)=g(x)-f(x)$.
We note the following facts about $h(x)$, which are readily verified:
$h(x)$ is an odd function, $h(\pm 1)=h(0)=0$, $h(x)>0$ if $0<x<1$, and
$h(x)$ has a single local maximum at $x_{max}=0.928\ldots\ $.
By Lemma \ref{lem:overlapArea} it suffices to show that
\[
a^2 h\left(\frac{a^2-b^2+c^2}{2ac}\right)
+ b^2 h\left(\frac{-a^2+b^2+c^2}{2bc}\right) \ge 0
\]
if $a\not =b$.

Let $P$ be the center of the first circle, $Q$ the center of the
second, and let $R$ be one of the points where the circles
intersect. Let $\theta$ be the angle of the triangle $PQR$ at $P$, and
let $\phi$ be the angle at $Q$. Applying the law of cosines for $PQR$
to the above inequality, we get
\[
a^2 h\left(\cos\theta\right) + b^2 h\left(\cos\phi\right) \ge 0
\]
Note that if $\theta=0$ or $\pi$ then $\phi=0$ or $\pi$ and
vice-versa. Since $h(1)=h(-1)=0$ the lemma is true in either case. So
suppose that $\theta$, $\phi\in(0,\pi)$. Note that since $\theta$ and
$\phi$ are angles of a (possibly flat) triangle this also implies that
$\theta+\phi\in(0,\pi)$. Then we can apply the law of sines to the
above inequality to get
\[
\left(\frac{c\sin\phi}{\sin(\theta+\phi)}\right)^2
h\left(\cos\theta\right)
+ \left(\frac{c\sin\theta}{\sin(\theta+\phi)}\right)^2
h\left(\cos\phi\right) \ge 0
\]
Letting $x=\cos\theta$ and $y=\cos\phi$, it suffices to show that
\[
(1-y^2)h(x)+(1-x^2)h(y) \ge 0
\]
where $x<1$, $y<1$, and $x>-y$ (since $\theta<\pi-\phi$). Let $H(x,y)$
equal the left-hand side of the above inequality. Note that
$H(x,1)=H(1,y)=0$ since $h(1)=0$, and $H(x,-x)=0$ for $x\in [-1,1]$
since $h(x)$ is odd. Hence the lemma will be proved if we can show
that $H(x,y)$ has a non-negative value at any local minimum in the
interior of the triangle bounded by the lines $x=1$, $y=1$, and
$x+y=0$. Since $h(x)\ge 0$ if $0\le x\le 1$, $H(x,y)\ge 0$ whenever
$0\le x\le 1$ and $0\le y\le 1$. Also, $H(x,y)=H(y,x)$. Therefore it
suffices to examine local minima in the interior of the triangle
bounded by the lines $x=1$, $y=0$, and $x+y=0$. We now show that there
are no such local minima in that region.

Computing the gradient, we see that if $(x_0,y_0)$ is a local extreme
point of $H(x,y)$ then
\begin{eqnarray*}
(1-{y_0}^2) h'(x_0) &=& 2x_0 h(y_0) \\
(1-{x_0}^2) h'(y_0) &=& 2y_0 h(x_0)
\end{eqnarray*}

If $(x_0,y_0)$ is in the interior of the triangle described above then
$0<x_0<1$ and $-1<y_0<0$. Then from the first equation, since
$h(y_0)<0$ we must have $h'(x_0)<0$, and from the second equation
since $h(x_0)>0$ we must have $h'(y_0)<0$. Therefore $x_{max}<x_0<1$,
and $-1<y_0<-x_{max}$. Now suppose further that $(x_0,y_0)$ is a local
minimum. Then we must have $\partial^2 H/\partial x^2>0$ at that
point, which implies that
\[
h''(x_0)>\frac{2 h(y_0)}{1-{y_0}^2}
\]
But it can be readily computed that the maximum value of $h''(x)$ on
the interval $[x_{max},1]$ is less than $-1.6$, while the minimum value
of $2h(y_0)/(1-{y_0}^2)$ on the interval $[-1,-x_{max}]$ is greater
than $-0.6$, contradicting the above inequality. This completes the
proof of the lemma.\qed\

\section{Bounding the parameter space}

In Section 5 we will show that a one-cusped hyperbolic $3$-manifold
with volume less than or equal to $2.848$ has a geometric
Mom-$2$ or Mom-$3$ structure. This proof involves constructing volume
estimates in terms of the first three non-trivial elements of the
Euclidean spectrum $e_2$, $e_3$, and $e_4$, and then evaluating those
formulas with rigorous computer assistance. This requires us to
restrict our attention from the space of all possible values of
$(e_2,e_3,e_4)$ to a compact subset of that space. Doing so is the
purpose of this section.

Specifically, we wish to prove the following:
\begin{proposition}
Let $N$ be a one-cusped hyperbolic $3$-manifold with
$\Vol(N)\le 2.848$. Then $e_2\le 1.4751$; furthermore, one of the
following statements must be true:
\begin{itemize}
\item $e_3 \le 1.5152$, or
\item $N$ contains a geometric Mom-$2$ structure consisting only of
  $(1,1,2)$-triples and $(1,2,2)$-triples.
\end{itemize}
\label{prop:param_bounds}\end{proposition}
Note that this effectively provides upper bounds on $e_2$ and $e_3$
for manifolds with $\Vol(N)\le 2.848$; lower bounds are trivially
provided by $e_k\ge e_1=1$. This proposition does not provide an upper
bound for $e_4$, but we can do without as described in the next
section.

The proof of Proposition \ref{prop:param_bounds} depends on several
complicated estimates of area and volume and spans multiple cases.
Similar arguments, presented in less detail for brevity, will be used
for the proofs in Section 5, and hence this section should perhaps be
considered a ``warm-up'' for those results.

%
%

We begin with the following simple result:
\begin{lemma}
The volume of a one-cusped hyperbolic three-manifold $N$ is at least
\[
\frac{{e_2}^{4}\sqrt{3}}{2}-\pi\left(  {e_2}^2-1-2\log{e_2}\right)
\]
\label{lem:first_vol}
\end{lemma}

\begin{corollary}
If $e_2>1.4751$ then $\Vol(N)>2.848$. \label{cor:first_vol}
\end{corollary}

\noindent\emph{Proof:} Let $\{p_{i}\}\subset\partial B_{\infty}$ be the set of
orthocenters of horoballs belonging to $\mathcal{O}(1)$. According to lemmas
\ref{lem:no_bad_triples} and \ref{lem:eucl_dist}, $d_{E}(p_{i},p_{j})\ge 
e_2$ for all $i\not =j$. Therefore $\partial B_{\infty}$ can be packed by
circles of radius $e_2/2$ centered at each point $p_{i}$. There are two
$H$-orbits of such points by lemma \ref{lem:two_orbits}, and therefore the
area of $\partial B_{\infty}/H$ is at least $\pi({e_2}^2)/2$. Since our
packing is by circles of uniform radius, we can immediately improve this bound
by a factor of $\sqrt{12}/\pi$ (the density of the hexagonal circle packing)
to get
\[
\operatorname{Area}(\partial B_{\infty}/H)\ge {e_2}^2\sqrt{3}%
\]

However we want to estimate volume, not area. The volume of $B_{\infty}/H$ is
exactly $\operatorname{Area}(\partial B_{\infty}/H)/2$, and this is a lower
bound on the volume of $N$, but this lower bound only accounts for the volume
inside the cusp neighborhood. We would like our estimate to count some of the
volume outside the cusp neighborhood as well.

We do this by ``inflating'' the cusp neighborhood in a uniform fashion. In the
upper half-space model, this corresponds to replacing $B_{\infty}$ with a new
horoball $C$ which is centered at infinity but whose Euclidean height is some
positive number which may be less than one. In particular we choose $C$ to
have height $1/e_2$. Then
\begin{align*}
\Vol(C/H)  &  = \frac{{e_2}^2}{2}\operatorname{Area}%
(\partial B_{\infty}/H)\\
&  = \frac{{e_2}^{4}}{2}\sqrt{3}%
\end{align*}
but of course $C/H$ is no longer embedded in $N$; it is only immersed. To get
a valid lower bound for $\Vol(N)$ we must ``cut off'' $C/H$ by
subtracting the volumes of those regions where $C/H$ extends beyond the Ford
domain corresponding to $B_{\infty}/H$. The maximum height of a Ford face
corresponding to $A\in\mathcal{O}(2)$ is just $1/e_2$; therefore $C/H$ only
extends past the Ford faces corresponding to $\mathcal{O}(1)$-horoballs and
not $\mathcal{O}(2)$-horoballs (or $\mathcal{O}(n)$-horoballs for $n>2$).
Hence at most we must subtract twice the volume of the region where the
horoball $C$ intersects the half-space corresponding to some $B\in
\mathcal{O}(1)$. That volume is given by the lessvol function (see
Definition \ref{def:lessvol}); specifically,
\begin{align*}
\Vol(N)  &  \ge\Vol(C/H) - 2\lessvol(1,e_2)\\
&  = \frac{{e_2}^4\sqrt{3}}{2} - \pi\left(  {e_2}^2-1-2\log{e_2}\right)
\end{align*}
This is the desired result.

The corollary follows since the above function is easily verified to be
increasing in $e_2$. \qed\ 

\bigskip
In practice, the area of the cusp torus $\partial B_{\infty}/H$ will
usually be far greater than our crude estimate of ${e_2}^2\sqrt{3}$,
because the density of the packing of $\partial B_{\infty}$ by circles
around the orthocenters of the $\mathcal{O}(1)$-horoballs is typically
nowhere near optimal. One strategy to try and capture the extra area
between the circles is simply to use larger circles: specifically, use
circles of radius $e_3/2$ instead of $e_2/2$ to estimate the
area. This has the advantage of incorporating $e_3$ into the area
estimate, which will lead to the upper bound on $e_3$ that we
seek. The disadvantage, of course, is that by increasing the radius we
can no longer guarantee that each circle will be embedded in $\partial
B_{\infty}/H$: the larger circles may overlap.

However, such overlaps imply the existence of triples which may be
part of a geometric Mom-$n$ structure. If two circles of radius
$e_3/2$ overlap then the distance between their centers is less than
$e_3$, which implies by Lemma \ref{lem:eucl_dist} that the
corresponding horoballs, together with $B_{\infty}$, must form a
$(1,1,2)$-triple. If there are no such triples, then our new expanded
circles will not overlap. This leads to the following:

\begin{lemma}
Let $N$ be a one-cusped hyperbolic three-manifold whose horoball diagram
contains no triples of type $(1,1,2)$. Then
\begin{equation}
\Vol(N)\ge \frac{{e_2}^2{e_3}^2\sqrt{3}}{2}-\pi\left(
{e_2}^2-1-2\log{e_2}\right)  \label{eq:2nd_vol_122}%
\end{equation}
\label{lem:2nd_vol_122}
\end{lemma}

\begin{corollary}
Under the assumptions of the lemma, if $e_3>1.8135$ then
$\Vol(N)>2.848$. \label{cor:2nd_vol_122}
\end{corollary}

\noindent\emph{Proof:} Let $\{p_{i}\}\subset\partial B_{\infty}$ be defined as
in the proof of Theorem \ref{lem:first_vol}. Since the horoball diagram for
$N$ contains no $(1,1,2)$ handles, $d_{E}(p_{i},p_{j})\ge e_3$ for all
$i\not =j$, and $\partial B_{\infty}$ can be packed by circles of radius
$e_3/2$ centered at each point $p_{i}$. The rest of the proof proceeds just
as in the proof of Theorem \ref{lem:first_vol}.

To prove the corollary, note that the right-hand side of equation
(\ref{eq:2nd_vol_122}) is both increasing in $e_2$ for fixed $e_3$ and
vice-versa. When $e_3=1.8135$ and $e_2=1$, $\Vol(N)>2.848$ and
the result follows. \qed\ 

\bigskip
The upper bound on $e_3$ obtained from Corollary \ref{cor:2nd_vol_122}
is unfortunately too large to be useful. The next step is to improve
upon it by accounting for the effect of the horoballs in
$\mathcal{O}(2)$ upon $\Area(\partial B_{\infty})$. If $e_2$ is small,
then we expect to be able to construct circles around the orthocenters
of these horoballs which will be almost as large as the circles around
the orthocenters of the horoballs in $\mathcal{O}(1)$, increasing our
estimate. While if $e_2$ is large, then $e_3$ must be large as well,
increasing the area of the circles constructed in the proof of Lemma
\ref{lem:2nd_vol_122}.

Thus we wish to estimate $\Area(\partial B_{\infty})$
using two sets of circles: one set corresponding to the elements of
$\mathcal{O}(1)$ as before, and another set corresponding to the elements of
$\mathcal{O}(2)$. But we need to decide how large the circles in this new
second set will be. Also, whatever new circles we
construct may intersect the previously constructed circles corresponding to
the $\mathcal{O}(1)$ horoballs. In particular, if $N$ possesses a
$(1,2,2)$-triple then $\mc{O}(2)$-horoballs and $\mc{O}(1)$-horoballs
will be close enough for these circles to possibly intersect. This can be
accounted for using the function $\overlapArea(a,b,c)$ from Lemma
\ref{lem:overlapArea}, leading to the following:

\begin{lemma}
Let $N$ be a one-cusped hyperbolic $3$-manifold containing
no triples of type $(1,1,2)$ and at most one triple of type
$(1,2,2)$. Furthermore, assume that $e_2\le 1.4751$ and $e_3\le
1.8135$. Then
\begin{align*}
\Area(\partial B_{\infty}/H) \ge & 2 \pi \left(\frac{e_3}{2}\right)^2
+ 2 \pi \left(\frac{e_3}{e_2}-\frac{e_3}{2}\right)^2 \\
& - 2 \overlapArea\left(\frac{e_3}{e_2}-\frac{e_3}{2}, \frac{e_3}{2},
1\right) \\
& - \overlapArea\left(\frac{e_3}{e_2}-\frac{e_3}{2},
\frac{e_3}{e_2}-\frac{e_3}{2}, \frac{1}{{e_2}^2}\right)
\end{align*}
and furthermore
\begin{equation}
\Vol(N) \ge \Area(\partial B_{\infty}/H) \frac{{e_3}^2}{2}
- \pi\left({e_3}^2 - 1 - 2\log e_3 + \frac{{e_3}^2}{{e_2}^2} - 1
           - 2\log{\frac{e_3}{e_2}}\right)
\label{eq:3rd_vol_122}\end{equation}
\label{lem:3rd_vol_122}\end{lemma}

\begin{corollary}
Under the assumptions of Lemma \ref{lem:3rd_vol_122} if $e_3>1.4751$
then $\Vol(N)>2.848$.
\label{cor:3rd_vol_122}\end{corollary}

\noindent\emph{Proof:} Let $\{p_{i}\}$ be the set of orthocenters of
$\mathcal{O}(1)$ horoballs as before, and let $\{q_{i}\}$ be the set of
orthocenters of $\mathcal{O}(2)$ horoballs. Assume that we have already
constructed circles of radius $e_3/2$ around each $p_{i}$ as in Lemma
\ref{lem:2nd_vol_122}. As before, since there is not a $(1,1,2)$-triple these
circles will not overlap.

Now in addition to this, construct circles of radius $e_3/e_2-e_3/2$
around each of the points $q_{i}$. These circles may conceivably overlap
either the previous circles or each other. Suppose that the circle around
$q_{i}$ overlaps the circle around $p_{j}$ for some $i$ and $j$, and let
$B_{i}$ and $B_{j}$ be the corresponding horoballs. Then clearly $d_{E}%
(q_{i},p_{j})< e_3/e_2$, which by Lemma \ref{lem:eucl_dist} implies that
$d(B_{i},B_{j})<o(3)$. Therefore $(B_{i},B_{j},B_{\infty})$ must be a triple
of type $(1,2,2)$ (or a triple of type $(1,1,2)$, but we're explicitly
excluding that case right now).

Or suppose that the circles around $q_{i}$ and $q_{j}$ overlap for some $i$
and $j$, and let $B_{i}$ and $B_{j}$ be the corresponding horoballs. If
$d(B_{i},B_{j})\ge o(3)$, then by Lemma \ref{lem:eucl_dist} we must have
$d_{E}(q_{i},q_{j})\ge e_3/{e_2}^2$. But it is easy to show that
$2(e_3/e_2-e_3/2)\le e_3/{e_2}^2$ for all positive $e_2$ and
$e_3$; therefore if the circles are indeed overlapping then $d(B_{i}%
,B_{j})<o(3)$. Therefore $(B_{i},B_{j},B_{\infty})$ must again be a triple of
type $(1,2,2)$ (recall that a triple of type $(2,2,2)$ is impossible by Lemma
\ref{lem:no_bad_triples}).

In summary, any overlaps between the new circles and the old ones, or between
the new circles and each other, arise due to the presence of a triple of type
$(1,2,2)$. And we have supposed that there is no more than one such triple up
to the action of $G$.

After taking the quotient by the action of $H$ we are left with two
new circles of radius $e_3/e_2-e_3/2$ and up to three new cases where
one circle overlaps another. First, a $(1,2,2)$-triple implies that
$d_{E}(q_{i},q_{j})=1/{e_2}^2$ for some $i$ and $j$; therefore in
$\partial B_{\infty}/H$ we may see either the two new circles
overlapping each other or else we may see one of the new circles
overlapping itself. Secondly, a $(1,2,2)$-triple implies that
$d_{E}(q_{i},p_{j})=1$ for some $i$ and $j$; in $\partial
B_{\infty}/H$ we may see up to two instances of a new circle being
overlapped by an old one.

Therefore the area of $\partial B_{\infty}/H$ is at least the area of the two
new circles plus the two old circles minus the three possible
overlaps. Using Lemma \ref{lem:overlapArea}, this proves the first
half of the lemma.

There is a technical issue that must be addressed: Lemma
\ref{lem:overlapArea} is not valid if the circles in question do not
overlap, or if one circle is contained in the interior of the
other. Since we're using Lemma \ref{lem:overlapArea} twice we need to
ensure that both of the following sets of inequalities hold:
\begin{eqnarray*}
\left|\frac{e_3}{e_2}-e_3\right| \le & 1 & \le \frac{e_3}{e_2} \\
0 \le & {e_2}^{-2} & \le \frac{2e_3}{e_2}-e_3 \\
\end{eqnarray*}
All of the above inequalities can be verified by elementary means when
$1\le e2\le 1.4751$ and $e2\le e3\le 1.8135$. Hence our use of Lemma
\ref{lem:overlapArea} is valid. Note this also confirms that the new
circles actually contribute to the area of $\partial B_{\infty}/H$.

To find a lower bound of the volume of $N$, we inflate the cusp
neighborhood as in the proof of Theorem \ref{lem:first_vol}. This time
we obtain a horoball $C$ centered at infinity with Euclidean height
$1/e_3$.  Then $\Vol(C/H)=\frac{1}{2}\operatorname{Area}(\partial
B_{\infty}/H){e_3}^2$. And while $C/H$ extends beyond the Ford
domain corresponding to $B_{\infty}/H$, at worst it only extends past
the Ford faces corresponding to $\mathcal{O}(1)$-horoballs and
$\mathcal{O}(2)$-horoballs.  Hence,
\[
\Vol(N) \ge \Area(\partial B_{\infty}/H)\frac{{e_3}^2}{2}
 - 2\lessvol(1,e_3)-2\lessvol(e_2,e_3)
\]
This proves the second half of the lemma. To prove the corollary,
note that the resulting volume bound is increasing in
$e_3$ for fixed $e_2$ and decreasing in $e_2$ for fixed $e_3$.
Hence the minimum value of the volume bound over the domain $1\le
e_2\le 1.4751$, $1.4751\le e_3\le 1.8135$ occurs when
$e_2=e_3=1.4751$, and at the point $\Vol(N)>2.848$. \qed\

\bigskip
The next step is to perform the same analysis in the case where the horoball
diagram for $N$ contains exactly one triple of type $(1,1,2)$ and no triples
of type $(1,2,2)$. (Note that if $N$ had at least one of each type of
triple then $N$ would have a geometric Mom-$2$.) This analysis is
similar to that of the previous case, and therefore in what follows some
details are omitted.

\begin{lemma}
Let $N$ be a one-cusped hyperbolic $3$-manifold
whose horoball diagram contains one triple of type $(1,1,2)$.
Then
\[
\Area(\partial B_{\infty}/H) \ge 2\pi\left(\frac{e_3}{2}\right)^2
-\overlapArea\left(\frac{e_3}{2}, \frac{e_3}{2}, e_2\right)
\]
and furthermore
\begin{equation}
\Vol(N) \ge \Area(\partial B_{\infty}/H) \frac{{e_2}^2}{2}
-\pi \left({e_2}^2-1-2\log e_2\right)
\label{eq:2nd_vol_112}\end{equation}
\label{lem:2nd_vol_112}\end{lemma}

\begin{corollary}
Under the same assumptions as in Lemma \ref{lem:2nd_vol_112}, if
$e_3>2.1491$ then $\Vol(N)>2.848$.
\label{cor:2nd_vol_112}
\end{corollary}

\noindent\emph{Proof:} As in the proof of Lemma \ref{lem:2nd_vol_122},
we wish to tile $\partial B_{\infty}/H$ with disks of radius $e_3/2$
centered at the centers of the $\mathcal{O}(1)$-horoballs. But since
in this case there is by assumption a single $(1,1,2)$-triple, such
disks will overlap exactly once. This, together with Lemma
\ref{lem:overlapArea}, proves the first part of the lemma. (It is a
trivial matter to confirm that \ref{lem:overlapArea} applies; the
relevant inequality is $0\le e2\le e3$.) To prove
the second part, we inflate the cusp to obtain a horoball $C$ centered
at infinity and with Euclidean height $1/e_2$, and then proceed just
as in Theorem \ref{lem:first_vol}. 

The resulting bound on $\Vol(N)$ is increasing in $e_3$ for fixed
$e_2$ and vice-versa, and when $e_2=1$ and $e_3=2.1491$ we get
$\Vol(N)>2.848$; this proves the corollary.\qed\

\bigskip
We now wish to improve the bound on $e_3$ by mimicking the argument
used in the previous case. That is, we wish to
construct additional circles in the horoball diagram of radius
$e_3/e_2-e_3/2$ corresponding to the $\mc{O}(2)$ horoballs as
before. As in the previous case, such circles are small enough that
the only overlaps between them and the circles of radius $e_3/2$ will
arise as a result of the presence of a $(1,1,2)$-triple (of which
we assume there is at most one) or a $(1,2,2)$-triple (of which we
will assume there are none at all).

Here we run into a problem, however, when we try to confirm that Lemma
\ref{lem:overlapArea} applies, or in other words when we  try to
confirm that the new circles are both overlapped by the old ones and
not completely contained within the old ones. According to Lemma
\ref{lem:eucl_dist}, in the presence of a $(1,1,2)$ triple the
orthocenter of a $\mc{O}(2)$ horoball will be at a distance of
$1/e_2$ from
the center of some $\mc{O}(1)$ horoball. Therefore Lemma
\ref{lem:overlapArea} applies if and only if the following
inequalities hold:
\[
\left|\frac{e_3}{e_2}-e_3\right| \le \frac{1}{e_2} \le \frac{e_3}{e_2}
\]
The right-hand inequality is trivially true since $e_3\ge 1$ but the
left-hand inequality simplifies to $e_3 \le (e_2-1)^{-1}$. This last
inequality is \emph{not} always true in the region $1\le e_2\le
1.4751$, $e_2\le e_3\le 2.1491$. So we have to be a little more
clever. Note that the inequality $e_3 \le (e_2-1)^{-1}$ fails when
$e_3$ and $e_2$ are both large, but the worst case in the previous
lemma occurred when $e_2$ was small. Hence a more sophisticated version
of Corollary \ref{cor:2nd_vol_112} is required, one which lets us
restrict our attention to a smaller region which does not intersect
the curve $e_3=(e_2-1)^{-1}$:

\begin{corollary}
Under the assumptions of Lemma \ref{lem:2nd_vol_112}, if
$e_3>2.1491-(e_2-1)$ then $\Vol(N)>2.848$.
\label{cor:2nd_vol_112_better}
\end{corollary}

\noindent\emph{Proof:} This follows from two observations. First, that
the volume estimate in Lemma \ref{lem:2nd_vol_112} is increasing in
$e_3$ for fixed $e_2$. Second, if $e_3=2.1491-(e_2-1)$ then
$\Vol(N)>2.848$, as can be readily verified by direct computation.\qed\

\bigskip
Note that the curve $e_3=(e_2-1)^{-1}$ lies above the line
$e_3=2.1491-(e_2-1)$. Thus with the preceding corollary in hand, we
can now prove the following:

\begin{lemma}
Let $N$ be a one-cusped hyperbolic $3$-manifold containing no triples
of type $(1,2,2)$ and one triple of type $(1,1,2)$. Furthermore,
assume that $e_2\le 1.4751$ and $e_3\le 2.1491-(e_2-1)$. Then
\begin{eqnarray*}
\Area(\partial B_\infty/H) &\ge& 2\pi\left(\frac{e_3}{2}\right)^2
+ 2\pi\left(\frac{e_3}{e_2}-\frac{e_3}{2}\right)^2
- \overlapArea\left(\frac{e_3}{2}, \frac{e_3}{2}, e_2\right)\\
&&- 2\overlapArea\left(\frac{e_3}{e_2}-\frac{e_3}{2}, \frac{e_3}{2},
\frac{1}{e_2}\right)
\end{eqnarray*}
and furthermore
\[
\Vol(N)\ge \Area(\partial B_\infty/H) \frac{{e_3}^2}{2}
- \pi\left({e_3}^2-1-2\log e_3 + \frac{{e_3}^2}{{e_2}^2} -1
- 2\log\frac{e_3}{e_2} \right)
\]
\label{lem:3rd_vol_112}\end{lemma}
\begin{corollary}
Under the assumptions of Lemma \ref{lem:3rd_vol_112}, if $e_3>1.5152$
then $\Vol(N)>2.848$.
\label{cor:3rd_vol_112}\end{corollary}

\noindent\emph{Proof:} As promised this proof will be somewhat light
on details due to the extreme similarity to the previous case.
Let $\{p_i\}$ be the orthocenters of the
$\mc{O}(1)$ horoballs and $\{q_i\}$ the orthocenters of the
$\mc{O}(2)$ horoballs as before. Construct circles of radius $e_3/2$
around the $p_i$'s. As in Lemma \ref{lem:2nd_vol_112} there is
one overlap between those two disks, caused by the $(1,1,2)$-triple.
In addition, we also have circles of radius $e_3/e_2-e_3/2$ centered
at the $q_i$'s. These circles do not intersect each other since there
are by assumption no $(1,2,2)$-triples in $N$; however since there is
one $(1,1,2)$-triple there will be two overlaps
between these circles and the circles of radius $e_3/2$. This together
with Lemma \ref{lem:overlapArea} and our observations after Corollary
\ref{cor:2nd_vol_112_better} prove the first part of the lemma. To
prove the second part, we inflate the cusp to obtain a horoball $C$
centered at infinity with Euclidean height $1/e_3$, which at worst
extends past the Ford faces corresponding to the
$\mathcal{O}(1)$-horoballs and $\mathcal{O}(2)$-horoballs, then apply
the lessvol function as before. The resulting bound on $\Vol(N)$ is
increasing in $e_3$ for fixed $e_2$, and decreasing in $e_2$ for fixed
$e_3$, and when $e_2=1.4751$ and $e_3=1.5152$ we get
$\Vol(N)>2.848$. This proves the corollary.\qed\

\bigskip
We now finally have enough tools to prove the main result of this
section:

\noindent\emph{Proof of Proposition \ref{prop:param_bounds}:}
Suppose that $N$ is such that $\Vol(N)\le 2.848$. Then
$e_2\le 1.4751$ by Corollary \ref{cor:first_vol}.

Suppose that $N$ does \emph{not} contain a geometric Mom-$2$ of
the type described in the proposition; we
wish to show that this implies $e_3\le 1.5152$. Our assumption
implies that $N$ cannot contain:
\begin{itemize}
\item Two or more $(1,1,2)$-triples
\item Two or more $(1,2,2)$-triples
\item Both a $(1,1,2)$-triple and a $(1,2,2)$-triple.
\end{itemize}
So if $N$ contains no $(1,1,2)$-triples, then it contains at most one
$(1,2,2)$-triple and Corollary \ref{cor:3rd_vol_122}
applies. If $N$ contains exactly one $(1,1,2)$-triple, then it
must contain no $(1,2,2)$-triples and Corollary \ref{cor:3rd_vol_112}
applies, completing the proof.\qed\

\section{Finding a geometric Mom-$n$, $n=2$ or $3$}

The goal of this section is to strengthen the results of the previous
section to obtain the following:
\begin{proposition}
Let $N$ be a one-cusped hyperbolic $3$-manifold with $\Vol(N)\le
2.848$. Then $N$ possesses a geometric Mom-$n$ structure for $n=2$
or $3$ which is torus-friendly.
\label{prop:no_mom3}\end{proposition}

The technique used to prove the above theorem is identical in
principle to the techniques used to prove Proposition
\ref{prop:param_bounds}. Namely, we assume that $N$ does not possess such a
geometric Mom-$n$ structure and then construct a lower bound on
the area of the cusp torus. Specifically we will construct circles
around the orthocenters of horoballs in $\mc{O}(1)$, $\mc{O}(2)$, and
$\mc{O}(3)$. The circles we will construct will be large enough that
they will overlap one another, but we can use the lack of a Mom-$2$ or
Mom-$3$ structure to carefully limit the number of such overlaps that will
occur. Once we have a bound on the area of the cusp torus we can use
that (along with a careful estimate of the volume of the manifold
outside the cusp neighborhood) to find a bound on the volume of
$N$. The difficulties, as we will see, are in the number of cases to
be considered and the complex nature of the volume bounds that
result. But for now, we begin with the following:

\begin{definition}
Let $e_{max}=\min(e_4,1.5152)$, and let
\[
A_0 = \sum_{i=1}^{3} 2 \pi \left(e_{max} ({e_i}^{-1}-\frac{1}{2})\right)^2
\]
\end{definition}

Recall that by the results of the previous section we may assume that
$e_3\le 1.5152<2$; therefore $e_{max}({e_i}^{-1}-1/2)$ is always positive for
$i\le 3$. Furthermore, if $p_i$ and $p_j$ are the orthocenters of
horoballs $B_i\in\mc{O}(i)$ and $B_j\in\mc{O}(j)$, where $i$,
$j\in\{1,2,3\}$, and if $d(B_i,B_j)\ge o(4)$, then
$d_E(p_i,p_j)\ge e_4/e_i e_j$ by Lemma \ref{lem:eucl_dist}. Then we
have the following:
\begin{eqnarray*}
\frac{e_4}{e_i e_j} - e_{max}\left({e_i}^{-1}-\frac{1}{2}\right)
- e_{max}\left({e_j}^{-1}-\frac{1}{2}\right)
&\ge& \frac{e_4}{e_i e_j}(1-e_i)(1-e_j)\\
&\ge& 0
\end{eqnarray*}
where the last line follows since $e_i\ge 1$ for all $i$.

Therefore if we place a circle of radius $e_{max}({e_i}^{-1}-1/2)$ around the
orthocenter of both of the horoballs in $\mc{O}(i)$, for
$i\in\{1,2,3\}$, and if the manifold $N$ does not possess \emph{any}
horoball triples involving $\mc{O}(1)$, $\mc{O}(2)$, and $\mc{O}(3)$,
then those six circles will have disjoint interiors and therefore
$A_0$ will be a lower bound for the area of the cusp torus
$\partial B_\infty /H$. Of course it is highly unlikely that $N$
will possess no such triples, and thus we must consider the
possibility that some of those six circles will overlap. This leads to
our next definition:

\begin{definition}
Let
\[
l_{i,j,k} = \overlapApprox\left(e_{max}({e_i}^{-1}-\frac{1}{2}),
e_{max}({e_j}^{-1}-\frac{1}{2}), c\right)
\]
where $c=\min(e_k/(e_i e_j), e_{max}({e_i}^{-1}+{e_j}^{-1}-1))$.
\end{definition}
Recall that $\overlapApprox(a,b,c)$ was defined in Definition
\ref{def:overlapApprox}, and by Lemma \ref{lem:overlapApprox} it
is always greater than or equal to the area of the intersection of a
circle of radius $a$ and a circle of radius $b$ whose centers are $c$
units apart, provided that $|a-b|\le c\le a+b$.

Hence if $p_i$, $p_j$ are respectively the orthocenters of
horoballs $B_i\in \mc{O}(i)$, $B_j\in \mc{O}(j)$, and if
$d(B_i,B_j)=o(k)$, where $i$, $j$, and $k$ are all in $\{1,2,3\}$,
then $l_{i,j,k}$ will be greater than or equal to the amount of
overlap between a circle of radius $e_{max}({e_i}^{-1}-1/2)$ around
$p_i$ and a circle of radius $e_{max}({e_j}^{-1}-1/2)$ around
$p_j$. We need to check that the condition $|a-b|\le c$ holds. Assuming
that $e_i\le e_j$, we have
\begin{eqnarray*}
c-|a-b| &\ge& \frac{e_k}{e_i e_j}
- e_{max}\left|{e_i}^{-1}-{e_j}^{-1}\right| \\
&=& \frac{1}{e_i e_j}\left(e_k-e_{max}(e_j-e_i)\right) \\
&\ge& \frac{1}{e_i e_j}\left(1-1.5152(1.5152-1)\right) \\
&\ge& 0
\end{eqnarray*}
where the second-to-last line used the fact that $1\le e_i\le e_j\le
1.5152$. Hence we do in fact have $|a-b|\le c$, i.e. it is not the
case that one circle lies entirely inside the other.

Note that it \emph{is} entirely possible that $e_k/(e_i e_j)$ will be
greater than $e_{max}({e_i}^{-1}+{e_j}^{-1}-1)$, i.e. that the two
circles don't overlap at all. In this case $l_{i,j,k}$ will simply
equal $0$ since $\overlapApprox(a,b,a+b)=0$ for all non-negative $a$
and $b$.

Now suppose that the manifold $N$ possesses a $(i,j,k)$-triple where
$i$, $j$, and $k$ are all elements of $\{1,2,3\}$. Then there will
exist orthocenters $p_i$ and $p_j$ corresponding to horoballs
$B_i\in\mc{O}(i)$ and $B_j\in\mc{O}(j)$ such that $d(B_i,B_j)=o(k)$,
and similarly for each cyclic permutation of $i$, $j$, and $k$. In
other words, an $(i,j,k)$-triple with $i$, $j$, $k\in\{1,2,3\}$ can
cause up to three different overlaps in $\partial B_\infty /H$. But
there will be no overlaps between the circles that we have constructed
that do not come from such a triple. In other words we have
established the following:

\begin{lemma}
Suppose that $N$ possesses horoball triples of the form
$(i_1,j_1,k_1)$, \ldots , $(i_s,j_s,k_s)$, where $i_r$, $j_r$,
$k_r\in\{1,2,3\}$ for each $r=1$,\ldots $s$. Furthermore, suppose
every such triple in $N$ appears in the above list exactly once. Then
\[\Area(\partial B_\infty /H) \ge A_0
- \sum_{r=1}^{s} l_{i_r,j_r,k_r}+l_{j_r,k_r,i_r}+l_{k_r,i_r,j_r}
\]
\label{lem:areaFns}
\end{lemma}

Once we have a bound on the area of $\partial B_\infty /H$, we can
compute the volume of $B_\infty/H$ and thereby
construct a bound on the volume of $N$ exactly as in the previous
section. Recall that in the previous section we ``inflated'' the cusp
neighborhood to get a larger neighborhood $C/H$ which was only
immersed in $N$, and then subtracted the volume of the regions where
$C/H$ extended beyond the Ford domain corresponding to $B_\infty
/H$. This allowed us to improve our volume bound by accounting for
some of the volume of $N-B_\infty /H$. We will do the same here,
``inflating'' our cusp neighborhood by replacing $B_\infty$ with a
new horoball $C$ centered at infinity but with Euclidean height
$1/e_{max}$. Then we have the following:

\begin{lemma}
Under the same assumptions as in Lemma \ref{lem:areaFns},
\begin{eqnarray*}
\Vol(N) &\ge& \frac{{e_{max}}^2}{2}
\left(A_0 - \sum_{r=1}^{s}
l_{i_r,j_r,k_r}+l_{j_r,k_r,i_r}+l_{k_r,i_r,j_r}
\right) \\
& &-\pi\left(-3+{e_{max}}^2\left(1+{e_2}^{-2}+{e_3}^{-2}\right)
+\log\left(\frac{{e_2}^2 {e_3}^2}{{e_{max}}^6}\right)\right)
\end{eqnarray*}
\label{lem:volFn}
\end{lemma}

Unfortunately one very quickly discerns that the volume bounds
obtained by the above lemma simply aren't large enough for our
purposes, particularly when $e_2$, $e_3$, and $e_4$ are small. For
example, if $e_2=e_3=e_4=1$ one can quickly determine that $A_0=3
\pi/2$ and that all the $l_{i,j,k}$'s are $0$, resulting in a volume
bound of $3 \pi/4$, which is less than $2.848$. Thus we need to find
more area in the case where $e_4$ is small. Since a small value of
$e_4$ implies that the $\mc{O}(4)$-horoballs will be close to
$B_\infty$, it is natural to try to construct circles around the
orthocenters of the $\mc{O}(4)$-horoballs to increase our estimate of
the area of $B_\infty/H$.

So suppose $e_4<1.5152$ (so $e_{max}=e_4$).  We construct a circle of
radius $1/(e_4 e_2)-e_4/e_2+e_4/2$ around the orthocenter of each of
the two horoballs in $\mc{O}(4)$. The reasoning behind the choice of
this particular radius is as follows. If $q_4$ is the orthocenter of
such a horoball, and if $p_i$ is the orthocenter of $B_i\in\mc{O}(i)$
for $i\in\{2,3,4\}$ (assume $p_i\not =q_4$), we want the circles
constructed around $q_4$ and $p_i$ to have disjoint interiors. Since
$d_E(p_i,q_4)\ge 1/(e_i e_4)$ by Lemma
\ref{lem:eucl_dist}, this means that when $i=2$ or $3$ we require
\[
\frac{1}{e_i e_4}
- \left(\frac{1}{e_4 e_2}-\frac{e_4}{e_2}+\frac{e_4}{2}\right)
-\left(\frac{e_4}{e_i}-\frac{e_4}{2}\right) \ge 0
\]
and when $i=4$ we require
\[
\frac{1}{{e_4}^2}
- 2\left(\frac{1}{e_4 e_2}-\frac{e_4}{e_2}+\frac{e_4}{2}\right) \ge 0
\]
The first inequality is equivalent to $({e_2}^{-1}-{e_i}^{-1})
(e_4-{e_4}^{-1})\ge 0$, which is clearly true for $i=2$ or $3$. The
second inequality simplifies to
\[
\frac{e_4-1}{{e_4}^2 e_2}
\left(2{e_4}^2+2{e_4}-e_2({e_4}^2+e_4+1)\right)
\ge 0
\]
But $1\le e_2\le e_4$, and the polynomial $2x^2+2x-x(x^2+x+1)$ is always
positive for $1\le x\le 1.5152$, and hence each factor on the
left-hand side of the above inequality is always non-negative, proving
the inequality.

Hence our new circles will not intersect each other, nor will they
intersect the previous circles that were created around the
orthocenters of the $\mc{O}(2)$ or $\mc{O}(3)$-horoballs. Now suppose
$q_4$ is the orthocenter of $B_4\in\mc{O}(4)$ as before, and suppose
$p_1$ is the orthocenter of a horoball $B_1\in \mc{O}(1)$. Suppose
further that $d(B_1,B_4)\ge o(2)$, and that hence $d_E(p_1,q_4)\ge
e_2/e_4$ by Lemma \ref{lem:eucl_dist}. We wish to know if the circles
we've constructed around $p_1$ and $q_4$ will overlap; that is, we
wish to verify the inequality
\begin{eqnarray*}
\frac{e_2}{e_4}
-\left(\frac{1}{e_4 e_2}-\frac{e_4}{e_2}+\frac{e_4}{2}\right)
-\frac{e_4}{2} &\ge& 0 \\
\Leftrightarrow\
\frac{e_2-1}{e_4 e_2} \left(e_2+1-{e_4}^2\right) &\ge& 0
\end{eqnarray*}
This inequality is \emph{not} always true in the domain $1\le e_2\le
e_4\le 1.5152$, particularly when $e_2$ is small and $e_4$ is
large. However, the inequality \emph{does} hold when we need it to
hold, namely when $e_2$ and $e_4$ are both small. Hence in what
follows we will make the additional assumption that $e_2+1\ge{e_4}^2$,
and take it for granted that if this assumption fails we don't need the
extra area anyway.

With the additional assumption, we now have circles constructed around
the orthocenters of the $\mc{O}(4)$-horoballs which don't intersect
each other, which don't intersect the circles around the centers of
the $\mc{O}(2)$-horoballs and $\mc{O}(3)$ horoballs, and which may
intersect the circles constructed around a $\mc{O}(1)$ horoball but
only if that horoball is less than $o(2)$ away from the
$\mc{O}(4)$-horoballs. In other words, the new circles do not intersect
themselves or any of the previous circles unless $N$ contains one or more
$(1,1,4)$-triples. If $N$ contains exactly one such triple, there will
be two overlaps to account for, while
if $N$ contains two or more such triples, we
have a geometric Mom-$2$ structure. Thus we can conclude the
following:

\begin{lemma}
If $e_4\le 1.5152$ and $e_2+1\ge {e_4}^2$, then either $N$ contains a
geometric Mom-$2$ structure consisting of two $(1,1,4)$-triples or
else the area estimate of Lemma \ref{lem:areaFns} can be increased by
\[
2\pi\left(
\frac{1}{e_4 e_2}-\frac{e_4}{e_2}+\frac{e_4}{2}
\right)^2-2 \overlapApprox(a,b,c)
\]
where $a=e_4/2$, $b=1/(e_4 e_2)-e_4/e_2+e_4/2$, and $c=1/e_4$.
\label{lem:areaFns2}
\end{lemma}

Note that $a$ and $b$ are the radii of the circles constructed around
the orthocenters of the $\mc{O}(1)$-horoballs and
$\mc{O}(4)$-horoballs respectively, while $c$ is the minimum possible
distance between those orthocenters in the presence of a
$(1,1,4)$-triple. Also we are implicitly assuming that the conditions
of Lemma \ref{lem:overlapApprox} are met, namely that $|a-b|\le c\le
a+b$. Fortunately the left-hand inequality simplifies to
$(e_2+1-{e_4}^2)/(e_2 e_4)\ge 0$, which is true by assumption, while
the right-hand inequality simplifies to $({e_4}^2-1)(e_2-1)/(e_2
e_4)\ge 0$.

We can also update our volume bounds:

\begin{lemma}
If $e_4\le 1.5152$ and $e_2+1\ge {e_4}^2$, then either $N$ contains a
geometric Mom-$2$ structure consisting of two $(1,1,4)$-triples or
else the volume estimate of Lemma \ref{lem:volFn} can be increased by
\[
\frac{{e_4}^2}{2} \left(
2\pi\left(
\frac{1}{e_4 e_2}-\frac{e_4}{e_2}+\frac{e_4}{2}
\right)^2-2 \overlapApprox(a,b,c)\right)
\]
where $a=e_4/e_2$, $b=1/(e_4 e_2)-e_4/e_2+e_4/2$, and $c=1/e_4$.
\label{lem:volFn2}
\end{lemma}

We now have in principle a procedure for proving Proposition
\ref{prop:no_mom3}. First, enumerate every possible combination of
triples which does \emph{not} include a geometric Mom-$2$ or
Mom-$3$ structure which is torus-friendly. Second, for each combination
construct the volume bounds from Lemma \ref{lem:volFn} and Lemma
\ref{lem:volFn2} above, which will be a function of $e_2$, $e_3$, and
$e_4$. Third, show that each such volume bound never attains a value
below $2.848$.

The first step is the simplest, and in fact we can make it even
simpler: we do not need to enumerate every possible combination of
triples that do not include a Mom-$n$, merely the maximal ones, since
adding additional triples only decreases the area bounds in Lemma
\ref{lem:areaFns} and Lemma \ref{lem:areaFns2}. There are a total of
eighteen maximal combinations of triples that must be considered:
\begin{itemize}
\item $(1,1,2)$, $(1,1,3)$
\item $(1,1,2)$, $(1,3,3)$
\item $(1,1,2)$, $(2,2,3)$
\item $(1,1,2)$, $(2,3,3)$
\item $(1,2,2)$, $(1,1,3)$
\item $(1,2,2)$, $(1,3,3)$
\item $(1,2,2)$, $(2,2,3)$
\item $(1,2,2)$, $(2,3,3)$
\item $(1,1,3)$, $(2,2,3)$
\item $(1,1,3)$, $(2,3,3)$
\item $(1,3,3)$, $(2,2,3)$
\item $(1,3,3)$, $(2,3,3)$
\item $(1,1,2)$, $(1,2,3)$, $(1,2,3)$
\item $(1,2,2)$, $(1,2,3)$, $(1,2,3)$
\item $(1,1,3)$, $(1,2,3)$, $(1,2,3)$
\item $(1,3,3)$, $(1,2,3)$, $(1,2,3)$
\item $(2,2,3)$, $(1,2,3)$, $(1,2,3)$
\item $(2,3,3)$, $(1,2,3)$, $(1,2,3)$
\end{itemize}

The second step of the procedure is also simple; it can in fact be
automated by a few lines of Mathematica code (\cite{mil}).

The third step, however, is daunting due to the complicated nature of
the volume bounds that result from Lemmas \ref{lem:volFn} and
\ref{lem:volFn2}. The expressions resulting from these two lemmas
defy analysis by hand.

To handle these complicated expressions, we resort to computer
assistance. Specifically, we use rigorous floating-point arithmetic of
the type used in \cite{gmt}. This is not standard interval
arithmetic, so we take a moment to review the techniques involved
here.

\begin{definition}
An \emph{affine 1-jet} $F=(f_0; f_1,f_2,f_3; f_\epsilon)$ consists of the
set of all
functions $f:[-1,1]^3\rightarrow \BR$ such that $|f(x_1,x_2,x_3)-(f_0 +
\Sigma f_i x_i)|\le f_\epsilon$ for all $(x_1,x_2,x_3)\in[-1,1]^3$. (Note
that we require $f_\epsilon\ge 0$.)
\label{defn:1_jet}
\end{definition}

Note that in \cite{gmt} complex 1-jets were used; however we only
require real numbers here. \cite{gmt} showed how, given two affine
1-jets $F$ and $G$, a computer which meets IEEE standards for
floating-point arithmetic can compute an affine 1-jet $H=(h_0;
h_1,h_2,h_3; h_\epsilon)$ which
``equals $F+G$'', in the sense that $f+g\in H$ for each $f\in F$ and
$g\in G$. Specifically, let $h_i=f_i+g_i$ for $i=0$, $1$, $2$, and
$3$, and let
\[
h_\epsilon = (1+\epsilon_a) (\epsilon_t+\epsilon_f)
\]
where
\begin{eqnarray*}
\epsilon_t &=& f_\epsilon + g_\epsilon \\
\epsilon_f &=& \frac{EPS}{2}
\left((|f_0+g_0|+|f_1+g_1|)+(|f_2+g_2|+|f_3+g_3|)\right) \\
\epsilon_a &=& 3 EPS
\end{eqnarray*}
Here $\epsilon_t$ accounts for the ``Taylor error'', i.e. the maximum
possible sup-norm distance between $f+g$ and the linear function
represented by $h$ as $f$ and $g$ vary over $F$ and $G$
respectively. Since we're only computing a linear function of $f$ and
$g$ anyway, $\epsilon_t$ only needs to account for the error carried
over from the operands $F$ and $G$; for more complicated operations
this term will be more significant. The quantity $\epsilon_f$
accounts for the floating-point error that may accrue from the calculation of
$h_0$ through $h_3$. Here $EPS$ is a (small) computer-dependant constant that
measures the granularity of the set of real numbers that the computer
is capable of representing. Roughly, $1+EPS$ will be the smallest real
number strictly greater than $1$ which has a floating-point
representation on the computer in question. Finally $\epsilon_a$
accounts for the floating-point error that may accrue from calculating
$\epsilon_t+\epsilon_f$. For any operation involving affine $1$-jets,
$\epsilon_a$ will always be of the form $n EPS$ where $n$ is an
integer roughly proportional to the base-$2$ logarithm of the number
of arithmetic operations necessary to compute
$\epsilon_t+\epsilon_f$. Constructed in this way, the error term
$h_\epsilon$ will be large enough to account
for the original error terms $f_\epsilon$ and $g_\epsilon$ and for the
floating-point error that might accrue from calculating the terms of
$H$. Similar constructions exist for ``$-F$'', ``$FG$'', and ``$F/G$''
provided the range of $g$ does not contain $0$ for any $g\in G$.

For more specific details on this process, we refer the reader to
Sections 5 and 6 of \cite{gmt} which describes the theory behind
these formulas and provides numerous examples.


Now suppose that we have a rational polynomial
$p(x_1,x_2,x_3)$ and we wish to compute the range of possible values
of $p$ over a box $I_1\times I_2\times I_3$ in $\BR^3$. It is not
possible to compute the exact range of possible values by computer due
to floating-point error. However affine 1-jets do make it possible to
rigorously determine an interval which must contain the range of
possible values, as follows. Note that for clarity's sake we start
with a simplified version of the procedure and fill in certain
troublesome details later.

Define $X_1$ to be the affine
1-jet $(a_1;b_1,0,0;0)$ where $x\mapsto a_1+b_1x$ is the unique
nondecreasing linear bijection from
$[-1,1]$ to the interval $I_1$. Define $X_2=(a_2;0,b_2,0;0)$ and
$X_3=(a_3;0,0,b_3;0)$ similarly. Then compute the affine 1-jet
$p(X_1,X_2,X_3)$ using the constructions in \cite{gmt} in place of the
usual arithmetic operations, and let $P=(p_0;p_1,p_2,p_3;p_\epsilon)$
denote the result. If $(x_1,x_2,x_3)\in I_1\times I_2\times I_3$, then
trivially
\[
(x_1,x_2,x_3)=(f_1(u,v,w),f_2(u,v,w),f_3(u,v,w))
\]
where $(u,v,w)\in[-1,1]^3$ and $f_i\in X_i$ for $i=1$, $2$, and
$3$. (Specifically, we may choose $f_i(x)=a_i+b_i x$.) Therefore
$p(x_1,x_2,x_3)$ must lie in the range of $p(f_1,f_2,f_3)$, which is
an element of $P$. And therefore
\[
p(x_1,x_2,x_3) \in [p_0-|p_1|-|p_2|-|p_3|-p_\epsilon,
p_0+|p_1|+|p_2|+|p_3|+p_\epsilon]
\]

In practice, there are several complications. First, the above
construction implicitly assumes that the coefficients $a_i$ and $b_i$
have exact binary representations. If they do not, then the affine
1-jets $X_i$ must be modified to have a non-zero $\epsilon$-term,
representing the sup-norm distance between $a_i+b_i x$ and the actual
unique nondecreasing linear bijection from $[-1,1]$ to $I_i$. In
practice it is usually simpler to replace $I_i$ with a slightly larger
interval whose endpoints do in fact have exact binary representations.
Second, computing whether or not $p(x_1,x_2,x_3)$ lies in the interval
described is itself a floating-point operation, and may introduce
error. This can be dealt with using similar techniques to
those used in calculating ``F+G'' in the first place. One final
complication is that this whole technique will only be effective if
$p$ is very simple or if the intervals $I_i$ are very small. If $p$ is
complicated (and the functions we're interested in are very
complicated) then it is usually necessary to subdivide $I_1\times
I_2\times I_3$ into much smaller sub-boxes to achieve any kind of
accuracy.

Nevertheless, we wish to apply the above techniques to compute a range
of possible values for the volume bounds produced by Lemmas
\ref{lem:volFn} and \ref{lem:volFn2}; if the computed ranges never
include values less than or equal to $2.848$ then we'll be
done. Unfortunately those volume bounds are not expressed as rational
polynomials; we also need to be able to rigorously compute both
natural logarithms and the minimum function, two operations that were
never implemented in \cite{gmt}.

Fortunately while the natural logarithm is not a rational polynomial
all of its derivatives are. This makes it possible to rigorously
compute logarithms by using a Taylor approximation to $\log x$ at
$x=1$, and using Taylor's theorem to compute an exact upper bound on
the difference between $\log x$ and the polynomial
approximation. Since Taylor's theorem expresses this difference in
terms of the derivatives of $\log x$ it is possible to incorporate it
into the ``Taylor error term'' $\epsilon_t$ of the resulting affine
1-jet and still determine correct upper bounds for the corresponding
$\epsilon_f$ and $\epsilon_a$ terms. Repeated use of the relation
$\log (ax)=\log a+\log x$ (where $a\not =1$ is just an arbitrary
positive constant; $9/8$ was used in our implementation) allows the
program to restrict its use of the Taylor approximation to regions
close to $x=1$, where the approximation is the most accurate.

It is tempting at this point to try and implement other non-polynomial
functions in this way. For example, consider the function
$\overlapArea(a,b,c)$ from Definition \ref{def:overlapArea}: it is
composed of polynomials, square roots, and the function
$\acos(x)$. Square roots were successfully implemented using affine
1-jets in \cite{gmt}, and the derivatives of $\acos(x)$ are square
roots of rational polynomials, so in theory it is possible to
implement the function $\overlapArea$ with affine 1-jets. In
practice, unfortunately, this works poorly. The function
$f(x)=\acos(x)-x\sqrt{1-x^2}$ which is used in the calculation of
$\overlapArea$ is composed of two functions which have vertical
tangent lines at the points $x=\pm 1$. When computing a function using
affine 1-jets, the size of the resulting error term will always be
proportional to the derivative of the function being computed; as a
consequence, computing $f(x)$ for affine 1-jets which include $x=\pm
1$ in their range in practice causes the error term to grow so large
as to make the entire 1-jet useless. This is the entire reason for the
existence of Lemma \ref{lem:overlapApprox}. While
$\overlapApprox(a,b,c)$ is only an approximation to
$\overlapArea(a,b,c)$, as a polynomial function it is far more useful
in this kind of computation.

Continuing, there remains one function to be implemented: a
function $Max0(F)$ such that if $f\in F$ and $H=Max0(F)$ then
$\max(f,0)\in H$. (Then the identity $\min(f,g)=f-\max(f-g,0)$ can be
used to compute minimums.) If the computer can rigorously determine that $0$
does not lie in the range of $f$ for any $f\in F$ then $Max0(F)$
equals either $F$ itself if $f_0>0$ or else $Max0(F)=(0;0,0,0;0)$ if
$f_0<0$. If the computer cannot exclude the possibility that $0$ lies
in the range of some $f\in F$ then $Max0(F)$ is defined to be equal to
the 1-jet $(s;0,0,0;s)$ where
\[
s = \frac{1}{2}(1+3 EPS)(f_0+((|f_1|+|f_2|)+(|f_3|+f_\epsilon))
\]
(The factor of $(1+3 EPS)$ is there to
account for floating-point error that may accrue during the
calculation of the rest the expression.) Note that while this
definition is technically correct in that $\max(f,0)$ will lie in this
affine $1$-jet for all $f\in F$, from a practical standpoint it is a
terrible definition as $s$ is almost guaranteed to be much, much larger than
the original error term $f_\epsilon$. Fortunately (and unlike the case
with $\acos(x)$ and square roots, above) in practice the
cases where the $Max0$ function had to be called on affine 1-jets that
may have contained $0$ in their ranges were rare enough not to cause
significant problems.

With these tools in hand we now consider the following:

\begin{theorem}
Suppose $N$ is a one-cusped hyperbolic 3-manifold with $\Vol(N)$ $\le
2.848$. Then $e_2\le 1.4751<1.5152$ and furthermore $N$ contains
a geometric Mom-$n$ structure which either
\begin{itemize}
\item is a geometric Mom-$2$ structure incorporating only the orthoclasses
  $\mc{O}(1)$ and $\mc{O}(2)$ and triples constructed from those three
  orthoclasses, or
\item incorporates only the orthoclasses $\mc{O}(1)$, $\mc{O}(2)$, and
  $\mc{O}(3)$, and triples constructed from those three orthoclasses,
  and is torus-friendly, and in addition we have $e_3\le 1.5152$, or
\item incorporates the orthoclasses $\mc{O}(1)$ and $\mc{O}(4)$ and a
  pair of $(1,1,4)$-triples, and in addition we have $e_4\le 1.5152$.
\end{itemize}
\label{thrm:comb_existence}
\end{theorem}

The proof of Theorem \ref{thrm:comb_existence}, using rigorous
computer assistance, proceeds as follows. If $e_2>1.4751$ then
$\Vol(N)>2.848$ by the results of the previous chapter. So, suppose
otherwise. Furthermore, suppose that $N$ does not contain a geometric
Mom-$n$ of one of the three types described. In particular, by
Proposition \ref{prop:param_bounds} we may assume that $e_3\le
1.5152$, and we may assume that $N$ contains at most one
$(1,1,4)$-triple unless $e_4>1.5152$.

Then let $(i_1,j_1,k_1)$, \ldots, $(i_s,j_s,k_s)$ be the complete list
of triples in $N$ satisfying the condition that $i_r$, $j_r$,
$k_r\in\{1,2,3\}$ for all $r$. This list of triples must be a
subcollection of one of the eighteen collections listed above. For
each of those eighteen collections of triples, we can construct a
lower bound on the volume of $N$ via Lemma \ref{lem:volFn} which we
will call $f_1(e_2,e_3,e_4)$. In addition, Lemma \ref{lem:volFn2}
provides a second, stricter lower bound in the case where $e_4\le
1.5152$ and $e_2+1\ge {e_4}^2$; call this function
$f_2(e_2,e_3,e_4)$. These bounds depend solely on the parameters
$e_2$, $e_3$, and $e_4$ and  can be computed using only the four basic
arithmetic operations, logarithms, the minimum function, and some
constants. Furthermore we may assume that $(e_2,e_3,e_4)$ lies in a
compact subset of $\BR^3$. Technically we have not established an
upper bound on $e_4$, but if $e_4>1.5154$ then we may simply replace
$e_4$ with $1.5152$ in $f_1$ and still obtain a valid lower bound on
the volume of $N$. We cannot do the same with $f_2$, but in practice
$f_2$ is only needed to improve the volume bound for small values of
$e_4$. Subdividing this compact region into a sufficiently small
number of pieces and using the rigorous floating-point arithmetic
techniques described above, we can rigorously demonstrate that
$\max(f_1,f_2)$ never attains a value less than or equal to $2.848$.

This approach has been successfully implemented (\cite{mil}) and used
to prove Proposition \ref{prop:no_mom3}. The resulting program
requires approximately 80 minutes to establish that
$\max(f_1,f_2)>2.848$ on the parameter space $1\le e_2\le 1.4751$,
$e_2\le e_3 \le e_4\le 1.5152$ in each of the 18 cases. To do so the
program subdivides each dimension of the parameter space into $2^8$
subintervals. The program also establishes that $f_1>2.848$ whenever
$1\le e_2\le 1.4751$, $e_2\le e_3\le 1.5152$, and $e_4=1.5152$,
requiring approximately 10 seconds to do so; the number of dimensions
in the parameter space has a tremendous affect on the running
time. The programs themselves consist of just
under 4800 lines of C++ code (close to half of which was generated
automatically by a short Mathematica program), running in Redhat Linux
on a four-processor PC.\qed

\section{Embedding a geometric Mom-$n$}

Theorem \ref{thrm:comb_existence} establishes conditions under which
we can assume that a one-cusped hyperbolic $3$-manifold $N$ possesses
a geometric Mom-$2$ or Mom-$3$ structure. Now what we wish to do
is thicken the cellular complex $\Delta$ associated to that structure
to obtain a full topological internal Mom-$n$ structure as defined in
\cite{gmm2}. Examining that definition, we see that there are three
obstacles to this.  The first is
that there is no guarantee that $\Delta$ is embedded in $N$; this is
the obstacle we will address in this section. The second obstacle is
that the complement of $\Delta$ in $N$ may have components whose
boundaries are not tori, which violates the definition of an internal
Mom-$n$ structure; this obstacle will be addressed in the next
section.  The final obstacle, namely determining whether or not
$\Delta$ is ``full'', will be tackled in Section 8. We proceed with
the question of embeddedness now; specifically we will to prove the
following:

\begin{theorem}
In Theorem \ref{thrm:comb_existence}, we may assume that the
geometric Mom-$n$ structure obtained has the property that the
cellular complex $\Delta=T\cup\{\lambda_i\}\cup\{\sigma_j\}$ is
embedded, where $\{\lambda_i\}$ is the set of 1-cells which are the
projections in $N$ of the shortest arcs joining the horoballs which
form the orthopair classes of the Mom-$n$, and where $\{\sigma_j\}$ is
the set of 2-cells which are the projections in $N$ of the totally
geodesic surfaces spanning the triples of horoballs in the
geometric Mom-$n$ structure and the arcs between them.
\label{thrm:comb_embedded}
\end{theorem}

For sake of notation, define the \emph{$\mathcal{O}(n)$-edge} to be
the image in $N$ of the shortest arc between any two horoballs which
constitute an element of $\mathcal{O}(n)$.  All of the lemmas in this
chapter will assume that $N$ is a one-cusped hyperbolic $3$-manifold.

\begin{lemma}
If the $\mc{O}(n)$-edge intersects the cusp torus for any $n$, then
$e_n>1.5152$.
\end{lemma}

\noindent\emph{Proof:} Any edge which intersects the cusp torus lifts
to an arc contained in a line which intersects the boundary of
$B_{\infty}$, neither of whose endpoints are at infinity. Such a line
clearly must be a half-circle of diameter $\ge 2$. But such a line
must also join the centers of two horoballs and hence by Lemma
\ref{lem:eucl_dist} we must have $e_{k}/(e_{m}e_{n})\ge 2$ for some
$k\in\{1,2,3\}$ and for some $m$ and $n$.  This implies that $e_{k}\ge
2>1.5152$. \qed\

\begin{lemma}
If any two-cell in $\Delta$ corresponding to an $(n,m,k)$-triple
intersects the cusp torus then one of $e_n$, $e_m$, or $e_k$ is
greater than $1.5152$.
\label{lem:ball_intersect}
\end{lemma}

\noindent\emph{Proof:} A two-cell which intersects the cusp torus
lifts to a totally geodesic two-cell contained in an ideal hyperbolic
triangle which intersects $B_{\infty}$, whose ideal vertices are the
centers of three horoballs $A$, $B$, $C$, with the property that
$d(A,B)=o(n)$, $d(B,C)=o(m)$, and $d(C,A)=o(k)$. There are two
possibilities. If one of the three edges of this ideal triangle also
intersects $B_{\infty}$, then the previous lemma applies and we're
done. So suppose none of the three edges intersect $B_{\infty}$, but
instead the triangle intersects $B_{\infty}$ at some point in its
interior.

Let $p$ be the highest point of the ideal triangle as viewed in the
upper half-space model (here ``highest'' refers to the Euclidean or
visual height), and let $q$ be the highest point of the entire
hyperbolic plane containing this triangle. Clearly if $p\not =q$ then
$p$ must be some point on an edge of the ideal triangle as close to
$q$ as possible, contradicting our supposition. So $p=q$ and $q$ is
contained in the ideal triangle. Now view the entire picture from the
point at infinity: from above, the ideal triangle projects onto a
Euclidean triangle in $B_{\infty}$ joining the orthocenters of $A$,
$B$, and $C$, and $q$ projects to the circumcenter of this
triangle. Since the projection of $q$ lies in the interior of the
triangle, and since the circumradius of the triangle is $\ge 1$, one
of the sides of the triangle must have length greater than or equal to
$\sqrt {2-2\cos\frac{2\pi}{3}}=\sqrt{3}$. Wolog, assume that the side
``from B to C'' has this property. Therefore by Lemma
\ref{lem:eucl_dist} as before we have $e_{k}/(e_{b}e_{c})
\geq\sqrt{3}$ for some $k\in\{1,2,3\}$ and some $b$ and $c$, and hence
$e_k\ge \sqrt{3}>1.5152$.\qed\

\bigskip
The next step is to establish whether or not any of the
$\mathcal{O}(n)$-edges contained in $\Delta$ can intersect one
another. Some of the following lemmas are stronger than necessary to
prove that this is not the case, but the stronger results will be used
later in the argument. We begin with some definitions to simplify
notation later:

\begin{definition} Let $\delta_1=0.15$,
$\delta_2=\cosh^{-1}(1.5152^{-1}+1.5152^{-2})=0.4337\ldots$, and
$\delta_3=\cosh^{-1}(2/1.5152)=0.7800\ldots$.
\end{definition}

\bigskip
Note that the value of $\delta_1$ comes from Lemma
\ref{lem:dist_line_self}. The values of $\delta_2$ and $\delta_3$ are
motivated by the following application of Corollary
\ref{cor:two_line_dist}:

\begin{lemma}
Suppose $A$, $B$, $C$, $D$, and the various distances between them are
defined as in Corollary \ref{cor:two_line_dist}. In particular suppose
that $d(A,C)=o(m)$, $d(B,D)=o(n)$, and suppose that $x$ is the
shortest distance between the line joining the centers of $A$ and $C$
and the line joining the centers of $B$ and $D$. Suppose further that
$\max(e_m,e_n)\le 1.5152$.  If $x<\delta_2$ then all four of the
distances $o(h)=d(A,B)$, $o(j)=d(B,C)$, $o(k)=d(C,D)$, and
$o(l)=d(D,A)$ are strictly less than $\min(o(m),o(n))$. If
$x<\delta_3$ then at least three of those four distances are strictly
less than $\min(o(m),o(n))$.\label{lem:one_cell_dist}
\end{lemma}

\noindent\emph{Proof:} Without loss of generality suppose that
$o(m)\le o(n)$, and that therefore $e_m\le e_n\le 1.5152$.

Suppose that $x<\delta_2$ and that
$e_{h}\geq e_{m}$. Then we have the following:
\begin{align*}
e_{m}+1  &  \leq e_{h}e_{k}+e_{j}e_{l}\\
&  =e_{m}e_{n}\cosh x\\
&  <(1.5152^{-1}+1.5152^{-2})e_{m}e_{n}%
\end{align*}
where the middle step used the result of Corollary \ref{cor:two_line_dist}.
Dividing both sides by $e_{m}$ and using the fact that $e_{m}\leq e_{n}%
\le 1.5152$ we get%
\[
1+\frac{1}{1.5152}\le 1+\frac{1}{e_{m}}<(1.5152^{-1}+1.5152^{-2})e_{n}
\le 1+\frac{1}{1.5152} \]
which is a contradiction. Hence $e_{h}<e_{m}$, and similarly for $e_{j}$,
$e_{k}$, and $e_{l}$.

Now suppose that $x<\delta_3$. If $e_{h}\geq e_{m}$ and $e_{k}\geq
e_{m}$, then we have%
\begin{align*}
e_{m}^2+1  &  \leq e_{h}e_{k}+e_{j}e_{l}\\
&  =e_{m}e_{n}\cosh x\\
&  <\frac{2}{1.5152}e_{m}e_{n}%
\end{align*}
Again, divide both sides by $e_{m}$. Then using $e_{n}\le 1.5152$ and the
AM-GM\ inequality we get
\[
2\leq e_{m}+\frac{1}{e_{m}}<\frac{2}{1.5152}e_{n}\leq2
\]
which is a contradiction. On the other hand, if $e_{h}\geq e_{m}$ and
$e_{j}\geq e_{m}$, then we have%
\begin{align*}
e_{m}+e_{m}  &  \leq e_{h}e_{k}+e_{j}e_{l}\\
&  =e_{m}e_{n}\cosh x\\
&  <\frac{2}{1.5152}e_{m}e_{n}%
\end{align*}
Dividing both sides by $e_{m}$ gives $2<\frac{2}{1.5152}e_{n}\leq2$, a
contradiction. Hence by symmetry no two of $e_{h}$, $e_{j}$, $e_{k}$,
and $e_{l}$ can be greater than or equal to $e_{m}$; therefore at
least three of them are strictly less than $e_{m}$. This proves the
lemma.\qed

\begin{lemma}
The $\mathcal{O}(1)$-edge does not intersects or pass within less than
$\delta_3$ of any $\mathcal{O}(n)$-edge where $e_n\le
1.5152$.\label{lem:o1_oN}
\end{lemma}

\noindent\emph{Proof:} Apply Lemma \ref{lem:one_cell_dist} with $m=1$
and $x<\delta_3$ to conclude that
either three of $\{o(h),o(j),o(k),o(l)\}$ are less than $o(1)=0$,
which is ridiculous, or else $e_n> 1.5152$, a contradiction.\qed

\begin{lemma}
If the $\mathcal{O}(2)$-edge intersects the $\mathcal{O}(3)$-edge, or
passes within less than $\delta_2$ of it, and if $e_3\le 1.5152$, then
there exists a geometric Mom-$2$ structure consisting only of
$(1,1,2)$-triples.\label{lem:o2_o3}
\end{lemma}

\noindent\emph{Proof:} Applying Lemma \ref{lem:one_cell_dist} we get
that $e_h=e_j=e_k=e_l=e_1$.  Hence $d(A,C)=d(A,D)=d(B,C)=d(B,D)=0$ and
therefore $\{A,B,C\}$ and $\{A,C,D\}$ are both $(1,1,2)$-triples. They
cannot be equivalent under the action of $\pi_1(N)$, because any
element of the group which maps one triple to the other would have to
map the pair $\{A,C\}$ to itself, and hence be either elliptic or the
identity, a contradiction. Therefore these two triples are distinct
and constitute a geometric Mom-$2$ structure, completing the
proof.\qed

\begin{lemma}
If $e_n\le 1.5152$ then the $\mc{O}(n)$-edge neither intersects itself
nor passes within $\delta_1$ of itself.\label{lem:oN_oN}
\end{lemma}

\noindent\emph{Proof:} This is a direct corollary of Lemma
\ref{lem:dist_line_self}.\qed\

\bigskip
The three preceding lemmas show that if $N$ contains a geometric
Mom-$2$ or Mom-$3$ structure of one of the types described in Theorem
\ref{thrm:comb_existence} then the $\mathcal{O}(n)$-edges of
the complex $\Delta$ are embedded in $N$ and do not intersect one
another, \emph{unless} the Mom-$n$ structure includes the
$\mathcal{O}(3)$-edge. In that case, either the $\mathcal{O}(n)$-edges
of the complex $\Delta$ are embedded and do not intersect or else
there exists \emph{another} geometric Mom-$n$ structure,
specifically a Mom-$2$ structure consisting only of $(1,1,2)$-triples.
In other words either the edges in question are
embedded and do not intersect, or else we can find a simpler Mom-$n$
structure; hence we may argue by induction that $N$ must possess a
geometric Mom-$n$ structure with embedded edges. This type of
induction argument will be repeated several times throughout this
section.

To prove Theorem \ref{thrm:comb_embedded}
we must also show that the $2$-cells of $\Delta$ corresponding to
$(n,m,k)$-triples are embedded and do not intersect one another.
Since the $2$-cells are simply connected and totally geodesic, and
since the $\mc{O}(n)$-edges in their boundary are geodesic arcs
perpendicular to the horospheres at their endpoints, it is
straightforward to show that if two $2$-cells intersect then a
$mc{O}(n)$-edge in the boundary of one $2$-cell must intersect
the other $2$-cell. The previous lemmas imply that this intersection
must occur in the interior of the $2$-cell. Therefore to complete the
proof of Theorem \ref{thrm:comb_embedded} it is sufficient to prove
the following:

\begin{proposition}
Suppose $N$ has a geometric Mom-$n$ structure of one of the types
described in Theorem \ref{thrm:comb_existence}. If any of the
$\mc{O}(n)$-edges of the complex $\Delta$ intersect any of the
$2$-cells in $\Delta$, then there must exist a simpler geometric
Mom-$n$ structure in $N$, which is also of one of the types described
in Theorem
\ref{thrm:comb_existence}.
\label{prop:edge_cell}
\end{proposition}

\noindent\emph{Proof: }The idea of the proof is as follows: for every
possible case, use the previous lemmas in this section to either
demonstrate a contradiction or reduce the problem to a strictly
previous case. But before enumerating the cases, consider the diagram
in
\begin{figure}[tb]
\begin{center}
\includegraphics{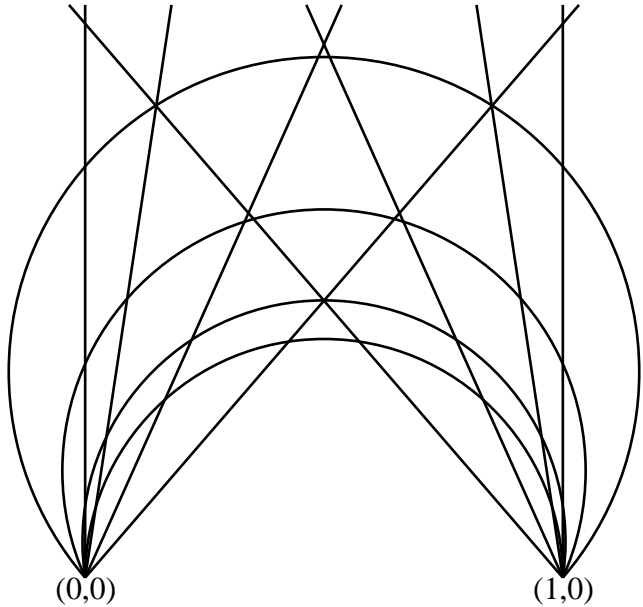}
\end{center}
\caption{An ideal triangle, together with three equidistant curves for
each side.}
\end{figure}
figure 5.

This diagram shows an ideal triangle in the upper half-space model of
$\BH^2$ with vertices at $0$, $1$, and $\infty$, together with nine
curves that are each equidistant from one of the three sides. For each
of the three sides of the triangle there are three curves: one at a
distance of $\delta_1$ from the side, one at a distance of $\delta_2$,
and one at a distance of $\delta_3$.
(To construct this diagram it is sufficient to note the following. A curve in
the upper half-space model which is at a constant distance $r$ from the line
from $0$ to $\infty$ is just a line passing through $0$ with slope $\pm\left(
\sinh r\right)  ^{-1}$; a curve which is equidistant from the line from $1$ to
$\infty$ is constructed similarly. A curve at constant distance $r$ from the
line from $0$ to $1$ is a circular arc passing from $0$ to $1$ through the
point $(\frac{1}{2},\frac{1}{2}e^{\pm r})$. See for example \cite{thu}.)

From the diagram, the following result is immediate:

\begin{lemma}
If $p$ is a point in the interior of the ideal triangle with vertices at $0$,
$1$, and $\infty$, and if $\lambda_1$, $\lambda_2$, $\lambda_3$ are the
sides of the triangle in any order, then:

\begin{enumerate}
\item  If $d(p,\lambda_1)\geq\delta_3$ and
$d(p,\lambda_2)\geq\delta_3$, then $d(p,\lambda_3)<\delta_1$.

\item  If $d(p,\lambda_1)\geq\delta_3$ and $d(p,\lambda_2)\geq
\delta_1$, then $d(p,\lambda_3)<\delta_3$.

\item (Corollary to the previous part) If
$d(p,\lambda_1)\geq\delta_3$, then either one of $d(p,\lambda_2)$,
$d(p,\lambda_3)$ is $<\delta_2$ or else both of $d(p,\lambda_2)$,
$d(p,\lambda_3)$ are $<\delta_3$.

\item If $d(p,\lambda_i)\ge\delta_2$ for all $i\in\{1,2,3\}$ then
$d(p,\lambda_j)<\delta_3$ for at least two $j\in\{1,2,3\}$.

\end{enumerate}\label{lem:triangle_dists}
\end{lemma}\qed\

\bigskip
The conclusion in the third part of the above lemma is annoyingly
weak. But note that the region in the
ideal triangle where $d(p,\lambda_1)\ge\delta_3$,
$d(p,\lambda_2)\ge\delta_2$, and $d(p,\lambda_3)\ge\delta_2$ is very
small; if $\delta_3$ is replaced with even a slightly smaller number
then a stronger conclusion would result. This turns out to be useful
enough that we do so now; the proof of the following lemma, while not
following immediately from the diagram, is elementary enough that we
omit it for brevity:

\begin{lemma}
If $d(p,\lambda_1)\ge 0.9$, then one of $d(p,\lambda_2)$,
$d(p,\lambda_3)$ is $<\delta_2$.
\label{lem:triangle_dists_2}
\end{lemma}\qed\

\bigskip
Now we begin enumerating the various cases to prove Proposition
\ref{prop:edge_cell}, according to which edge is involved and which
triple corresponds to the $2$-cell involved. In each case, let $p$
be the point where the edge intersects the $2$-cell, and note that the
two-cell is contained in an ideal triangle which is isometric to that used in
Lemma \ref{lem:triangle_dists}. We will refer to the sides of the
triangle as $\lambda_j$ for $j\in\{1,2,3\}$ and let
$d_j=d(p,\lambda_j)$. Also we will let $A_1$, $A_2$, and $A_3$ denote
the horoballs such that $\lambda_1$ goes from $A_1$ to $A_2$,
$\lambda_2$ goes from $A_2$ to $A_3$, and $\lambda_3$ goes from $A_3$
to $A_1$. Finally let $B_1$ and $B_2$ be the horoballs at the ends of
the $\mathcal{O}(n)$-edge which passes through $p$. See
\begin{figure}[tb]
\begin{center}
\includegraphics{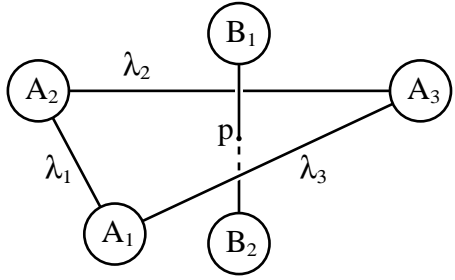}
\end{center}
\caption{The horoballs and edges used in the proof of Proposition
\ref{prop:edge_cell}.}
\end{figure}
figure 6.

\noindent\emph{Case 1:\ }$\mathcal{O}(1)$\emph{-edge, any triple.} Part 1 of
Lemma \ref{lem:triangle_dists} implies that $d_j<\delta_3$
for some $j$. Now apply Lemma \ref{lem:o1_oN} to obtain a contradiction.

\noindent\emph{Case 2:\ }$\mathcal{O}(2)$\emph{-edge, }$(1,1,2)$%
-\emph{triple.} Suppose that $\lambda_1$ and $\lambda_2$ contain
lifts of the $\mathcal{O}(1)$-edge while $\lambda_3$ contains a lift of the
$\mathcal{O}(2)$-edge. Part 1 of Lemma \ref{lem:triangle_dists} implies that
either $d_1<\delta_3$, $d_2<\delta_3$, or $d_3<\delta_1$. In the first
two cases Lemma \ref{lem:o1_oN} applies while in the third case Lemma
\ref{lem:oN_oN} applies. Either possibility leads to a contradiction.

\noindent\emph{Case 3:\ }$\mathcal{O}(2)$\emph{-edge, }$(1,2,2)$%
-\emph{triple.} Suppose that $\lambda_1$ contains a lift of the
$\mathcal{O}(1)$-edge while $\lambda_2$ and $\lambda_3$ contain lifts
of the $\mathcal{O}(2)$-edge. If $d_1<\delta_3$ then Lemma
\ref{lem:o1_oN} applies and we're done. If either $d_2$ or $d_3$ is
less than $\delta_2$ then the same argument used in the proof of
Lemma \ref{lem:o2_o3} proves the existence of a geometric Mom-$2$
structure involving only $(1,1,2)$-triples; either this new Mom-$2$ is
embedded or else we reduce to case 2.

By part 3 of Lemma \ref{lem:triangle_dists} the only remaining
possibility is that $d_1\ge\delta_3$, $d_2<\delta_3$, and
$d_3<\delta_3$. Since $d_2<\delta_3$, by Lemma
\ref{lem:one_cell_dist}, at least three of the four horoball pairs
$\{A_2,B_1\}$, $\{B_1,A_3\}$, $\{A_3,B_2\}$, $\{B_2,A_2\}$ are
elements of $\mc{O}(1)$. Similarly, since $d_3<\delta_3$ at least
three of the four pairs $\{A_3,B_1\}$, $\{B_1,A_1\}$, $\{A_1,B_2\}$,
and $\{B_2,A_3\}$ lie in $\mc{O}(1)$. But $\{A_1,A_2\}$ is already in
$\mc{O}(1)$; hence by Corollary \ref{cor:no_bad_triples} at least one
of $\{A_1,B_1\}$, $\{A_2,B_1\}$ must not lie in $\mc{O}(1)$, and
similarly one of $\{A_1,B_2\}$, $\{A_2,B_2\}$ must not lie in
$\mc{O}(1)$. Up to symmetry, we must have $\{A_1,B_1\}$,
$\{A_2,B_2\}\in\mc{O}(1)$, and $\{A_2,B_1\}$,
$\{A_1,B_2\}\not\in\mc{O}(1)$. Hence $d(A_2,B_1)\ge o(2)$ and
$d(A_1,B_2)\ge o(2)$. Now apply Corollary \ref{cor:two_line_dist} to
$A_1$, $A_2$, $B_1$, and $B_2$ to obtain
\begin{eqnarray*}
{e_2}^2+1 & \le & e_2 \cosh d_1 \\
\Rightarrow\ d_1 & \ge & \cosh^{-1}(e_2+{e_2}^{-1}) \\
& \ge & \cosh^{-1} 2
\end{eqnarray*}

Note that $\cosh^{-1} 2=1.3169\ldots$ so Lemma
\ref{lem:triangle_dists_2} applies. Hence one of $d_2$, $d_3$ is less
than $\delta_2$, and a preceding argument applies. This completes the
proof in this case.

We can save a bit of time at this point by noting that the argument in case 3
still works if we replace the $(1,2,2)$-triple with a $(1,2,3)$-triple, a
$(1,3,3)$-triple, or (with minor modifications) a $(1,1,3)$-triple;
henceforth we will assume that case 3 encompasses all of these
possibilities.

\noindent\emph{Case 4: }$\mathcal{O}(3)$\emph{-edge, }$(1,1,2)$%
-\emph{triple.} Suppose that $\lambda_1$ and $\lambda_2$ contain lifts
 of the $\mc{O}(2)$-edge while $\lambda_3$ contains a lift of the
 $\mc{O}(2)$-edge. If $d_1<\delta_3$ or $d_2<\delta_3$ then Lemma
 \ref{lem:o1_oN} applies, giving a contradiction. Otherwise by part 1 of Lemma
 \ref{lem:triangle_dists} we have $d_3<\delta_1$. Since
 $\delta_1<\delta_2$, we can use the argument from Lemma \ref{lem:o2_o3}
 to find a geometric Mom-$2$ structure using only $(1,1,2)$ triples,
 reducing the problem to case 2.

\noindent\emph{Case 5: }$\mathcal{O}(3)$\emph{-edge,
 }$(1,1,3)$-\emph{triple.} This case is nearly identical to the
 previous one, except we use Lemma \ref{lem:oN_oN} to obtain a
 contradiction instead of using Lemma
 \ref{lem:o2_o3} to reduce to case 2.

\noindent\emph{Case 6: }$\mathcal{O}(3)$\emph{-edge,
 }$(1,2,2)$-\emph{triple.} Use the argument in case 3 with trivial
 modifications.

\noindent\emph{Case 7:\ }$\mathcal{O}(2)$\emph{-edge, }$(2,2,3)$%
-\emph{triple.} Suppose that $\lambda_1$ and $\lambda_2$ contain
lifts of the $\mathcal{O}(2)$-edge while $\lambda_3$ contains a lift of the
$\mathcal{O}(3)$-edge. If $d_i<\delta_2$ for any $i\in\{1,2,3\}$ then
we can use the argument from Lemma \ref{lem:o2_o3} as before to
construct a geometric Mom-$2$ structure with just $(1,1,2)$-triples and
reduce to case 2. So suppose $d_i\ge\delta_2$ for all
$i\in\{1,2,3\}$; then by part 4 of Lemma \ref{lem:triangle_dists} we
have $d_j<\delta_3$ for at least two $j\in\{1,2,3\}$. Suppose that
$d_1<\delta_3$ and $d_2<\delta_3$; the other possibilities can be
handled similarly. Then by Lemma \ref{lem:one_cell_dist}
at least three of the pairs $\{A_1,B_1\}$,
$\{B_1,A_2\}$, $\{A_2,B_2\}$, and $\{B_2,A_1\}$ lie in $\mc{O}(1)$,
and at least three of the pairs $\{A_2,B_1\}$, $\{B_1,A_3\}$,
$\{A_3,B_2\}$, and $\{B_2,A_2\}$ lie in $\mc{O}(1)$.

Suppose all three of the pairs $\{A_1,B_i\}$, $\{A_2,B_i\}$, and
$\{A_3,B_i\}$ were in $\mc{O}(1)$ for some $i\in\{1,2\}$. Then both of
the triples
$\{A_1,A_2,B_i\}$ and $\{A_2,A_3,B_i\}$ are $(1,1,2)$-triples. If
there is an isometry $g\in\pi_1(N)$ taking one triple to another, then
that isometry must be a parabolic isometry fixing $B_i$ and taking
$A_1$ to $A_2$ and $A_2$ to $A_3$. But this would imply that the arc
from $A_1$ to $A_3$ intersects either $B_i$, $A_2$, or the arc from
$B_i$ to $A_2$, contradicting our assumptions to this point. If there
is no such isometry $g$, then we have a geometric Mom-$2$
structure with two
distinct $(1,1,2)$-triples and can reduce to case 2.

So suppose that at least one of the pairs $\{A_1,B_i\}$, $\{A_2,B_i\}$,
and $\{A_3,B_i\}$ does not lie in $\mc{O}(1)$, both for $i=1$ and for
$i=2$. Up to symmetry, the only way this can happen given our earlier
assumptions is if the pairs $\{A_1,B_1\}$, $\{B_1,A_2\}$,
$\{A_2,B_2\}$, and $\{B_2,A_3\}$ lie in $\mc{O}(1)$ while the pairs
$\{A_1,B_2\}$ and $\{A_3,B_1\}$ do not. If $\{A_1,B_2\}$ lies in
$\mc{O}(2)$, then the triples $\{A_1,A_2,B_2\}$ and $\{B_1,B_2,A_2\}$
form a geometric Mom-$2$ structure with one $(1,1,2)$-triple and one
$(1,2,2)$-triple, allowing us to reduce the problem to either case 2
or case 3. So suppose $d(A_1,B_2)\ge o(3)$. For similar reasons we can
suppose that $d(A_3,B_1)\ge o(3)$. Then applying Corollary
\ref{cor:two_line_dist} to $A_1$, $A_3$, $B_1$, and $B_2$ we obtain
\begin{eqnarray*}
{e_3}^2+1 & \le & e_2 e_3\cosh d_3 \\
          & \le & {e_3}^2\cosh d_3 \\
\Rightarrow\ d_3 & \ge & \cosh^{-1}(1+{e_3}^{-2})
\end{eqnarray*}
Note that if $e_3\le 1.5152$ then this implies that $d_3 \ge 0.9$ and
hence Lemma \ref{lem:triangle_dists_2} applies. Therefore at least one
of $d_1$, $d_2$ is less than $\delta_2$, and a preceding argument
applies. This completes the proof in this case.

\noindent\emph{Case 8:\ }$\mathcal{O}(2)$\emph{-edge,
}$(2,3,3)$-\emph{triple.} An almost identical argument to the one in
case 7 applies, except that the problem may reduce to case 5 instead
of case 2.

\noindent\emph{Case 9: }$\mc{O}(3)$\emph{-edge,
}$(2,2,3)$\emph{-triple}. Suppose that $\lambda_1$ and $\lambda_2$
contain lifts of the $\mc{O}(2)$-edge while $\lambda_3$ contains a
lift of the $\mc{O}(3)$-edge. If $d_3<\delta_1$ then Lemma
\ref{lem:oN_oN} provides a contradiction, while if $d_i<\delta_2$ for
$i=1$ or $2$ then the argument from Lemma \ref{lem:o2_o3} will produce
a simpler geometric Mom-$2$ structure. Hence suppose that $d_3\ge\delta_1$
and $d_i\ge\delta_2>\delta_1$ for $i=1$ and $2$. Part 2 of Lemma
\ref{lem:triangle_dists} then implies that $d_j<\delta_3$ for at least
two different $j\in\{1,2,3\}$. 

Suppose $d_1<\delta_3$ and $d_3<\delta_3$. Since $d_1<\delta_3$, by
Lemma \ref{lem:one_cell_dist} at least three of $\{A_1,B_1\}$,
$\{B_1,A_2\}$, $\{A_2,B_2\}$, $\{B_2,A_1\}$ lie in $\mc{O}(1)$. Wolog
assume that $\{A_1,B_1\}$ and $\{B_1,A_2\}$ lie in $\mc{O}(1)$. If
$\{B_1,A_3\}$ lies in $\mc{O}(1)$ then the triples $\{A_1,B_1,A_2\}$
and $\{A_2,B_1,A_3\}$ are both $(1,1,2)$-triples, and form a simpler
geometric Mom-$2$ structure unless they are equivalent by the action of
parabolic element $g\in\pi_1(N)$ which fixes $B_1$; but this would
imply that the arc from $B_1$ to $A_3$ intersects the arc from $A_1$ to
$A_2$, a simpler case. If $\{B_1,A_3\}$ lies in $\mc{O}(2)$ then
$\{A_1,B_1,A_2\}$ is a $(1,1,2)$-triple while $\{A_2,B_1,A_3\}$ is a
$(1,2,2)$-triple and we have a simpler Mom-$2$. So suppose
$\{B_1,A_3\}\not\in\mc{O}(1)\cup\mc{O}(2)$. Then since $d_3<\delta_3$,
by Lemma \ref{lem:one_cell_dist} the pairs $\{B_2,A_1\}$ and
$\{B_2,A_3\}$ must lie in $\mc{O}(1)\cup\mc{O}(2)$. Now consider
$\{B_2,A_2\}$; if this pair lies in $\mc{O}(1)\cup\mc{O}(2)$ then the
triples $\{A_1,A_2,B_1\}$ and $\{A_1,A_2,B_2\}$ will form a simpler
geometric Mom-$2$ structure, so suppose
$\{B_2,A_2\}\not\in\mc{O}(1)\cup\mc{O}(2)$. We have reached a point
where $d(B_1,A_3)\ge o(3)$ and $d(B_2,A_2)\ge o(3)$, so applying Corollary
\ref{cor:two_line_dist} to the four horoballs $A_2$, $A_3$, $B_1$,
$B_3$ we get
\begin{eqnarray*}
{e_3}^2+1 & \le & e_2 e_3 \cosh d_2 \\
& \le & {e_3}^2 \cosh d_2 \\
\Rightarrow\ d_2 & \ge & \cosh^{-1}(1+{e_3}^{-2})
\end{eqnarray*}
Note that if $e_3\le 1.5152$ then this implies that $d_2 \ge 0.9$ and
hence Lemma \ref{lem:triangle_dists_2} implies that at least one of
$d_1$, $d_3$ is less than $\delta_2$. But if $d_1<\delta_2$ then by
Corollary \ref{cor:two_line_dist} we would have $\{B_2,A_2\}\in\mc{O}(1)$,
while if $d_3<\delta_2$ then $\{B_1,A_3\}\in\mc{O}(1)\cup\mc{O}(2)$;
both conclusions contradict our assumptions to this point.

The argument when $d_2<\delta_3$ and $d_3<\delta_3$ is identical to
this one by symmetry, and the the argument when $d_1<\delta_3$ and
$d_2<\delta_3$ is similar.

\noindent\emph{Case 10:\ }$\mathcal{O}(3)$\emph{-edge, }$(1,2,3)$%
-\emph{triple.} (This is by far the trickiest case.) Suppose that
 $\lambda_1$ contains a lift of the $\mc{O}(1)$-edge, $\lambda_2$
 contains a lift of the $\mc{O}(2)$-edge, and $\lambda_3$ contains a
 lift of the $\mc{O}(3)$-edge. If $d_1<\delta_3$, then Lemma
 \ref{lem:o1_oN} leads to a contradiction. If $d_2<\delta_2$, then as
 before we can repeat the argument of Lemma \ref{lem:o2_o3} to obtain
 a simpler Mom-$2$ and reduce to case 2. While if $d_3<\delta_1$ then
 Lemma \ref{lem:oN_oN} provides a contradiction.

So suppose that $d_1\ge\delta_3$, $d_2\ge\delta_2$, and
$d_3\ge\delta_1$. Note that by part 2 of Lemma
\ref{lem:triangle_dists} this implies that $d_2<\delta_3$ and
$d_3<\delta_3$. By Lemma \ref{lem:one_cell_dist}, at least three of
the four pairs $\{A_2,B_1\}$, $\{B_1,A_3\}$, $\{A_3,B_2\}$, and
$\{B_2,A_1\}$ lie in $\mc{O}(1)$, while at least three of the four
pairs $\{A_1,B_1\}$, $\{B_1,A_3\}$, $\{A_3,B_2\}$, $\{B_2,A_1\}$ lie
in $\mc{O}(1)\cup\mc{O}(2)$. Unfortunately this is not quite enough
information to construct a simpler Mom-$2$ or Mom-$3$, so we must dig
deeper.

Let $\sigma$ be the $2$-cell spanning $A_1$, $A_2$,
$A_3$, and the arcs between them. Since $\{A_3,A_1\}$ and
$\{B_1,B_2\}$ are both elements of the orthopair class $\mc{O}(3)$,
there exists a group element $g\in\pi_1(N)$ which sends $\{A_3,A_1\}$
to $\{B_1,B_2\}$. Furthermore since $\sigma$ is totally geodesic,
$g(\sigma)\cap\sigma$ must contain a geodesic line segment with one
endpoint $p$ in the arc from $B_1$ to $B_2$. Let $q$ be the other
endpoint of this line segment. What are the possible locations of $q$?

If $q$ lies in the interior of $\sigma$ then either $g(A_j)$
intersects $\sigma$ for some $j$, in which case we can apply Lemma
\ref{lem:ball_intersect}, or else either $g(\lambda_1)$ or
$g(\lambda_2)$ intersect $\sigma$, in which case the problem reduces
to case 1 or case 3 respectively. So suppose $q$ lies on the boundary
of $\sigma$. If $q$ lies in the interior of $\lambda_1$, then
$g^{-1}(\lambda_1)$ intersects $\sigma$ and we can reduce to case
1. Similarly if $q$ lies in the interior of $\lambda_2$ then we can
reduce to case 3.

If $q$ lies in $\sigma\cap A_2$, then $g$ must be a
parabolic element of $\pi_1(N)$ fixing $A_2$. But note that the
projections of $\lambda_3$ and $g(\lambda_3)$ onto the surface of
$A_2$ clearly intersect; if $g$ is parabolic fixing $A_2$, then the
only way this can happen is if $\lambda_3$ and $g(\lambda_3)$
themselves intersect. This contradicts Lemma \ref{lem:oN_oN}.

If $q$ lies in $\sigma\cap A_3$, then $g$ must send the triple
$\{A_1,A_2,A_3\}$ to the triple $\{A_3,B_1,B_2\}$; in particular
$\{A_3,B_1,B_2\}$ is a $(1,2,3)$-triple. Wolog, assume that
$\{A_3,B_2\}$ lies in $\mc{O}(2)$ and $\{A_3,B_1\}$ lies in
$\mc{O}(1)$. Since $d_2\le\delta_3$, this implies that
$\{A_2,B_1\}$ and $\{A_2,B_2\}$ lie in $\mc{O}(1)$ as well. Thus
$\{A_2,B_1,A_3\}$ is a $(1,1,2)$-triple while $\{A_2,B_2,A_3\}$ is a
$(1,2,2)$-triple, forming a simpler geometric Mom-$2$ structure and
reducing the problem to either case 1, 2, or 3.

There are two remaining possibilities. Suppose $q$ lies in
$\sigma\cap A_1$; then $g$ must send $\{A_1,A_2,A_3\}$ to
$\{A_1,B_1,B_2\}$. Note in particular that we must have
$g(A_2)=A_1$. Wolog assume that $\{A_1,B_1\}$ lies in $\mc{O}(1)$
while $\{A_1,B_2\}$ lies in $\mc{O}(2)$. Applying Lemma
\ref{lem:no_bad_triples} to the triple $\{A_1,A_2,B_1\}$ we see that
$\{A_2,B_1\}$ cannot be in $\mc{O}(1)$. Since $d_2<\delta_3$, we must
have all of $\{B_1,A_3\}$, $\{A_3,B_2\}$, $\{B_2,A_2\}$ in $\mc{O}(1)$
instead. So $\{A_1,B_1,A_3\}$ is a $(1,1,3)$-triple, while
$\{A_1,B_2,A_3\}$ is a $(1,2,3)$-triple. Now consider
$g^{-1}(\sigma)$. Comparing $\sigma$ and $g(\sigma)$, we see that
$g^{-1}(\sigma)$ must intersect $\lambda_3$, and must be bounded at
one corner by $A_2=g^{-1}(A_1)$. Specifically
$g^{-1}(\lambda_1)\in\mc{O}(1)$ and $g^{-1}(\lambda_3)\in\mc{O}(3)$
must have one endpoint on $A_2$.
If any of the four arcs spanning the horoball pairs $\{A_1,B_1\}$,
$\{B_1,A_3\}$, $\{A_3,B_2\}$, or $\{B_2,A_1\}$ (all of which are in
$\mc{O}(1)\cup\mc{O}(2)$) intersect $g^{-1}(\sigma)$ then we can
reduce to a previous case, so suppose this doesn't happen.
We still must have $g^{-1}(\sigma)$ intersecting $\lambda_3$. If
either $g^{-1}(\lambda_1)$ or $g^{-1}(\lambda_2)$ intersect either of
the two-cells spanning the triples $\{A_1,B_1,A_3\}$ or
$\{A_1,B_2,A_3\}$, then we can still reduce the problem to a previous
case. The only way $g^{-1}(\sigma)$ can intersect $\lambda_3$ without
such an intersection occurring (and without $g^{-1}(\sigma)$
intersecting the interior of a horoball, which contradicts Lemma
\ref{lem:ball_intersect}) is if $g^{-1}(A_1)=B_2$ and
$g^{-1}(\lambda_3)$ intersects the two-cell spanning
$\{A_1,B_1,A_3\}$. But this merely lets us reduce to case 5 instead.

Thus we come to the final possibility: suppose that $q$ lies in the
interior of $\lambda_3$. Then the image of the arc from $p$ to $q$
under $g^{-1}$ must be another geodesic line segment going from
$g^{-1}(p)$ on $\lambda_3$ to $g^{-1}(q)$ in the interior of
$\sigma$. Note that $p\not =g^{-1}(q)$, otherwise $g$ would be
elliptic of order $2$. Since $g^{-1}(q)$ is the point where $\sigma$
intersects the arc from $g^{-1}(A_1)$ to $g^{-1}(A_3)$, and since this
arc is also a lift of the $\mc{O}(3)$-edge, the same arguments that
apply to $p$ also apply to $g^{-1}(q)$. Specifically: if
$d(g^{-1}(q),\lambda_1)<\delta_3$, $d(g^{-1}(q),\lambda_2)<\delta_2$,
or $d(g^{-1}(q),\lambda_3)<\delta_1$, then we can apply Lemma
\ref{lem:o1_oN}, \ref{lem:o2_o3}, or \ref{lem:oN_oN} respectively to
get either a contradiction or a reduction to a simpler case. So assume
that none of those three inequalities hold. Note that this implies
that $d(g^{-1}(q),\lambda_3)<\delta_3$ just as our previous
assumptions implied that $d_3=d(p,\lambda_3)<\delta_3$.

Now consider the pairs $\{B_1,B_2\}$ and
$\{g^{-1}(A_1),g^{-1}(A_3)\}$. These cannot be the same pair; if they
were, then $g^2$ would fix the pair $\{B_1,B_2\}$ and hence be either
elliptic or the identity, a contradiction. So the two pairs are either
completely disjoint or else intersect in a single element. Suppose
they intersect in a single element; wolog, assume in particular that
$B_1=g^{-1}(A_1)$. Note that we can't have $g(A_1)=B_1$; if we did
then $g^2$ would fix $A_1$ and hence so would $g$, a contradiction. So
$g(A_1)=B_2$, and $g(A_3)=B_1$. To summarize we have $g(B_1)=A_1$ and
$g^2(B_1)=B_2$, while at the same time we must have
$g^{-1}(B_1)=A_3$. Arrange the upper half-space model of
$\BH^3$ so that $B_1$ is centered at the point at infinity,
$A_1=g(B_1)$ and $A_3=g^{-1}(B_1)$ are centered on the real line, and
$A_2$ is centered at some point with positive imaginary part. By
assumption, the arcs from $B_1$ to $B_2$ and from $B_1$ to
$g^{-1}(A_3)$ both intersect the interior of the two-cell spanning the
triple $\{A_1,A_2,A_3\}$; hence the centers of $B_2=g^2(B_1)$ and
$g^{-1}(A_3)=g^{-2}(B_1)$ also have positive imaginary part. But this
is impossible: by direct calculation in $\mathrm{PSL}(2,\BC)$, if $g$
is an element such that $g(\infty)$ and $g^{-1}(\infty)$ are both real
then either $g^2(\infty)$ and $g^{-2}(\infty)$ are also both real or
else their imaginary parts have opposite signs. Hence, the pairs
$\{B_1,B_2\}$ and $\{g^{-1}(A_1),g^{-1}(A_3)\}$ must be completely
disjoint, or in other words the arcs $g(\lambda_3)$ and
$g^{-1}(\lambda_3)$ do not share a horoball at their endpoints.

Now suppose that $d(\lambda_3,g(\lambda_3))\ge\delta_2$ and
$d(\lambda_3,g^{-1}(\lambda_3))\ge\delta_2$. Examining figure 5 (and
remembering our previous assumptions), this implies that both $p$ and
$g^{-1}(q)$ must lie in the small region bounded by the curves
$d(\cdot ,\lambda_1)=\delta_3$, $d(\cdot ,\lambda_2)=\delta_2$, and
$d(\cdot ,\lambda_3)=\delta_2$. The diameter of this region can be
readily computed (it's not a triangle but it is contained inside one),
and it is far less than $\delta_1$. Hence
$d(g(\lambda_3),g^{-1}(\lambda_3))<\delta_1$, and hence Lemma
\ref{lem:oN_oN} applies, producing a contradiction. So therefore one
of $d(\lambda_3,g(\lambda_3))$ and $d(\lambda_3,g^{-1}(\lambda_3))$
must be less than $\delta_2$, and clearly if one is then both are by
isometry. Furthermore, the diameter of the region in figure 5 bounded
by the curves $d(\cdot ,\lambda_1)=\delta_3$, $d(\cdot
,\lambda_2)=\delta_2$, and $d(\cdot ,\lambda_3)=\delta_1$ can also be
computed to be less than $\delta_3$. Then by applying Lemma
\ref{lem:one_cell_dist} repeatedly, we get the following:
\begin{itemize}
\item All four of the pairs $\{A_1,B_1\}$, $\{B_1,A_3\}$,
$\{A_3,B_2\}$, and $\{B_2,A_1\}$ lie in $\mc{O}(1)\cup\mc{O}(2)$.
\item All four of the pairs $\{A_1,g^{-1}(A_1)\}$,
$\{g^{-1}(A_1),A_3\}$, $\{A_3,g^{-1}(A_3)\}$, and
$\{g^{-1}(A_3),A_1\}$ lie in $\mc{O}(1)\cup\mc{O}(2)$.
\item At least three of the pairs $\{B_1,g^{-1}(A_1)\}$,
$\{g^{-1}(A_1),B_2\}$, $\{B_2,g^{-1}(A_3)\}$, and
$\{g^{-1}(A_3),B_1\}$ lie in $\mc{O}(1)\cup\mc{O}(2)$.
\end{itemize}
So the six horoballs $A_1$, $A_3$, $B_1$, $B_2$, $g^{-1}(A_1)$, and
$g^{-1}(A_3)$ form a picture like the one in
\begin{figure}[tb]
\begin{center}
\includegraphics{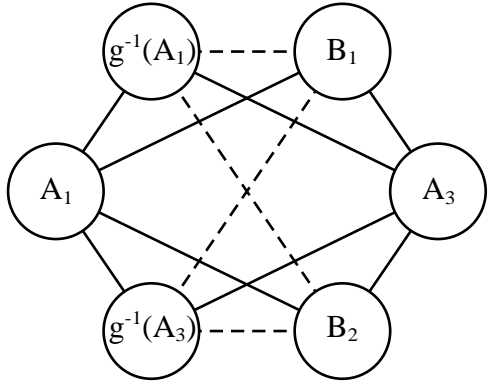}
\end{center}
\caption{A possible collection of horoballs used in case 10; the solid
lines (and three of the four dashed lines) indicate pairs which lie in
$\mc{O}(1)\cup\mc{O}(2)$.}
\end{figure}
figure 7, where each of solid edges and three of the four dashed edges is
a lift of either the $\mc{O}(1)$-edge or the
$\mc{O}(2)$-edge. Suppose, for sake of example, that the dashed edge
from $B_1$ to $g^{-1}(A_1)$ is the only edge which is not a lift of
either the $\mc{O}(1)$-edge or the $\mc{O}(2)$-edge. Then consider the
triple $\{A_1,B_2,g^{-1}(A_3)\}$. This triple must be an
$(a,b,b)$-triple where $a$ and $b$ are $1$ and $2$ in some order
(i.e., it's either a $(1,1,2)$-triple or a $(1,2,2)$-triple). But
this triple also has the property that each edge of the triple is
shared by another triple in the diagram which must also be either a
$(1,1,2)$-triple or a $(1,2,2)$-triple. It is straightforward to see
that some pair of these four triples, specifically the pair which
shares the edge which is a lift of the $\mc{O}(a)$-edge, must be
distinct under the action of $\pi_1(N)$ and hence form a geometric
Mom-$2$ structure. A similar argument holds for each of the other
dashed edges in the diagram.

This at last completes the proof in this case.

\noindent\emph{Case 11: }$\mc{O}(3)$\emph{-edge,
}$(1,3,3)$\emph{-triple}. Suppose $\lambda_3$ contains a lift of the
$\mc{O}(1)$-edge while $\lambda_1$ and $\lambda_2$ contain lifts of
the $\mc{O}(3)$-edge. If $d_3<\delta_3$ then Lemma \ref{lem:o1_oN}
provides a contradiction, while if $d_i<\delta_1$ for $i=1$ or $2$
then Lemma \ref{lem:oN_oN} provides a contradiction. So suppose
$d_1\ge\delta_1$, $d_2\ge\delta_1$, and $d_3\ge\delta_3$. Lemma
\ref{lem:triangle_dists} implies that $d_1$ and $d_2$ are both less
than $\delta_3$.

Suppose that in fact $d_1<\delta_2$. By Lemma \ref{lem:one_cell_dist},
all four pairs $\{A_1,B_1\}$, $\{B_1,A_2\}$, $\{A_2,B_2\}$,
$\{B_2,A_1\}$ and at least three of the four pairs $\{A_2,B_1\}$,
$\{B_1,A_3\}$, $\{A_3,B_2\}$, $\{B_2,A_2\}$ lie in
$\mc{O}(1)\cup\mc{O}(2)$. Wolog, assume that
$\{B_1,A_3\}$ $\in\mc{O}(1)\cup\mc{O}(2)$. We wish to construct a
geometric Mom-$2$ or Mom-$3$ structure that only uses triples
considered in
previous cases. However we have to be careful to ensure that we do not
construct a Mom-$3$ which is not torus-friendly; this can be ensured by not
selecting any
geometric Mom-$3$ structure that contains exactly two $(1,2,3)$-triples. To
start with, note that the triple $\{A_2,B_1,A_3\}$ is of type
$(1,1,2)$ or $(1,2,2)$. Now consider the four triples
$\{A_1,B_1,A_2\}$, $\{A_1,B_2,A_2\}$, $\{B_1,A_1,B_2\}$, and
$\{B_1,A_2,B_2\}$. Each of these triples is of type $(1,2,3)$,
$(1,1,3)$, or $(2,2,3)$.  If all four of these triples are of type
$(1,2,3)$, then no two triples can be equivalent under the action of
$\pi_1(N)$ because each two of those triples share a common
``edge''. So any three of those triples will form a geometric
Mom-$3$ structure which is torus-friendly. If between one and three of
these triples
are of type $(1,2,3)$, choose one triple of type $(1,2,3)$, one triple
not of type $(1,2,3)$, and $\{A_2,B_1,A_3\}$ to get a simpler
geometric Mom-$3$ structure which is torus-friendly. And if none of the four
triples is of type $(1,2,3)$, then all four must be of type $(a,a,3)$
for a fixed $a\in\{1,2\}$. Pick two such triples, say
$\{A_1,B_1,A_2\}$ and $\{A_1,B_2,A_2\}$; these triples can't be
equivalent under the action of $\pi_1(N)$ since that group has no
elliptic elements; hence they form a simpler geometric Mom-$2$ structure.

So instead suppose that $d_1\ge\delta_2$ and for similar reasons
suppose that $d_2\ge\delta_2$. Suppose now that
$\{A_2,B_1\}\not\in\mc{O}(1)\cup\mc{O}(2)$. Since $d_1$ and $d_2$ are
both less than $\delta_3$, Lemma \ref{lem:one_cell_dist} implies that
all of $\{A_1,B_1\}$, $\{A_1,B_2\}$, $\{A_3,B_1\}$, $\{A_3,B_2\}$, and
$\{A_2,B_2\}$ are in $\mc{O}(1)\cup\mc{O}(2)$. Then the triples
$\{A_1,B_1,A_3\}$ and $\{A_1,B_2,A_3\}$ are each of type $(1,1,2)$ or
$(1,2,2)$, and the triples $\{A_1,B_2,A_2\}$ and $\{A_2,B_2,A_3\}$ are
each of type $(1,1,3)$, $(1,2,3)$, or $(2,2,3)$. If 
$\{A_1,B_2,A_2\}$ and $\{A_2,B_2,A_3\}$ are equivalent due to the
action of $g\in\pi_1(N)$, then $g$ must be parabolic fixing $B_2$, and
hence $\lambda_3$ must intersect the arc from $B_2$ to $A_2$, a
previous case. So $\{A_1,B_2,A_2\}$, $\{A_2,B_2,A_3\}$, and
$\{A_1,B_2,A_3\}$ are all distinct triples and form a simpler Mom-$3$
which furthermore is torus-friendly, \emph{unless} $\{A_1,B_2,A_2\}$
and $\{A_2,B_2,A_3\}$ are both of type $(1,2,3)$. That is only
possible if $\{A_1,B_2\}$ and $\{B_2,A_3\}$ are both in $\mc{O}(2)$
(the other possibility, that they are both in $\mc{O}(1)$, makes
$\{A_1,B_2,A_3\}$ a $(1,1,1)$-triple, which is impossible). In this
case $\{A_1,B_2,A_3\}$ is of type $(1,2,2)$ and cannot be equivalent
to $\{A_1,B_1,A_3\}$ under the action of $\pi_1(N)$ since a group
element sending one triple to another would have to fix the pair
$\{A_1,A_3\}$; thus we have a simpler geometric Mom-$2$ structure.

So we may assume that $\{A_2,B_1\}\in\mc{O}(1)\cup\mc{O}(2)$, and
by symmetry we may assume that
$\{A_2,B_2\}\in\mc{O}(1)\cup\mc{O}(2)$. If the same holds for both
$\{A_1,B_1\}$ and $\{A_1,B_2\}$ then we may proceed exactly as if
$d_1<\delta_2$. So assume that one of those two pairs is not in
$\mc{O}(1)\cup\mc{O}(2)$, and similarly assume that one of the pairs
$\{A_3,B_1\}$, $\{A_3,B_2\}$ is not in
$\mc{O}(1)\cup\mc{O}(2)$. Now applying Corollary \ref{cor:two_line_dist}
to the four horoballs $A_1$, $A_3$, $B_1$, and $B_3$, we get either
\[
{e_3}^2+1 \le e_3 \cosh d_3
\]
or
\[
2e_3 \le e_3 \cosh d_3
\]
In either case, $d_3\ge\cosh^{-1}2>0.9$, and hence by Lemma
\ref{lem:triangle_dists_2} one of $d_1$, $d_2$ must be less than
$\delta_2$, contradicting our assumptions to this point and completing
this case.

\noindent\emph{Case 12: }$\mc{O}(3)$\emph{-edge,
}$(2,3,3)$\emph{-triple}. Suppose $\lambda_1$ and $\lambda_2$ contain
lifts of the $\mc{O}(3)$-edge while $\lambda_3$ contains a lift of the
$\mc{O}(2)$-edge.

Note that if $d_1<\delta_1$ or $d_2<\delta_1$ then we get a
contradiction from Lemma \ref{lem:oN_oN}, and if $d_3<\delta_2$ then
we can produce a simpler geometric Mom-$2$ structure just as in Lemma
\ref{lem:o2_o3}. So assume $d_1\ge\delta_1$, $d_2\ge\delta_1$, and
$d_3\ge\delta_2>\delta_1$. By Lemma \ref{lem:triangle_dists} this implies
that at least two of $d_1$, $d_2$, and $d_3$ must be less than
$\delta_3$, and by Lemma \ref{lem:one_cell_dist} that in turn implies
that at least two of the following statements are true:
\begin{itemize}
\item At least three of $\{A_1,B_1\}$, $\{B_1,A_2\}$, $\{A_2,B_2\}$,
$\{B_2,A_1\}$ are in $\mc{O}(1)\cup\mc{O}(2)$.
\item At least three of $\{A_2,B_1\}$, $\{B_1,A_3\}$, $\{A_3,B_2\}$,
$\{B_2,A_2\}$ are in $\mc{O}(1)\cup\mc{O}(2)$.
\item At least three of $\{A_3,B_1\}$, $\{B_1,A_1\}$, $\{A_1,B_2\}$,
$\{B_2,A_3\}$ are in $\mc{O}(1)\cup\mc{O}(2)$.
\end{itemize}

Suppose for a moment that \emph{all four} of $\{A_1,B_1\}$,
$\{B_1,A_3\}$, $\{A_3,B_2\}$, $\{B_2,A_1\}$ are in
$\mc{O}(1)\cup\mc{O}(2)$. I.e., suppose that $\{A_1,B_1,A_3\}$ and
$\{A_1,B_2,A_3\}$ are both of type $(1,1,2)$ or $(1,2,2)$. If these
two triples are not equivalent under $\pi_1(N)$ then they form a
simpler geometric Mom-$2$ structure. If they are equivalent due to
$g\in\pi_1(N)$, then either $g$ is elliptic and fixes $\{A_1,A_3\}$ (a
contradiction), $g$ is parabolic fixing one of $A_1$ or $A_3$ (in
which case the one-cell from $B_1$ to $B_2$ will intersect another
one-cell, reducing the problem to a previous case), or else up to
symmetry we may assume that $\{A_1,B_1\}$, $\{A_3,B_2\}$ are in
$\mc{O}(1)$ and $\{A_1,B_2\}$, $\{A_3,B_1\}$ are in $\mc{O}(2)$. In
the latter case, assume wolog that
$\{B_1,A_2\}\in\mc{O}(1)\cup\mc{O}(2)$ (we know this must be true for
one of $\{B_1,A_2\}$, $\{B_2,A_2\}$). Then one of the two triples
$\{A_1,A_2,B_1\}$ and $\{A_2,A_3,B_1\}$ must be of type $(2,2,3)$ or
$(1,1,3)$; this triple along with $\{A_1,B_1,A_2\}$ and
$\{B_1,B_2,A_1\}$ form a simpler geometric Mom-$3$ structure where each
triple is of a different type (and hence none are equivalent under
$\pi_1(N)$) and exactly one triple is of type $(1,2,3)$ (so the
Mom-$3$ is torus-friendly).

Hence we may assume that \emph{at most} three of $\{A_1,B_1\}$,
$\{B_1,A_3\}$, $\{A_3,B_2\}$, $\{B_2,A_1\}$ are in
$\mc{O}(1)\cup\mc{O}(2)$.

If all four of $\{A_1,B_1\}$, $\{B_1,A_2\}$, $\{A_2,B_2\}$,
$\{B_2,A_1\}$ lie in $\mc{O}(1)\cup\mc{O}(2)$ (or by symmetry all four
of $\{A_2,B_1\}$, $\{B_1,A_3\}$, $\{A_3,B_2\}$, $\{B_2,A_2\}$) then by
a similar argument to the one used in the beginning of the previous
case we can also construct a simpler Mom-$2$ or a simpler
geometric Mom-$3$ structure which is torus-friendly. So assume this does not
happen either. In summary, we can assume that for none of the three
sets of horoball pairs listed above do all four pairs in
the set lie in $\mc{O}(1)\cup\mc{O}(2)$. (In particular this implies
that $d_i\ge\delta_2$ for all $i\in\{1,2,3\}$, by Lemma
\ref{lem:one_cell_dist}.)

It is then straightforward to check that at least one of the three
statements in the above list must in fact be \emph{false}: you cannot
choose three pairs from each set without choosing all four pairs from
at least one set. In other words, exactly two of the three statements
in the above list are true.

Now consider the three pairs $\{A_1,B_1\}$, $\{A_2,B_1\}$, and
$\{A_3,B_1\}$. Suppose all three are in $\mc{O}(1)\cup\mc{O}(2)$. Then
$\{A_1,B_1,A_3\}$ is of type $(1,1,2)$ or $(1,2,2)$, while
$\{A_1,B_1,A_2\}$ and $\{A_2,B_1,A_3\}$ are each of type $(1,1,3)$,
$(2,2,3)$, or $(1,2,3)$. Note that if these last two triples are
equivalent due to $g\in\pi_1(N)$ then $g$ must be parabolic fixing
$B_1$, in which case $\lambda_3$ must intersect the arc from $B_1$ to
$A_2$, a previous case. So assume these triples are not equivalent
under $\pi_1(N)$. Then the three triples $\{A_1,B_1,A_2\}$,
$\{A_2,B_1,A_3\}$, and $\{A_3,B_1,A_1\}$ form a simpler geometric
Mom-$3$ structure, which is torus-friendly unless
$\{A_1,B_1,A_2\}$ and $\{A_2,B_1,A_3\}$ are both of type
$(1,2,3)$. Note that this is only possible if $\{A_1,B_1\}$ and
$\{A_3,B_1\}$ are in $\mc{O}(1)$ while $\{A_2,B_1\}$ is in
$\mc{O}(2)$. Now note that $\{A_j,B_2\}$ must be in
$\mc{O}(1)\cup\mc{O}(2)$ for at least one $j\in\{1,2,3\}$. If $j=2$
then $\{A_2,B_1,B_2\}$ is either a third $(1,2,3)$-triple (and not
equivalent to either $\{A_1,B_1,A_2\}$ or $\{A_2,B_1,A_3\}$ since it
shares an ``edge'' with both) or else it
is of type $(1,1,3)$ or $(2,2,3)$; either way we get a geometric
Mom-$3$ structure with either one or three triples of type $(1,2,3)$, which
therefore is torus-friendly. Suppose $j=1$ or $3$; by symmetry assume
$j=1$. Then by our previous assumptions $\{A_2,B_2\}$ and
$\{A_3,B_3\}$ must not be in $\mc{O}(1)\cup\mc{O}(2)$. Now we can
apply Corollary \ref{cor:two_line_dist} to the horoballs $A_2$, $A_3$,
$B_1$, and $B_2$ to get
\begin{eqnarray*}
e_3 e_2 + e_3 & \le & e_3 e_2 \cosh d_2 \\
\Rightarrow\ d_2 & \ge & \cosh^{-1}(1+{e_2}^{-1})
\end{eqnarray*}
Note that $e_2\le e_3\le 1.5152$ then implies $d_2>0.9$, and hence by
Lemma \ref{lem:triangle_dists_2} one of $d_1$, $d_3$ must be less than
$\delta_2$, contradicting our assumptions up to this point.

So assume at least one of $\{A_1,B_1\}$, $\{A_2,B_1\}$, $\{A_3,B_1\}$
is not in $\mc{O}(1)\cup\mc{O}(2)$, and similarly for $B_2$ instead of
$B_1$. Now go back to the three statements listed above; we know
exactly one of them is false. Suppose it is the first statement which
is false and the others true (the other possibilities can be handled
similarly). The only way to reconcile this with the sentence at the
beginning of this paragraph, is if $\{A_1,B_1\}$ and $\{A_2,B_2\}$ are
not in $\mc{O}(1)\cup\mc{O}(2)$, or the same but with the roles of
$B_1$ and $B_2$ reversed. But in either case, applying Corollary
\ref{cor:two_line_dist} to $A_1$, $A_2$, $B_1$, $B_2$ yields
\begin{eqnarray*}
{e_3}^2+1 & \le & {e_3}^2 \cosh d_1 \\
\Rightarrow\ d_1 & \ge & \cosh^{-1}(1+{e_3}^{-2})
\end{eqnarray*}
And once more, if $e_3\le 1.5152$ then this implies $d_1>0.9$, so by
Lemma \ref{lem:triangle_dists_2} one of $d_2$, $d_3$ must be less than
$\delta_2$, contradicting our assumptions. This completes the proof in
this case.

\noindent\emph{Case 13:\ }$\mathcal{O}(4)$\emph{-edge, }$(1,1,4)$%
-\emph{triple.} The exact same argument as in case 2 applies,
with the obvious modifications.

This, finally, completes the proof of Proposition
\ref{prop:edge_cell}, and in turn the proof of Theorem
\ref{thrm:comb_embedded}.\qed


\section{Torus-friendly Mom-$n$'s}

Having established Theorem \ref{thrm:comb_embedded}, we now have an
embedded cellular complex $\Delta$ corresponding to a geometric
Mom-$n$ structure in the cusped manifold $N$. The next
step in upgrading $\Delta$ to an internal Mom-$n$ structure of the
type defined in \cite{gmm2} is to prove the following:

\begin{theorem}
Suppose $N$ is a one-cusped hyperbolic 3-manifold with $\Vol(N)$ $\le
2.848$ and suppose $\Delta$ is the embedded cellular complex
corresponding to the geometric Mom-$n$ structure produced by
Theorems \ref{thrm:comb_existence} and \ref{thrm:comb_embedded}. Then
the components of $N-\Delta$ each have torus boundary, or else there
exists a simpler geometric Mom-$n$ structure which is also of the
type described in Theorems \ref{thrm:comb_existence} and
\ref{thrm:comb_embedded}.
\label{thrm:comb_notfalse}
\end{theorem}

\noindent\emph{Proof:} From Theorem \ref{thrm:comb_existence} we know
that $n=2$ or $3$ and that the geometric Mom-$n$ structure thus
constructed is torus-friendly. Recall that this means either $n=2$, or
$n=3$ and the Mom-$3$ structure does \emph{not} have exactly two
triples of type $(p,q,r)$ where $p$, $q$, and $r$ are distinct indices.
We wish to prove that being torus-friendly implies that the
boundary consists of a collection of tori.

Let $M$ be a thickened
neighborhood of $\Delta$; note that by construction $\chi(\partial
M)=2\chi(\Delta)=0$. If $\partial M$ consists of nothing but tori then
we're done. If $\partial M$ contains components which are not tori,
then one of those components must be a sphere. Hence we wish to show
that $\partial M$ does not contain any spherical components. Since $N$
is hyperbolic this is equivalent to showing that $N-\Delta$ does not
have any components which are 3-balls.

Suppose $n=2$ and that one of the components of \(N-\Delta\) is a
3-ball $B$. In the universal cover \(\BH^3\), $B$ lifts to a
collection of 3-balls; choose one of them and call it
\(\tilde{B}\). The boundary of \(\tilde{B}\) consists of two types of
``faces''. First, there are totally geodesic faces which are the lifts
of $2$-cells of $\Delta$ corresponding to triples.  Note that no such
$2$-cell can possibly contribute more
than two faces to $\partial\tilde{B}$.  Second, $\tilde{B}$ has
horospherical faces which are lifts of pieces of the cusp torus. There
can be an arbitrary number of such faces, but each such horospherical
face must only be adjacent to totally geodesic faces.  (In particular
the number of totally geodesic faces of \(\tilde{B}\) must be greater
than zero.)  We can extrude \(\tilde{B}\) in the direction of the
horospherical faces (or equivalently, extrude $B$ in the direction of
the cusp of $N$ before lifting to \(\BH^3\)) to obtain an ideal
hyperbolic polyhedron which contains \(\tilde{B}\), and whose ideal
triangular faces each contain a unique totally geodesic face of
\(\tilde{B}\).

But a geometric Mom-$2$ structure only has two triples, and hence
\(\tilde{B}\) cannot have more than four totally geodesic faces. Thus,
the ideal polyhedron containing \(\tilde{B}\) must be an ideal
tetrahedron. But the same argument can be made for any component of
\(N-\Delta\) which is not the cusp neighborhood, and we've already used up all
of the available totally geodesic faces. Therefore $B$ and the cusp
neighborhood must be the only components of \(N-\Delta\). This is
impossible, since \(\chi(\partial M)=0\). Therefore $N-\Delta$ cannot
contain any $3$-ball components.

Now suppose $n=3$ and that one of
the components of \(N-\Delta\) is a 3-ball $B$. As before, $B$ lifts
to a 3-ball \(\tilde{B}\) in \(\BH^3\), which is in turn contained in
an ideal hyperbolic polyhedron whose ideal triangular faces each
contain a unique totally geodesic face of \(\tilde{B}\). This time,
there are six totally geodesic faces available; since a polyhedron
with triangular faces must have an even number of faces, the number of
totally geodesic faces of \(\tilde{B}\) must be either $4$ or $6$. If
the number is $6$, then as before this implies that $B$ and the cusp
are the only components of \(N-\Delta\) which is impossible. Therefore
\(\tilde{B}\) has $4$ totally geodesic faces, and hence is contained
in an ideal hyperbolic tetrahedron. This implies that $B$ is a
truncated ideal hyperbolic tetrahedron, where the faces arising
from the truncation are horospherical instead of geodesic.

Now consider the $2$-cells in $\Delta$ which correspond to triples
from the geometric Mom-$3$ structure. Since there are only three of them,
some pair of faces of $B$ must arise from two sides of the same
$2$-cell. Suppose the corresponding triple is of type $(a,b,c)$ where
$a$, $b$, and $c$ are distinct integers. The two faces of $B$ must
share a common edge. Lifting up to $\BH^3$, we see that $\tilde{B}$
has two geodesic faces which project down to the same $2$-cell in
$\Delta$, and that these two faces have a common edge. There must be a
group element $g\in\pi_1(N)$ which sends one face to the other and
fixes the common edge, either fixing the horoballs at each end of the
edge or swapping them. Therefore $g$ is the identity or is elliptic of
order $2$; either result is a contradiction.

Therefore the  triple in question is of type $(a,a,b)$ for some
distinct $a$ and $b$. (Lemma \ref{lem:no_bad_triples} excludes the
possibility that it is of type $(a,a,a)$.) Again, the two
corresponding faces must share a common edge. This edge cannot be the
$\mc{O}(b)$-edge (as defined in Section 6) by the same argument as in
the previous paragraph; hence the common edge must be a
$\mc{O}(a)$-edge. Choose an orientation for this edge; i.e., make the
corresponding $1$-cell a directed arc. Lifting to $\BH^3$ again, this
induces an orientation on at least three of the edges of $\tilde{B}$,
namely all of those edges which project down to the
$\mc{O}(a)$-edge. There are now two possibilities; see
\begin{figure}[tb]
\begin{center}
\includegraphics{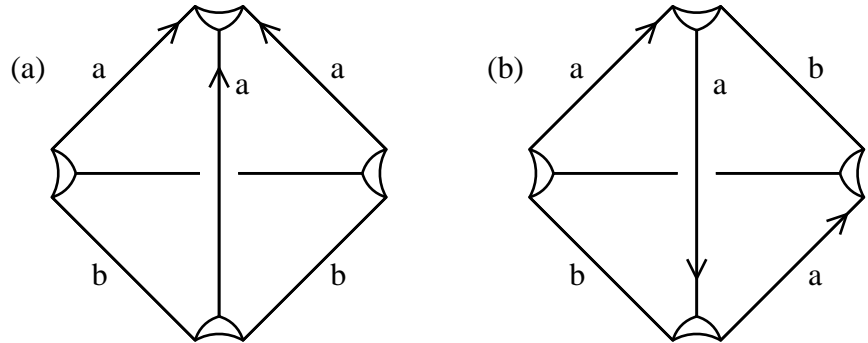}
\end{center}
\caption{Configuration of the three-cell \(\tilde{B}\), both (a) invalid
and (b) valid}
\end{figure}
figure 8. First, all three of
those edges may be oriented toward the same horospherical face of
\(\tilde{B}\), as in figure 7(a). Note that there is a group element
\(g\in\pi_1(N)\) which sends one face of type $(a,a,b)$ to the other, and that
this group element preserves the orientation of the $\mc{O}(a)$-edges;
therefore $g$ must be a parabolic
element fixing the horosphere toward which those edges point.
But a parabolic non-trivial element of \(\pi_1(N)\) which takes a
horosphere to itself must act on that horosphere by a translation.
This would imply that all three edges labelled $a$ in figure 7(a) lie
in the same hyperbolic plane, and that \(\tilde{B}\) therefore is
flat, which contradicts our assumption that $\Delta$ was embedded in
$N$.

The other possibility is that some two of the edges are oriented toward
different horospherical faces. Some thought will show that such a configuration
must look like figure 7(b); we now assume that this is the situation that
pertains.

In this situation, we have labelled all but one of the six geodesic
edges of $\tilde{B}$; now we turn our attention to the last
edge. Suppose that this edge projects down to the $\mc{O}(a)$-edge or
the $\mc{O}(b)$-edge. Then at least two of the three triples in our
Mom-$3$ only incorporate the orthopair classes $\mc{O}(a)$ and
$\mc{O}(b)$; throwing away the third triple will leave us with a
simpler geometric Mom-$2$ structure as desired. Hence, suppose that the last
edge projects down to the $\mc{O}(c)$-edge, where $c$ is the remaining
index used in the Mom-$3$. Thus the other two geodesic faces of
$\tilde{B}$ both project to $2$-cells corresponding to triples of type
$(a,b,c)$. Note that these last two faces cannot project down to the
same $2$-cell, since that would imply the existence of a non-trivial
$g\in\pi_1(N)$ which is elliptic or the identity, just as
before. Therefore the Mom-$3$ contains exactly two triples of type
$(a,b,c)$ where $a$, $b$, and $c$ are distinct, contradicting the
assumption that the Mom-$3$ was torus-friendly. This proves the
theorem.\qed

\bigskip
It is worth pointing out that embedded geometric Mom-$3$'s which are
not torus-friendly, where some component of their complement is a $3$-ball,
do exist; we give an example below. It is also often the case that a
manifold can have both a Mom-$3$ which is torus-friendly and a Mom-$3$
which is not. There are several hyperbolic manifolds which give rise to an
embedded ``Mom-like'' cellular complex with three $1$-cells and
\emph{four} $2$-cells, such that two or three of the $2$-cells are of
type $(a,b,c)$. If one discards one of the four $2$-cells, then the
resulting geometric Mom-$3$ structure may or may not be torus-friendly
depending on which $2$-cell is discarded. It is easy to see,
however, that if triples of type $(a,b,c)$, $(a,b,c)$, and $(a,a,b)$
have already been found, then any possible fourth triple can be
combined with some two of the first three triples to produce a Mom-$3$
which is torus-friendly.

As an example of a manifold containing a geometric Mom-$3$ which is not
torus-friendly, consider the manifold known as m170 in the SnapPea
census; a cusp diagram for this manifold is shown in
\begin{figure}[tb]
\begin{center}
\includegraphics{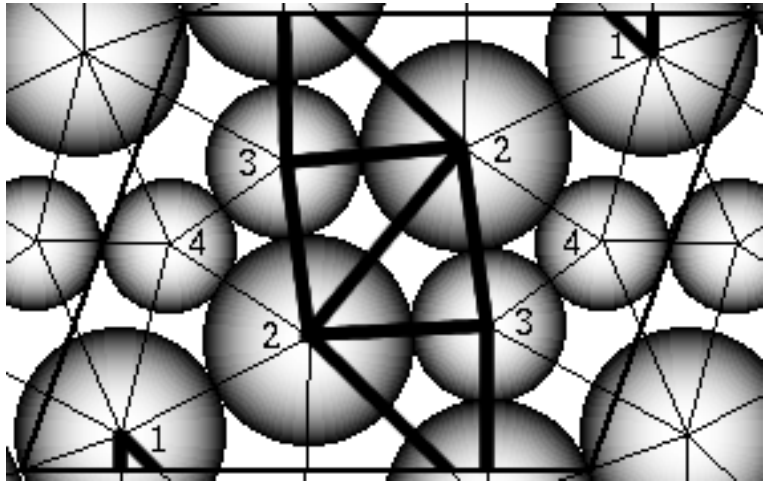}
\end{center}
\caption{The manifold m170 with a Mom-3 structure which is not
torus-friendly. The indices denote the orthoclasses of the horoballs.}
\end{figure}
figure 9. The highlighted triangles in the diagram are all corners of
the same ideal hyperbolic simplex. This simplex is bounded by only
three faces and three edges in the triangulation of m170; those faces
and edges together make up a cellular complex corresponding to a
geometric Mom-$3$ structure, with one triple of type $(2,2,3)$ and
two distinct triples of type $(1,2,3)$. This is therefore a
geometric Mom-$3$ which is not torus-friendly. At the same time the
interior of the highlighted simplex is a component of the complement
of this cellular complex, i.e. there is a component of $N-\Delta$
which is a $3$-ball. Hence this geometric Mom-$3$ cannot be turned
into a topological Mom-$3$ of the type described in \cite{gmm2}.
It is worthwhile to note, however, that this manifold does possess
other geometric Mom-$3$ structures which are torus-friendly.
Specifically there exists a fourth triple, of type
$(1,1,2)$, which can be used to construct such a Mom-$3$.

%

\section{Fullness}

Based on the result of Theorem \ref{thrm:comb_notfalse}, throughout
this section we assume that $N$ possesses a geometric Mom-$n$
structure where $n=2$ or $3$, that the corresponding cellular complex
$\Delta$ is embedded, and that the components of $N-\Delta$ which are
not cusp neighborhoods have torus boundary.  At this point we switch
from discussing cellular complexes to discussing handle structures as
follows: thicken $T$ to $T\times I$ where $T\times 0$ faces the cusp,
thicken each 1-cell of $\Delta$ to a 1-handle from $T\times 1$ to
itself, and thicken each 2-cell of $\Delta$ to a 2-handle which runs
over $T\times 1$ and three 1-handles counting multiplicity. The reason
for this change in focus is solely to take advantage of the language
and conclusions of \cite{gmm2}. Specifically, $T\times I$ and
the newly constructed 1-handles and 2-handles form a handle
decomposition of a submanifold $M\subset N$. We will abuse notation
and allow $\Delta$ to also refer to the handle decomposition of $M$;
it should always be clear in context whether we are referring to a
cellular complex or a handle decomposition. By Theorem
\ref{thrm:comb_notfalse} we may assume that the boundary of $M$ is a
collection of tori. Also, since $M$ is a subset of a hyperbolic
manifold and contains both a cusp torus and a geodesic arc from the
cusp torus to itself, $i_{\ast}\pi_1(M)$ cannot be abelian, where
$i:M\rightarrow N$ is the inclusion map. In other words,
$i:M\rightarrow N$ is a \emph{non-elementary} embedding.  Therefore
$(M,T,\Delta)$ is a topological internal Mom-$n$ structure according
to \cite{gmm2}.  As in \cite{gmm2}, we will adopt the terminology of
Matveev and refer to the intersection of the 1-handles
(resp. 2-handles) with $T\times 1$ as \emph{islands}
(resp. \emph{bridges}), and the complement in $T\times 1$ of the
islands and bridges will be called \emph{lakes}. The \emph{valence} of
an island will be defined to be equal to the valence of the
corresponding 1-handle, or equivalently the number of ends of bridges
lying on the island. We assume that these valences are at least two;
of any 1-handle has valence one simply remove both it and the 2-handle
adjacent to it to obtain a simpler Mom-$n$ structure.

Clearly each 1-handle in $\Delta$ contributes two islands while each
2-handle contributes three bridges. Suppose $\sigma$ is a 2-handle
corresponding to a triple of type $(a,a,b)$ where the type is defined
as in Section 6. In other words some lift
$\tilde{\sigma}$ of $\sigma$ in the universal cover of $N$ is adjacent
to three horoballs $\{A,B,C\}$ such that the orthopairs $\{A,B\}$ and
$\{B,C\}$ are in $\mc{O}(a)$ while $\{C,A\}\in\mc{O}(b)$. Let $a_0$
and $a_1$ denote the islands which are the endpoints of the 1-handle
around the $\mc{O}(a)$-edge, and define $b_0$ and $b_1$
similarly. Then the intersection of $\tilde{\sigma}$ with $\partial B$
projects down to a bridge whose endpoints both lie on islands in the
set $\{a_0,a_1\}$.

\begin{definition}
If this bridge described above joins $a_0$ to $a_1$ then we will say
$\sigma$ is a \emph{loxodromic 2-handle}. If instead this bridge joins
$a_i$ to itself for $i=1$ or $2$ then we will say $\sigma$ is a
\emph{parabolic 2-handle}.
\label{defn:parabolic_2handle}
\end{definition}

In either case there exists $g\in\pi_1(N)$ which sends $\{A,B\}$ to
$\{B,C\}$ since those are in the same orthopair class; furthermore $g$
is uniquely defined. If $\sigma$ is a parabolic 2-handle then we must
in fact have $g(A)=C$ and $g(B)=B$ (i.e. $g$ is a parabolic group
element), and the bridge from $a_i$ to itself must follow a straight
closed path in the cusp torus which corresponds to $g$. Wolog suppose
$i=0$, i.e. $a_0$ is joined to itself by a bridge. Then it is not hard
to see that the other two bridges corresponding to $\sigma$ must join
$a_1$ to $b_0$ and $b_1$ respectively. These bridges have equal length
by Lemma \ref{lem:eucl_dist}, and the angle between the bridges at
$a_1$ must equal the angle between the two bridge ends at $a_0$, since
both angles are equal to the angle at which $\sigma$ meets itself along
the $\mc{O}(a)$-edge. Clearly this angle is a straight angle at $a_0$,
and therefore at $a_1$ as well. Thus the configuration of islands and
bridges resulting from $\sigma$ is as in
\begin{figure}[tb]
\begin{center}
\includegraphics{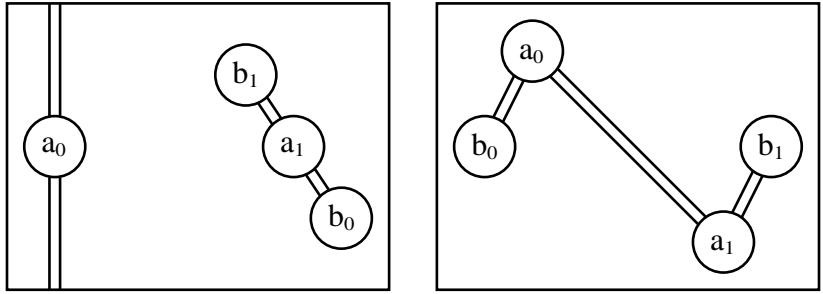}
\end{center}
\caption{The islands and bridges belonging to a parabolic (left) and
loxodromic (right) 2-handle. Note that the islands labelled $a_0$ and
$a_1$ are endpoints of a single 1-handle (and similarly for $b_0$ and
$b_1)$.}
\end{figure}
the left side of figure 10.

If $\sigma$ is instead a loxodromic 2-handle then we must have
$g(A)=B$ and $g(B)=C$ (which implies that $g$ is a loxodromic group
element although we will not prove that here). Again it is not hard to
see that wolog the other two bridges corresponding to $\sigma$ must
join $a_0$ to $b_0$ and $a_1$ to $b_1$ respectively. Again, these two
bridges must have the same length by Lemma \ref{lem:eucl_dist} and the
angle formed by the bridge ends at $a_0$ equals the angle formed by
the bridge ends at $a_1$. Given that $N$ is orientable the
configuration of islands and bridges arising from $\sigma$ must look
like the right side of figure 10.

The handle structure $\Delta$ is called \emph{full} if all of the
lakes are simply connected.  If $\Delta$ is a full handle structure,
then Theorem 4.1 of \cite{gmm2} applies to $N$. Our goal in this
section is to either prove that $\Delta$ is full or else construct a
new topological internal Mom-$n$ structure in $N$ which is full.
Lemma 4.5 of \cite{gmm2} does this in the case that
$n=2$. Specifically it shows that a topological internal Mom-2
structure that is not full can be replaced with one which is. Hence we
can combine Theorem \ref{thrm:comb_embedded}, Theorem
\ref{thrm:comb_notfalse}, Theorem 4.1 of \cite{gmm2}, and Lemma 4.5 of
\cite{gmm2} to conclude the following:

\begin{proposition}
If $N$ has a geometric Mom-2 structure, wolog that structure can
be thickened into a topological internal Mom-2 structure
$(M,T,\Delta)$ that is both full and not false. Consequently, $M$ is
hyperbolic and $N$ can be recovered by Dehn surgery on $M$.
\label{prop:full_mom2}\qed
\end{proposition}

The assertion in the second sentence is one of the primary results of
\cite{gmm2}.

Unfortunately, the topological argument used in that proof does not
extend easily to the $n=3$ case. Fortunately in this context we may
take advantage of the fact that our Mom-$3$ is more than just a
topological object; by construction the 1-handles and 2-handles of
$\Delta$ have geodesic cores, an assumption that is not made in
\cite{gmm2}. Call a Mom-$n$ structure with this additional property a
\emph{geodesic} internal Mom-$n$ structure. Our goal now is to prove
the following:

\begin{theorem}
If $(M,T,\Delta)$ is an embedded geodesic internal Mom-3 structure in a
hyperbolic manifold $N$ which is not false, with 1-handles
corresponding to the orthopair classes $\mc{O}(1)$, $\mc{O}(2)$, and
$\mc{O}(3)$, then $N$ has a full topological internal Mom-$k$
structure where $k\le 3$.
\label{thrm:full_mom3}
\end{theorem}

\noindent\emph{Proof:} Suppose $\Delta$ is not full; i.e. suppose that
$\Delta$ has one or more lakes which are not simply connected. The
possible shapes for such a lake are the following: a torus with one or
more holes, a disk with one or more holes, or an essential annulus
with zero or more holes. In each case we wish to either find a
contradiction or else construct a topological internal Mom-$2$
structure.

Suppose $T\times1$ contains a lake which is a torus with holes. Then
there is a simple closed loop $\gamma$ which bounds a disk in
$T\times1$ containing all of the islands and lakes. Push that disk
into $T\times I$ to obtain a compressing disk for $M$ which separates
$M$ into two pieces $M_{1}$ and $M_{2}$, such that $M_{1}$ is
homeomorphic to $T\times I$ minus a $0$-handle and $M_{2}$ consists of
that $0$-handle together with the $1$-handles and $2$-handles of
$\Delta$. Since $(M,T,\Delta)$ is not false all of the boundary
components of $M$ are tori; since $M$ is the connected sum of $M_{1}$
and $M_{2}$, and since $M_{1}$ has torus boundary, one of the boundary
components of $M_{2}$ must be a $2$-sphere. Since $N$ is hyperbolic,
the only possibility is that $M_{2}$ is contained inside a
$3$-ball. Therefore $M$ is the connected sum of $T\times I$ and a
sub-manifold of $N$ contained inside a $3$-ball; this is impossible if
$(M,T,\Delta)$ is geodesic.

Now suppose $T\times I$ contains a lake which is a disk with
holes. Then let $\gamma$ be a simple closed curve parallel to the
boundary of that disk, such that all of the islands and lakes inside
the disk are also inside $\gamma$. Let $T_0$ denote the component of
$\partial M$ which contains $\gamma$.  As before, $\gamma$ bounds a
disk in $T\times I$ which is a compressing disk for $M$ although it
may not separate $M$. If it does separate $M$, then arguing as in the
previous case we may show that the component of the separated manifold
which does not contain $T\times0$ is contained in a $3$-ball in $N$,
which is a contradiction if $(M,T,\Delta)$ is geodesic. So suppose
that the compressing disk does not separate $M$; let $M_{1}$ be the
manifold obtained after the compression.  Note that the connectedness
of $M_1$ implies that there must be a $1$-handle in $\Delta$ with one
endpoint inside $\gamma$ and one endpoint outside; this implies that
$T_0-\gamma$ is connected, i.e. $\gamma$ is essential in $T_0$. Hence
the compression turns $T_0$ into a connected 2-sphere boundary
component of $M_1$; call this sphere $S_1$.  $M_{1}$ also has a handle
structure $\Delta_{1}$ consisting of the $0$-handle carved out of
$T\times I$ by the compressing disk, $T\times I$ minus that $0$-handle
(which is homeomorphic to $T\times I$), and the $1$-handles and
$2$-handles of $\Delta$. Now $S_1$ must bound a $3$-ball in $N$, and
that $3$-ball must lie on the outside of $M_1$ since $M_1$ contains
$T\times 0$. Add this $3$-ball to $M_1$ as a $3$-handle to obtain a
new manifold $M_{2}\subset N$ with torus boundary and new handle
structure $\Delta_2$. Choose a $1$-handle which connects the
$0$-handle of $\Delta_2$ to $T\times I$; cancel that $1$-handle with
the $0$-handle, and cancel the $3$-handle with a $2$-handle to obtain
a handle structure $\Delta_3$ with only $1$-handles and $2$-handles.

We need to know that the embedding $i_2:M_2\rightarrow N$ is
non-elementary, but note that $M_2$ actually contains $M$: adding the
$3$-handle to $M_1$ restores the portion of $M$ that was removed by
the compression. So $i_2:M_2\rightarrow N$ is non-elementary since
$i:M\rightarrow N$ is.

Now using the methods of \cite{gmm2}, $(M_{2},T,\Delta_{3})$
can be simplified to obtain a new topological internal Mom-$k$
structure $(M_{4},T,\Delta_{4})$ on $N$. Since the construction of
$\Delta _{3}$ deleted $1$-handles and $2$-handles without adding new
ones, the complexity of $(M_{4},T,\Delta_{4})$ as defined in
\cite{gmm2} must be strictly less than the complexity of
$(M,T,\Delta)$; hence $k\leq2$. Then by Lemma 4.5\ of \cite{gmm2} we
may assume that $N$ contains a full topological internal Mom-$k$
structure with $k\leq2$.

If $T\times I$ contains a lake which is an essential annulus with one
or more holes, let $\gamma$ be a simple closed curve in the lake
parallel to the boundary of one of those holes. Then proceed just as
in the previous case.

The remaining possibility, and the one which will require the most
analysis, is that $T\times 1$ contains lakes which are essential annuli
without holes and lakes which are disks. Given six islands and nine
lakes, for Euler characteristic reasons there must be three disk lakes
and an unknown number of annulus lakes. Let $A_{1}$, $A_{2}$, \ldots ,
$A_{r}$ be the annulus lakes, and let $B_{1}$, $B_{2}$, \ldots ,
$B_{r}$ be the connected components of
$T\times1-(\cup_{i=1}^{r}A_{i})$.  Each $B_{i}$ is also an annulus,
composed of islands, bridges, and disk lakes; furthermore each $B_{i}$
must contain at least one island.

\begin{lemma}
If $r>1$ then $N$ must contain a full topological internal Mom-2
structure.
\end{lemma}

\noindent\emph{Proof:}
Suppose $r=2$. Choose simple paths $\mu_i$, $i\in\{1,2\}$, such that
$\mu_i$ crosses $B_i$ transversely for each $i$, $\mu_i$ does not
cross any island for either $i$, and such that the total number of
bridges crossed by $\mu_1$ and $\mu_2$ is minimal. Since there are
only three disk lakes total in $B_1$ and $B_2$, the number of bridges
crossed by $\mu_1$ and $\mu_2$ combined is at most five.

Slice $M$ open along $A_1\times I$ and $A_2\times I$ to obtain a new manifold
$M_1\subset N$; note that $M_1$ will still have torus boundary and the
inclusion $i_1:M_1\rightarrow N$ will still be non-elementary. $M_1$
consists of two thickened annuli, namely $B_i\times I$ for
$i\in\{1,2\}$, and the 1-handles and 2-handles of $\Delta$. Each
$\mu_i$ forms part of the boundary of a disk in $B_i\times I$; thicken
each disk to obtain a decomposition of $B_i\times I$ into a $1$-handle
(the thickened disk) and a $0$-handle (the complement of the thickened
disk). Thus we obtain a standard handle decomposition $\Delta_1$ of
$M_1$ consisting of two $0$-handles, two new $1$-handles (which we
also refer to as $\mu_1$ and $\mu_2$), and the original $1$-handles
and $2$-handles of $\Delta$. By construction the total valence of the
new $1$-handles $\mu_1$ and $\mu_2$ is at most five.

Now let $\lambda_1$, $\lambda_2$, and $\lambda_3$ denote the original
1-handles of $\Delta$, ordered in such a way that the valence of
$\lambda_1$ is maximal. Note that wolog we may assume that at least
one of the $\lambda_i$'s connects the two $B_i$'s; otherwise we can
simply throw away all of the $1$-handles and $2$-handles of $\Delta$
which are connected to $B_1$ to obtain a simpler geometric Mom-$n$ and
apply Proposition \ref{prop:full_mom2}. Suppose that $\lambda_1$ has
both endpoints on $B_1$ (or equivalently on $B_2$). Then in
$\Delta_1$, $\lambda_1$ has both endpoints on the same
$0$-handle. Drill out the core of $\lambda_1$ and this $0$-handle to
obtain a new torus boundary component $T_2$, and cancel the other
$0$-handle of $\Delta_1$ with one of the $\lambda_i$'s which connects
$B_1$ to $B_2$. The result is a new manifolds $M_2\subset N$ and an
internal topological Mom-$k$ structure $(M_2,T_2,\Delta_2)$ where
$\Delta_2$ consists of the remaining $\lambda_i$, $\mu_1$ and $\mu_2$,
and $T_2\times I$. Since the sum of the valences of the $\lambda_i$'s
equals 9, either the valence of $\lambda_1$ was at least 4 and the
valence of the cancelled $\lambda_i$ was at least 2, or else the
valence of all the $\lambda_i$'s equals 3. Either way the total valence of
the two $\lambda_i$'s removed to construct $\Delta_2$ is at least 6,
more than the total valence of the new $1$-handles $\mu_1$ and
$\mu_2$. Therefore the complexity of $(M_2,T_2,\Delta_2)$, as defined by
\cite{gmm2}, is less than the complexity of
$(M,T,\Delta)$. Consequently $k\le 2$. Then apply Lemma 4.5 of
\cite{gmm2} to complete the proof. If $\lambda_1$ connects $B_1$ to
$B_2$ but $\lambda_2$ had both endpoints on $B_1$, then drill out the
core of $\lambda_2$ and one $0$-handle and cancel $\lambda_1$ with the
other $0$-handle then proceed as above. If every $\lambda_i$ connects
$B_1$ to $B_2$ then drill out the cores of $\lambda_1$, $\lambda_2$,
and both $0$-handles at the same time to construct $T_2$, then proceed
as above. This completes the proof in this case.

If $r=3$ the proof is similar. In this case we can construct $\mu_1$,
$\mu_2$, and $\mu_3$ as paths which cross at most six bridges in
total, then split along each $A_i\times I$ and turn the $\mu_i$'s into
new $1$-handles with total valence at most 6. The $B_i$'s must all be
connected by $1$-handles; if $B_1$ were not connected to the others
then we could throw
away all of the 1-handles and 2-handles of $\Delta$ connected to $B_1$
to obtain a simpler geometric Mom-$n$ and apply Proposition
\ref{prop:full_mom2}. Hence at most one of the $\lambda_i$'s has both
endpoints on the same $B_i$. Drill out the cores of enough
$\lambda_i$'s and $0$-handles to construct a new torus boundary
component $T_2$ and cancel any remaining $1$-handles with the
remaining $0$-handles to obtain an internal topological Mom-$k$
$(M_2,T_2,\Delta_2)$. Since the $\lambda_i$'s have total valence 9 and
the $\mu_i$'s have total valence at most 6, we have reduced complexity
and hence $k\le 2$. Then again apply Lemma 4.5 of \cite{gmm2}.

If $r\ge 4$ then at least two $B_i$'s contain only one island
each. Since the bridges follow straight paths and since the $B_i$'s
are annuli, if $B_i$ contains only one island then it also contains
only one bridge. Such an island is the endpoint of a 1-handle
$\lambda_1$ which is connected to only a single $2$-handle in
$\Delta$; throw away the $1$-handle and the $2$-handle to obtain a
simpler geometric Mom-$n$, then apply Proposition
\ref{prop:full_mom2}. This completes the proof of the lemma.\qed

\bigskip
So suppose $r=1$. Let $A=A_1$ and $B=B_1$ for simplicity.
We wish to use an argument similar to the one in the above lemma to
obtain a topological Mom-$2$ structure. Specifically we wish to find a
path $\mu$ which crosses
$B$ transversely and which crosses no islands and as few bridges as
possible. Suppose $\mu$ can be chosen to cross fewer than
$v(\lambda_1)$ bridges, where $v(\lambda_1)$ is the valence of one of
the $1$-handles of $\Delta$. Then as in the lemma we can split
$T\times I$ along $A\times I$, then decompose $B\times I$ into a
$0$-handle and a $1$-handle where the $1$-handle is obtained by
thickening a disk with $\mu$ in its boundary. Then by drilling out the
cores of $\lambda_1$ and the $0$-handle we will obtain a topological
internal Mom-$2$ structure which is full by Lemma 4.5 of \cite{gmm2}.

Since $B$ contains only three disk lakes we can always choose $\mu$ to
cross at most four bridges. Thus we're done if there exists a
$1$-handle with valence five or more. However this is not always the case.
Let $\lambda_{1}$, $\lambda_{2}$, and $\lambda_{3}$ be the 1-handles of
$\Delta$. We now consider the possible values for the valences of these three
1-handles. As previously stated the valence must add up to nine, so
there are only three possibilities up to symmetry.

First suppose $v(\lambda_1)=5$ and $v(\lambda_2)=v(\lambda_3)=2$. Then
we're done by the above argument.

Next, suppose $v(\lambda_{1})=4$, $v(\lambda_{2})=3$, and
$v(\lambda_{3})=2$. If we can find a path $\mu$ connecting the
boundary components of $B$ which crosses no islands and no more than
three bridges then the usual splitting-and-drilling procedure as will
result in a strictly simpler internal Mom-$k$ structure. However it is
not immediately apparent that such a path must exist; further analysis
is required.

Let $\partial_{+}B$ and $\partial_{-}B$ denote the two
boundary components of $B$. Consider the case where $\partial_{+}B$
contains exactly one island, and hence there is a bridge connecting
that island to itself along a straight path in $T\times1$. That bridge
must be one corner of a parabolic $2$-handle, as in the left side of
figure 10. The existence of a parabolic $2$-handle turns out to have
strong geometric consequences as follows:

\begin{lemma}
Suppose $\partial_{+}B$ contains exactly one island, and that consequently
$(M,T,\Delta)$ contains a parabolic 2-handle. If
$(M,T,\Delta)$ contains a second parabolic 2-handle, then $N$ must have a full
geodesic internal Mom-2 structure.
\label{lem:two_parabolics}
\end{lemma}

\noindent\emph{Proof:}
As before, if $\lambda$ is the 1-handle corresponding to $\mc{O}(a)$
then we will denote the islands at the end of $\lambda$ by $a_{0}$ and
$a_{1}$ in some order. Also for brevity we will say that a 2-handle is
``of type $(a,b,c)$'' if it corresponds to a triple of type $(a,b,c)$.
A parabolic 2-handle by definition of type $(a,a,b)$ for some $a$ and
$b$; assume the island in $\partial_{+}B$ is the island $a_{0}$. If
the second parabolic 2-handle is of type $(a,a,b)$ or $(b,b,a)$, then
clearly we have a geometric Mom-2 structure which is necessarily
full, torus-friendly, embedded, and so forth.

Suppose the second parabolic 2-handle is of type $(a,a,c)$, where
$\mc{O}(c)$ is the remaining orthopair class in our Mom-3
structure. The only way this is possible is if the island $a_{1}$ is
also connected to itself by a bridge. But then there would be two
additional bridges which meet at the island $a_{0}$ at a straight
angle, which is impossible if that is the only island in
$\partial_{+}B$.

Suppose the second parabolic 2-handle is of type $(b,b,c)$. Then wolog
the island
$b_{0}$ is also connected to itself by a bridge; since all of the bridges are
contained in $B$, an annulus, the path this bridge follows must be in the same
homotopy class as the bridge in $\partial_{+}B$. Since the bridges are
geodesic, this means they must have the same length. Recall that if
$e_{a}$, $e_{b}$, and
$e_{c}$ are the elements of the Euclidean spectrum corresponding to
$\mc{O}(a)$, $\mc{O}(b)$, and $\mc{O}(c)$ then Lemma
\ref{lem:eucl_dist} implies that:
\[
\frac{e_{c}}{{e_b}^{2}}=\frac{e_{b}}{{e_a}^{2}}\ge 1
\]
Hence $e_{c}\geq e_{b}\geq e_{a}$, and wolog $c\ge b\ge a$. By the
assumptions of Theorem \ref{thrm:full_mom3} this implies that $a=1$,
$b=2$, and $c=3$.
Now consider
the island $a_{1}$. While there is not necessarily a
bridge connecting this island to itself, there is still a triple of horoballs
corresponding to this island and the translation $g\in\pi_{1}(T)$ which acts
in the direction of the annulus $B$. That triple must be of type $(a,a,k)$ for
some $k$. Therefore Lemma \ref{lem:eucl_dist} implies that
\[
\frac{e_{k}}{{e_a}^{2}}=\frac{e_{b}}{{e_a}^{2}}%
\]
So $e_{k}=e_{b}$. If $k\le b=2$ then this implies that $N$ has a geometric
Mom-2 structure. If $k\geq3$ then $e_{k}\geq e_{3}\geq e_{2}$, i.e.
$e_{3}=e_{2}$. But $e_{3}/{e_2}^{2}=e_c/{e_b}^2\geq1$, so
$e_{3}=e_{2}$ implies that $e_{3}=e_{2}=1$. This means that the
horoballs centered at the islands $b_{0}$, $b_{1}$, $c_{0}$, and
$c_{1}$ are all \emph{full-sized}, i.e. in the upper half-space model
they appear as Euclidean spheres of diameter 1 and are tangent to the
horoball at infinity. Since there is a parabolic 2-handle of type
$(b,b,c)$, by Lemma \ref{lem:eucl_dist} the distance from $c_0$ and
$c_1$ to $b_1$ is at most $1/e_c=1$, and hence the horoballs centered at
$c_0$ and $c_1$ abut the horoball at $b_1$. Similarly the length of
the bridge from $b_0$ to itself is $e_c/{e_b}^2=1$ so the horoball at
$b_0$ abuts itself; by symmetry, so does the horoball centered at
$b_1$. Note that this implies that shortest essential curve on the
cusp torus has length at most one, so any full-sized horoball must
abut itself. Finally the horoball at $b_1$ must abut the full-sized
horoball centered at $a_1$, thanks to the parabolic triple of type
$(a,a,b)$. This is a contradiction, as there is no way to arrange all
of these full-sized horoballs around the island $b_1$ while keeping
their interiors disjoint.

If the second parabolic 2-handle is of type $(c,c,a)$ then permute the
variables, replacing $c$ with $a$, $a$ with $b$, and $b$ with
$c$. Then proceed as in the previous case.

The remaining possibility is that the second parabolic 2-handle is of type
$(c,c,b)$. Wolog the island $c_{0}$ is connected to itself by a bridge
following a
path in the same homotopy class (and hence of the same length)\ as the bridge
in $\partial_{+}B$. Define $e_{a}$, $e_{b}$, and $e_{c}$ as before. Then%
\[
\frac{e_{b}}{e_{c}^{2}}=\frac{e_{b}}{e_{a}^{2}}\geq1
\]
Therefore $e_{b}\geq e_{c}=e_{a}$; hence $e_{2}=e_{1}$ and wolog
$a=1$, $c=2$, and $b=3$. Note that the length of the bridge connecting
$a_{0}$ to itself must be $e_{3}$ by Lemma \ref{lem:eucl_dist}.  Now
consider the island $a_{1}$; it is connected by bridges to the islands
$b_{0}$ and $b_{1}$, and these bridges both have length $1/e_{3}$ and
meet at $a_{1}$ in a straight angle. Similarly, the island $c_{1}$ is
connected to $b_{0}$ and $b_{1}$ by bridges of length $1/e_{3}$ which
meet at $c_{1}$ in a straight angle. The only way this can happen
inside the annulus $B$ is if the island $c_{1}$ lies exactly halfway
along the shortest geodesic path from $a_{1}$ to itself, and vice
versa. In other words, there are two geodesic paths of length
$e_{3}/2$ connecting $a_{1}$ to $c_{1}$. These paths correspond to
triples of horoballs of type $(1,2,k)$ and $(1,2,l)$ for some $k$ and
$l$ such that $e_{k}=e_{l}=e_{3}/2<e_{3}$; clearly $k,l\in
\{1,2\}$. These two triples cannot be equivalent under the action of
$\pi_1(N)$; since they involve the same two islands, if
$g\in\pi_{1}(N)$ mapped one triple to the other than $g$ would also
have to fix the cusp torus $T$, but the two bridges in question are
clearly not equivalent under the action of $\pi_{1}(T)$. Therefore the
triples $(1,2,k)$ and $(1,2,l)$ constitute a geometric
Mom-2 structure. This completes the proof of the lemma.\qed

\bigskip
So therefore we may assume that there is at most one island connected
to itself by a bridge; in particular $\partial_{-}B$ contains more
than one island.

Next consider a minimum-length sequence $\gamma=\{\iota_{0},\beta_{1}%
,\iota_{1},\ldots,\beta_{n},\iota_{n}\}$ of islands and bridges such that
$\iota_{0}$ is the sole island in $\partial_{+}B$ and $\iota_{n}$ is an island
in $\partial_{-}B$. Note there are at most four bridges in the sequence
(otherwise there would be more than six islands in total). If
$\gamma$ contains only one or two bridges, it is straightforward to
show that there must be a path $\mu$ from $\partial_{+}B$ to
$\partial_{-}B$ lying in a small neighborhood of $\gamma$ which
crosses at most three bridges. (Remember that at most two of the
islands in the sequence $\gamma$ can have valence 4.)  Suppose
$\gamma$ contains four bridges; then together $\gamma$,
$\partial_{+}B$, and $\partial_{-}B$ contain all six islands and at
least seven bridges. If the remaining bridges are not placed in such a
way that at least two bridge ends meet the islands in $\gamma$ on each
side of the sequence (see
\begin{figure}[tb]
\begin{center}
\includegraphics{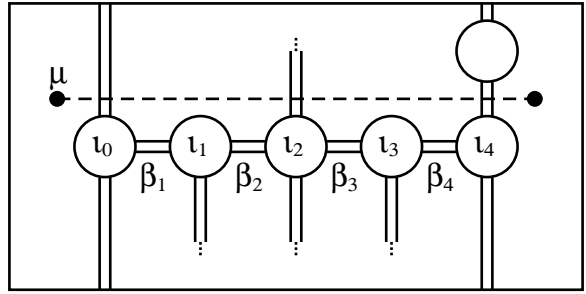}
\end{center}
\caption{Finding a path $\mu$ from $\partial_+ B$ to $\partial_- B$
which crosses at most three bridges.}
\end{figure}
figure 11), then there will be
a path $\mu$ crossing at most 3 bridges. So suppose there are at least
two additional bridge ends on each side of the sequence; this implies
that each of the two remaining bridges join an island in $\gamma$ to
another island in $\gamma$ (and not, say to the other island in
$\partial_{-}B$).
Now note that up to reordering of the indices there are only
two combinations of triples which result in a Mom-$3$ structure which is
torus-friendly and with these valences, and without containing a
Mom-$2$ structure as a subset: triples of type $(a,a,b)$,
$(a,a,c)$, and $(b,b,c)$ for some ordering $\{a,b,c\}$ of the indices
$\{1,2,3\}$, or triples of type $(a,a,c)$, $(a,b,b)$, and $(a,b,c)$.
There is no way to place the two remaining
bridges that (a) preserves the minimality of $\gamma$ (b) doesn't
imply the existence of a second parabolic 2-handle, (c) ensures two
islands of each valence $\{2,3,4\}$, and (d) ensures that the number
of bridges between each island matches the numbers produced by one of
the two combinations of triples described above. If $\gamma$ contains 3
bridges then $\gamma$, $\partial_{+}B$, and $\partial_{-}B$ contain at
least five islands and at least six bridges. Again, if there are not
two additional bridge ends on each side of $\gamma$ then we can find
an appropriate path $\mu$; so suppose there are at least two
additional bridge ends on each side. This implies that the missing
island is of valence 2, since otherwise it would account for too many
bridge ends. Again, there is no way to place an island of valence 2
and the three remaining bridges that satisfies the conditions (a),
(b),\ (c), and (d) above.

Therefore if $\partial_{+}B$ contains exactly one island then we can find a
path $\mu$ from $\partial_{+}B$ to $\partial_{-}B$ crossing at most 3 bridges.

Suppose then that $\partial_{+}B$ and $\partial_{-}B$ each contain at least
two islands. Define $\gamma$ as before; there are at most 3 bridges in the
sequence. If there are only 1 or 2 bridges, then we can find $\mu$
crossing at most three bridges as before,
so suppose $\gamma$ contains exactly 3 bridges. Then $\gamma$, $\partial_{+}%
B$, and $\partial_{-}B$ contain all six islands and at least seven bridges.
Furthermore the two missing bridges must each join an island in $\gamma$ to
another island in $\gamma$ or else there will be less than two additional
bridge ends on one side of the sequence and hence a path $\mu$ crossing at
most 3 bridges. But there is no way to place the two missing bridges that (a)
preserves the minimality of $\gamma$, (b) ensures each missing bridge joins
$\gamma$ to itself, (c) ensures two islands of each valence
$\{2,3,4\}$, and (d) ensures that the number of bridges between each
island matches the numbers produced by one of the two combinations of
triples described earlier.

This completes the proof in the case where $v(\lambda_{1})=4$, $v(\lambda
_{2})=3$, and $v(\lambda_{3})=2$.

The remaining case is where $v(\lambda_{1})=v(\lambda_{2})=v(\lambda_{3})=3$.
Define $\partial_{+}B$ and $\partial_{-}B$ as before. If $\partial_{+}B$ and
$\partial_{-}B$ each contain exactly one island, then Lemma
\ref{lem:two_parabolics} above shows $N$
has a full internal Mom-2 structure. Suppose $\partial_{+}B$ and $\partial
_{-}B$ each contain at least two islands. Consider the minimum-length sequence
$\gamma=\{\iota_{0},\beta_{1},\iota_{1},\ldots,\beta_{n},\iota_{n}\}$ as
before; $\gamma$ has no more than three bridges. If $\gamma$ has two or fewer
bridges, then since each island has valence 3 there must be a path $\mu$ from
$\partial_{+}B$ to $\partial_{-}B$ which crosses at most two bridges; split,
thicken, and drill just as before to construct a topological internal Mom-2
structure. If $\gamma$ has exactly three bridges, then $\gamma$, $\partial_{+}B$,
and $\partial_{-}B$ contain all six islands at at least seven bridges,
but there
is no way to place the missing two bridges that (a)\ preserves the minimality
of $\gamma$ and (b) ensures all six islands have valence 3, except for
configurations like the one shown in
\begin{figure}[tb]
\begin{center}
\includegraphics{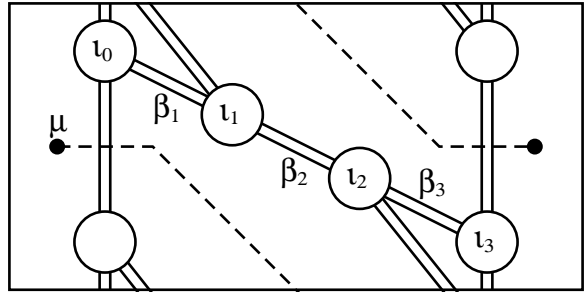}
\end{center}
\caption{If each island has valence $3$ and $\gamma$ has exactly three
bridges then we can always find a path $\mu$ crossing two bridges.}
\end{figure}
figure 12, in which there is clearly a path $\mu$ which crosses only two
bridges.

So wolog assume that $\partial_{+}B$ contains exactly 1 island and
$\partial_{-}B$ contains at least 2. If $\partial_{+}B$ contains
exactly one island then the
Mom-3 structure must contain a parabolic 2-handle of type $(a,a,b)$
for some $a$
and $b\in\{1,2,3\}$, $a\not =b$. We may assume the island in
$\partial_{+}B$ is the one denoted $a_{0}$, and hence the islands $b_{0}$
and $b_{1}$
are each joined to $a_{1}$ by bridges which meet at $a_{1}$ in a straight
angle. Consider the other two 2-handles in the Mom-3 structure; we may assume
they are either of type $(b,b,c)$ and $(c,c,a)$, or of type $(c,c,b)$
and $(a,b,c)$, where $c$
is the remaining element of $\{1,2,3\}$, as no other combination produces
1-handles of the given valences without also including a geometric Mom-2
structure. Suppose the other 2-handles are of type $(b,b,c)$ and
$(c,c,a)$. If either
of these 2-handles are parabolic then we can apply Lemma
\ref{lem:two_parabolics}; so
suppose neither is parabolic. A loxodromic 2-handle of type $(b,b,c)$
implies that there is a bridge joining $b_{0}$ to $b_{1}$; this bridge
together with the two-step path from $b_{0}$ to $a_{1}$ to $b_{1}$ must form a
homotopically non-trivial loop in $B$. Also, wolog the island $c_{i}$ is
joined to $b_{i}$ by a bridge for $i=0,1$. Moreover since $(M,T,\Delta)$ is a
geodesic structure, the angles subtended by the three bridges at $b_{0}$ must
equal the angles subtended by the three bridges at $b_{1}$, albeit in opposite
order since $N$ is orientable. Thus the bridges and islands must be arranged
as in
\begin{figure}[tb]
\begin{center}
\includegraphics{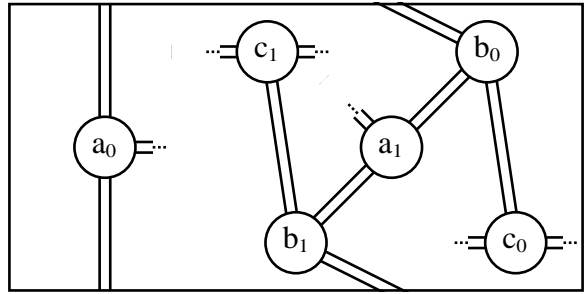}
\end{center}
\caption{The islands and bridges corresponding to a parabolic 2-handle
of type $(a,a,b)$ and a loxodromic 2-handle of type $(b,b,c)$. Note there
is no way to include a loxodromic 2-handle of type $(c,c,a)$ as well.}
\end{figure}
figure 13; note in particular that $c_{0}$ and $c_{1}$ must be
on opposite sides of the path $b_{0}\rightarrow a_{1}\rightarrow
b_{1}\rightarrow b_{0}$. Now the third 2-handle, which is of type $cca$ and
which we are assuming is non-parabolic, must imply the existence of a bridge
in $B$ which joins $c_{0}$ to $c_{1}$ without crossing any other bridge, which
is clearly impossible.

Suppose then that the other $2$-handles in the Mom-$3$ structure are
of type $(c,c,b)$ and $(a,b,c)$. Again, we may assume the $2$-handle
of type $(c,c,b)$ is a loxodromic $2$-handle, not a parabolic
one. Wolog as in Figure 10 there
exist bridges joining $b_{0}$ to $c_{0}$, $c_{0}$ to $c_{1}$, and
$c_{1}$ to $b_{1}$. These bridges together with the bridges joining
$b_{0}$ and $b_{1}$ to $a_{1}$ must again follow a homotopically
non-trivial loop in $B$. In order to form such a loop in the annulus
$B$, the angle at $b_0$ between the bridge to $a_1$ and the bridge to
$c_0$ must equal the angle at $b_1$ between the bridge to $a_1$ and
the bridge to $c_1$, and those angles must be in the same direction.
But since $(M,T,\Delta)$ is a
geodesic structure and since $M$ is orientable those angles must in
fact be equal in the opposite direction as well. Hence those two
angles are both straight angles; see
\begin{figure}[tb]
\begin{center}
\includegraphics{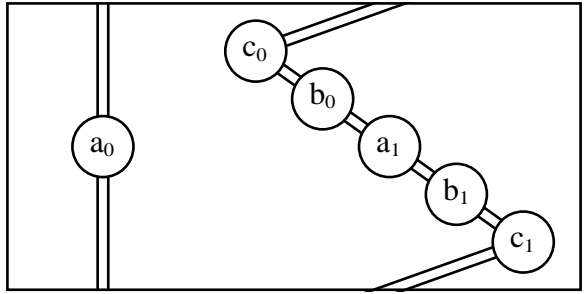}
\end{center}
\caption{The islands and bridges corresponding to a parabolic 2-handle
of type $(a,a,b)$ and a loxodromic 2-handle of type $(c,c,b)$. Note
there is no way to insert another bridge from $a_1$ to any of the
$b_i$'s or the $c_i$'s.}
\end{figure}
figure 14. The last 2-handle,
of type $(a,b,c)$, implies that there are three more bridges:\ one
joining $a_{i}$
to $b_{j}$ for some $i$ and $j$, another joining $b_{k}$ to $c_{l}$ for some
$k$ and $l$ where $k\not =j$, and a third joining $c_{s}$ to $a_{t}$ where
$s\not =l$ and $t\not =i$. In particular there must be a another
bridge joining $a_1$ to one of $b_0$, $b_1$, $c_0$, or $c_1$, but
there is no way to lay such a bridge without intersecting a
pre-existing bridge. This is a contradiction.

This completes the proof of the theorem.\qed

\section{Proof of Theorem \ref{thrm:big_thrm} and applications}

We begin by tying together the last three sections along with
\cite{gmm2} and completing the proof of Theorem \ref{thrm:big_thrm}.

Suppose $N$ is a one-cusped hyperbolic $3$-manifold with $\Vol(N)\le
2.848$. Then by Theorem \ref{thrm:comb_existence} $N$ contains a
geometric Mom-$2$ or Mom-$3$ structure (which is not false), and
by Theorem \ref{thrm:comb_embedded} we may assume that the corresponding
cellular complex $\Delta$ is embedded in $N$. By Theorem
\ref{thrm:comb_notfalse}, we may further assume that the components of
$N-\Delta$ which are not cusp neighborhoods have torus boundary.
Then by either Proposition \ref{prop:full_mom2} or Theorem
\ref{thrm:full_mom3} as appropriate, we may assume that the
submanifold $M$ and handle decomposition $\Delta$ obtained by
thickening the geometric Mom-$2$ or Mom-$3$ structure satisfy
the definition of a full topological internal Mom-$k$ structure as given
in \cite{gmm2}.

Then by Theorem 4.1 of \cite{gmm2} we can conclude that there
exists a full topological internal Mom-$k$ structure $(M,T,\Delta)$ in
$N$ where $k\le 3$ and $M$ is hyperbolic. (Note this Mom-$k$ structure
may bear little to no resemblance to the structure we started with!)
This implies that $N$ can be recovered from $M$ by a hyperbolic Dehn
filling on all but one of the cusps of $M$. 

Finally Theorem 5.1 of \cite{gmm2} enumerates the possible choices for
$M$. There are only 21 hyperbolic manifolds $M$ which can form part of
a full topological internal Mom-$k$ structure $(M,T,\Delta)$ for $k\le
3$, and these are precisely the manifolds listed in the table in
figure 1. This completes the proof of Theorem \ref{thrm:big_thrm}.\qed

\bigskip
Theorem \ref{thrm:big_thrm} does not in itself constitute an
enumeration of all one-cusped manifolds with volume less than
$2.848$. However it is possible to analyze the Dehn surgery spaces of
each of the 21 manifolds listed in figure 1. We use the following
theorem from \cite{fkp}:

\begin{theorem}(Futer, Kalfagianni, and Purcell): Let $M$ be a
complete, finite-volume hyperbolic manifold with cusps. Suppose $C_1$,
\ldots, $C_k$ are disjoint horoball neighborhoods of some subset of
the cusps. Let $s_1$, \ldots, $s_k$ be slopes on $\partial C_1$,
\ldots, $\partial C_k$, each with length greater than $2 \pi$. Denote
the minimal slope length by $l_{\mathrm{min}}$. If $M(s_1,\ldots,s_k)$
satisfies the geometrization conjecture, then it is a hyperbolic
manifold, and
\[
\Vol(M(s_1,\ldots,s_k)) \ge \left(1-\left(\frac{2\pi}
{l_{\mathrm{min}}}\right)^2\right)^{3/2} \Vol(M).
\]
\label{thrm:fkp}
\end{theorem}

Therefore if $M$ is one of the two-cusped manifolds listed in figure
1, and if we wish to enumerate all one-cusped manifolds with volume
less than or equal to $2.848$ that can be obtained by filling, then it
is only necessary to examine surgery coefficients with slope less than
or equal to
\[
2\pi
\left(\sqrt{1-\left(\frac{2.848}{\Vol(M)}\right)^{2/3}}\right)^{-1}
\]

As an example, suppose $M$ is the Whitehead link complement, known as
m129 in the SnapPea census. This manifold admits a symmetry which
exchanges the cusps, therefore it does not matter which cusp we choose
to fill in. (This is true for all of the manifolds listed in figure 1
except s785.) Using SnapPea, we see that the volume of m129 is
$3.6638\ldots$, which implies that we need only consider Dehn fillings
along slopes of length less than $15.99$. Given that a
maximal cusp torus around one cusp has a
longitude of length $2 \sqrt{2}$ and a meridian of length $\sqrt{2}$
at right angles to the longitude, we need only consider Dehn fillings
with coefficients $(a,b)$ where $a$ and $b$ are relatively prime
integers satisfying $2a^2+8b^2\le 256$, clearly a finite and
manageable number of cases. For each such filling, we can use SnapPea,
\cite{mos}, and other such tools to confirm rigorously whether or not
the resulting one-cusped manifold is hyperbolic and has volume less
than $2.848$.

For s776, the only three-cusped manifold in figure 1, we need to fill
in two cusps to obtain a one-cusped manifold. However Theorem
\ref{thrm:fkp} only provides an upper bound on one of the
corresponding slopes. Nevertheless, there are finite number of
possibilities for that one slope and therefore filling in that one
slope results in a finite number of two-cusped manifolds for which we
can repeat the above analysis. Note that s776 admits symmetries which
permute all of its cusps, so again it does not matter which cusps we
fill. It should also be pointed out that Martelli and Petronio have
already determined the complete list of fillings on s776 that result
in non-hyperbolic manifolds (\cite{mp}).

The results of the above analysis will be presented in detail in an
upcoming paper (\cite{mm}); for now we present the results of that
analysis without proof:

\begin{theorem}
The only one-cusped orientable hyperbolic 3-manifolds with volume less
than or equal to $2.848$ are the manifolds known in the SnapPea census
as m003, m004, m006, m007, m009, m010, m011, m015, m016, and m017.
\label{thrm:one_cusped}
\end{theorem}


Theorem \ref{thrm:big_thrm} can also be used to analyze closed
hyperbolic 3-manifolds. Lemma 3.1 of \cite{acs} states the following:

\begin{lemma}Suppose that $M$ is a closed
orientable hyperbolic $3$-manifold and that $C$ is a shortest geodesic
in $M$ such that $\operatorname{tuberad}(C)\ge (\log 3)/2$. Set
$N=\operatorname{drill}_C(M)$. Then $\Vol(N)<3.02 \Vol(M)$.
\label{lem:acs}
\end{lemma}

Here $\operatorname{tuberad}(C)$ refers to the \emph{tube radius},
i.e. the maximal radius of an embedded tubular neighborhood around
$C$, while $\operatorname{drill}_C(M)$ is simply the manifold $M-C$
equipped with a complete hyperbolic metric. In other words $N$ is a
one-cusped hyperbolic 3-manifold from which the closed manifold $M$
can be recovered by Dehn filling. Two remarks are in order at this
point. First, the Weeks manifold has volume less than $2.848/3.02$ and
is the smallest known closed orientable hyperbolic 3-manifold. Second,
according to \cite{gmt} if the condition that
$\operatorname{tuberad}(C)\ge (\log 3)/2$ fails then $M$ must have
volume greater than that of the Weeks manifold anyway.  Thus combining
the above lemma with Theorem \ref{thrm:one_cusped} yields the
following:

\begin{theorem} Suppose that $M$ is a closed orientable hyperbolic
$3$-manifold with volume less than that of the Weeks manifold. Then
$M$ can be obtained by a Dehn filling on one of the 10 one-cusped
manifolds listed in Theorem \ref{thrm:one_cusped}.
\end{theorem}

Clearly we can use this result to identify the minimum-volume
closed hyperbolic 3-manifold, using \cite{fkp} as before to limit the
number of Dehn fillings that need to be considered. The results of
such an analysis will also be presented in \cite{mm}.


Of future interest is the problem of strengthening the bound of
$2.848$ in Theorem \ref{thrm:big_thrm}, and thus improving the
classification of both closed and cusped low-volume hyperbolic
3-manifolds. The SnapPea census of cusped manifolds suggests that the
number $2.848$ should be far higher: the smallest known manifold which
does not possess an internal Mom-$2$ structure is m069, which has a
volume greater than $3.4$. It would not be unreasonable to attempt to
prove a stronger version of Theorem
\ref{thrm:big_thrm} which would apply to all one-cusped manifolds with
volume less than or equal to $3.7$; with such a result it would be
possible to determine the first infinite string of volumes of
one-cusped hyperbolic manifolds limiting on the volume of the
Whitehead link complement. Similarly it should be possible to
determine the first infinite string of volumes of closed manifolds
limiting on the figure-eight knot complement.


Finally we expect that Mom-technology can be applied directly to
closed manifolds, and not just indirectly via the use of Lemma
\ref{lem:acs}.  The obstacles to this are primarily geometrical
rather than conceptual. While the definition of a geometrical Mom-$n$
structure can easily be extended to closed manifolds (by considering
triples of geodesics rather than triples of horoballs), difficulties
arise when we consider the $\lessvol$ and $\overlapArea$ functions
defined in Lemmas \ref{lem:lessvol} and \ref{lem:overlapArea}. In the
cusped case the $\lessvol$ function uses the fact that the equidistant
surface between two horoballs is a plane. In the closed case the
equidistant surface between neighbourhoods of two geodesics is a more
complicated surface. Similarly the $\overlapArea$ function uses the
fact that the shadow of one horoball on the surface of another is a
Euclidean circle; in the closed case, the shadow of one tubular
neighbourhood on the surface of another is typically a
not-quite-elliptical region whose shape depends on the both the
distance and the angle between the corresponding core geodesics.
Developing a Mom-based theorey for closed manifold will require a more
sophisticated analysis of these two geometrical problems.


\end{document}